\documentclass[11pt]{article}
\usepackage[utf8]{inputenc}
\usepackage[T1]{fontenc}
\usepackage[frenchb,english]{babel} 
\usepackage{textcomp}
\usepackage{amsmath,amssymb}
\usepackage{amsthm}
                  
\usepackage{lmodern}

\usepackage[a4paper,margin=2cm]{geometry}

\usepackage{graphicx}             
\usepackage{xcolor}               
\usepackage{microtype,stmaryrd}          

\usepackage{hyperref}
\hypersetup{pdfstartview=XYZ}

\usepackage{bbm}
\usepackage{mathrsfs}
\usepackage{subfig}

\def\build#1_#2^#3{\mathrel{
\mathop{\kern 0pt#1}\limits_{#2}^{#3}}}
\def\llbracket{[\hspace{-.10em} [ }
\def\rrbracket{ ] \hspace{-.10em}]}

\newtheorem{theorem}{Theorem}
\newtheorem{proposition}[theorem]{Proposition}
\newtheorem{definition}[theorem]{Definition}
\newtheorem{lemma}[theorem]{Lemma}

\newtheorem{corollary}[theorem]{Corollary}

\def\w{\mathrm{w}}
\def\t{\mathcal{T}}
\def\f{\mathcal{F}}

\def\W{\mathcal{W}}
\def\S{\mathcal{S}}

\def\N{\mathbb{N}}

\def\D{\mathbb{D}}
\def\P{\mathbb{P}}
\def\E{\mathbb{E}}
\def\R{\mathbb{R}}
\def\z{\mathcal{Z}}

\def\ee{\mathcal{E}}
\def\ve{{\varepsilon}}
\def\la{\longrightarrow}
\def\da{\downarrow}
\def\ov{\overline}

\def\dd{\mathrm{d}}
\def\wh{\widehat}
\def\wt{\widetilde}
\def\tr{\mathrm{tr}}
\def\mm{\mathcal{M}}
\def\nn{\mathcal{N}}
\def\pp{\mathcal{P}}
\def\bp{\mathbf{p}}
\def\qq{\mathcal{Q}}
\def\H{\mathbb{H}}
\def\B{\mathbb{B}}
\def\BP{\mathbb{B}\mathbb{P}}
\def\ll{\mathcal{L}}
\def\rr{\mathcal{R}}

\def\rem{\noindent{\bf Remark. }}

\title{Spine representations for non-compact models\\
of random geometry\footnote{Supported by the ERC Advanced Grant 740943 {\sc GeoBrown}}}

\author{Jean-Fran\c cois Le Gall, Armand Riera}

\date{\small Universit\'e Paris-Saclay}
\begin{document}
\maketitle

\begin{abstract}
We provide a unified approach to the three main non-compact models of
random geometry, namely the Brownian plane, the infinite-volume Brownian disk,
and the Brownian half-plane. This approach allows us to investigate
relations between these models, and in particular to prove that
complements of hulls in the Brownian plane are infinite-volume Brownian disks.
\end{abstract}

\tableofcontents

\section{Introduction}

In the recent years, much work has been devoted to the continuous models 
of random geometry that arise as scaling limits of planar maps, which are discrete graphs embedded in the
sphere. The most famous model is the
Brownian map or Brownian sphere, which is the limit in the Gromov-Hausdorff sense of large 
planar maps with $n$ faces chosen uniformly at random in a suitable class and viewed as metric spaces for the graph distance rescaled by the factor
$n^{-1/4}$, when $n\to\infty$ (see in particular \cite{Abr,AA,BJM,Uniqueness,Mar,Mie}). The rescaling factor $n^{-1/4}$ is relevant
because the typical diameter of a random planar map with $n$ faces is of order $n^{1/4}$ when 
$n$ is large, and thus the rescaling leads to a compact limit. However, choosing
a rescaling factor that tends to $0$ at a slower rate than $n^{-1/4}$ yields a different limiting
space, which can be interpreted as an infinite-volume version of the Brownian sphere and is called the
Brownian plane \cite{Plane,CLG}. On the other hand, scaling limits of random planar maps with a boundary
have also been investigated \cite{AHS,BMR,Bet,BM,GM2,Mar2}. In that case, assuming that the planar map
has a fixed boundary size equal to $n$ and a volume (number of faces) of order $n^2$, rescaling 
the graph distance by the factor $n^{-1/2}$ again leads to a compact limiting space called the
Brownian disk. If however the volume grows faster than $n^2$, the same rescaling yields 
a non-compact limit which is the infinite-volume Brownian disk. The so-called Brownian half-plane
arises when choosing a rescaling 
factor that tends to $0$ at a slower rate than $n^{-1/2}$.
A comprehensive discussion of all possible scaling limits of large random quadrangulations
with a boundary, including the cases mentioned above, is given in the recent paper of Baur, Miermont and Ray \cite{BMR}.
In this discussion, the Brownian disk and the infinite-volume Brownian disk, the Brownian plane and the Brownian half-plane
play a central role. It is worth noting that the Brownian plane is also the scaling limit \cite{Plane,Bud} in the local Gromov-Hausdorff sense
of the random lattices called the uniform infinite planar triangulation (UIPT) and the uniform infinite
planar quadrangulation (UIPQ), which have been studied extensively since the
introduction of the UIPT by Angel and Schramm \cite{AS}. Similarly,
the Brownian half-plane arises as the scaling limit \cite{BMR,GM0} of the uniform half-plane quadrangulation, which has been 
introduced and 
studied in \cite{CC,CMie}. 
We finally mention that the preceding models of random geometry are closely related to
Liouville quantum gravity surfaces, and the Brownian disk, the Brownian plane and the Brownian
half-plane correspond respectively to the quantum disk, the quantum cone and the quantum wedge, see
\cite[Corollary 1.5]{MS}, and \cite{GM0} for the case of the Brownian half-plane. 

The main goal of the present article is to provide a unified approach to the three most important
non-compact models of random geometry, namely the Brownian plane, the infinite-volume Brownian disk and
the Brownian half-plane. As we will discuss below, it is remarkable that these three models
can all be constructed in a similar manner from the same infinite Brownian tree equipped with Brownian labels, 
subject to different conditionings --- the precise definition of these conditionings however requires some care especially
in the case of the infinite-volume Brownian disk. As an application of these constructions, we are able to get new relations
between the different models of random geometry. In particular, we prove that the complement of 
a hull in the Brownian plane, equipped with its intrinsic metric, is an infinite-volume Brownian disk (this may be
viewed as an infinite-volume counterpart of a property derived in \cite{Disks} for the Brownian sphere). The
latter property
was in fact a strong motivation for the present study, as it plays a
very important role in the forthcoming work \cite{Rie} concerning isoperimetric bounds in the Brownian plane.
We also prove that the ``horohull'' of radius $r$ in the Brownian plane, corresponding to the connected component 
containing the root of the
set of points whose ``relative distance'' to infinity is greater than $-r$, is a Brownian disk 
with height $r$ (here, a Brownian disk with height $r$ is obtained by conditioning a free pointed Brownian disk
on the event that the distinguished point is at distance exactly $r$ from the boundary).

Let us now explain our approach in more precise terms. The starting point of our construction is
an infinite ``Brownian tree'' $\mathfrak{T}_*$ that consists of a spine isometric to $[0,\infty)$ and two 
Poisson collections of subtrees grafted respectively to the left side and to the right side of the 
spine. For our purposes, it is very important to distinguish the left side and the right side because
we later need an order structure on the tree. The trees branching off the spine can be obtained
as compact $\R$-trees coded by Brownian excursions (so they are scaled versions of Aldous'
celebrated CRT). To be specific, in order to define the subtrees branching off the left side of the tree,
one may consider a Poisson point measure
$$\sum_{i\in I} \delta_{(t_i,e_i)},$$
with intensity $2\,\mathbf{1}_{[0,\infty)}(t)\,\dd t\,\mathbf{n}(\dd e)$, where $\mathbf{n}(\dd e)$ stands for the It\^o measure
of positive Brownian excursions, and then declare that, for every $i\in I$, the tree $\t_i$ coded by $e_i$
is grafted to the left side of the spine at level $t_i$. For subtrees branching off the right side, we proceed
in the same way, with an independent Poisson point measure. We equip $\mathfrak{T}_*$ with
the obvious choice of a distance (see Section \ref{coding-infinite} below). Then $\mathfrak{T}_*$ is a non-compact $\R$-tree, and, for every 
$v\in \mathfrak{T}_*$, we can define the geodesic line segment $\llbracket \rho,v\rrbracket$ between the
root $\rho$ (bottom of the spine) and $v$, and we use the notation $\rrbracket \rho,v\llbracket=\llbracket \rho,v\rrbracket
\backslash\{\rho,v\}$. The tree $\mathfrak{T}_*$ may be viewed as an ``infinite Brownian tree'' corresponding 
to process 2 in Aldous \cite{Ald}. 

We then introduce labels on $\mathfrak{T}_*$, that is, to each 
point $v$ of $\mathfrak{T}_*$ we assign a real label $\Lambda_v$. We let the labels on the spine be given
by a three-dimensional Bessel process $R=(R_t)_{t\geq 0}$ started from $0$. Then, conditionally on
$R$, the labels on the different subtrees are independent, and the labels on a given subtree $\t_i$
branching off the spine at level $t_i$ are given by Brownian motion indexed by $\t_i$ and started
from $R_{t_i}$ at the root of $\t_i$ (which is the point of the spine at level $t_i$). In other words,
labels evolve like linear Brownian motion when moving along a segment of a subtree branching
off the spine. 

We finally need a last operation, which ensures that we have only nonnegative labels. We let $\mathfrak{T}$
be the subset of $\mathfrak{T}_*$ that consists of all $v\in\mathfrak{T}_*$ such that
labels do not vanish along $\rrbracket \rho,v\llbracket$. So the spine is contained in $\mathfrak{T}$, but
some of the subtrees branching off the spine in $\mathfrak{T}_*$ are truncated at points where 
labels vanish. For each subtree $\t_i$, the theory of exit measures gives a way to define a quantity
$\z_0(\t_i)$ measuring the ``number'' of branches of $\t_i$ that are cut in the truncation procedure
(or in a more precise manner, the ``number'' of points $v$ of $\t_i$ such that $\Lambda_v=0$
but $\Lambda_w>0$ for every $w\in\rrbracket \rho,v\llbracket$), and we write $\z_0$ for the sum of
the quantities $\z_0(\t_i)$ for all subtrees $\t_i$ branching off the spine. We have in fact $\z_0=\infty$ a.s., 
but a key point of the subsequent discussion is to discuss conditionings of the
labeled tree $\mathfrak{T}$ that ensure that $\z_0<\infty$. 

We are now in a position to define the random metric  that will be used in the construction
of the non-compact models of random geometry of interest in this work. Set
$\mathfrak{T}^\circ=\{v\in \mathfrak{T}:\Lambda_v>0\}$, and for $v,w\in\mathfrak{T}^\circ$,
$$D^\circ(v,w)= \Lambda_v+ \Lambda_w -2\max\Big(\inf_{u\in[v,w]}\Lambda_u,\inf_{u\in[w,v]}\Lambda_u\Big),$$
where $[v,w]$ stands for the set of points visited when going from $v$ to $w$
clockwise around the tree (see Section \ref{coding-infinite} for a
more precise definition). We slightly modify $D^\circ(v,w)$ by setting $\Delta^\circ(v,w)=D^\circ(v,w)$
if the maximum in the last display is positive, and $\Delta^\circ(v,w)=\infty$ otherwise.
Finally, we let $(\Delta(v,w);v,w\in \mathfrak{T}^\circ)$ be the maximal symmetric 
function of $(v,w)\in \mathfrak{T}^\circ\times \mathfrak{T}^\circ$ that is bounded above by $\Delta^\circ$ and satisfies the triangle
inequality. It turns out that the function $(v,w)\mapsto \Delta(v,w)$ takes finite values and can be extended by continuity 
to a pseudo-metric on $\mathfrak{T}$, and we may thus consider the quotient space of $\mathfrak{T}$ for the equivalence relation defined by setting $v\simeq w$ if and 
only if $\Delta(v,w)=0$. The quotient space $\mathfrak{T}/\!\simeq$ equipped with the metric induced by $\Delta$ is:
\begin{itemize}
\item[1.] the Brownian plane under the special conditioning $\z_0=0$;
\item[2.] the infinite-volume Brownian disk with perimeter $z>0$ under the special conditioning $\z_0=z$;
\item[3.] the Brownian half-plane under no conditioning (then $\z_0=\infty$ a.s.).
\end{itemize}
The really new contributions of the present work are cases 2 and 3, because case 1 corresponds to the construction
of the Brownian plane in \cite{CLG} (which is different from the one in \cite{Plane}): in that case,
the conditioning on $\z_0=0$ turns the process of labels on the spine into a nine-dimensional
Bessel process $X=(X_t)_{t\geq 0}$ started from $0$, and the subtrees branching off the spine are conditioned to have
positive labels (see Section \ref{Br-plane} below). 

%The idea is to consider first the case when the spine is truncated at a given height $h>0$. Then the
%conditioning $\z_0=z$ makes sense, and furthermore the process $(R_t)_{0\leq t\leq h}$
%of labels
%on the spine under this conditioning has the same distribution as the process $X^{(r)}=(X_{\mathrm{L}_r+t}-r)_{0\leq t\leq h}$ where 
%$X$ is a nine-dimensional
%Bessel process as previously, $\mathrm{L}_r$ is the
%last passage time of $X$ at a level $r>0$, and the process $X^{(r)}$ is also conditioned on $\z'_0=z$,
%where (assuming that labels on the spine are now given by $X^{(r)}$ and that subtrees branching off the spine are conditioned not to 
%hit level $-r$) $\z'_0$ is the total exit measure at $0$. It is a remarkable fact that the preceding presentation 
%is valid independently of the choice of $r>0$. To summarize, in the case of a spine truncated at level $h$, the model conditioned on $\z_0=z$
%has an alternative presentation for which it is easier to pass to the limit $h\to\infty$ in order
%to get the relevant model for the infinite-volume Brownian disk. 

A remarkable feature of the preceding constructions is the fact that labels on $\mathfrak{T}$ have 
a nice geometric interpretation in terms of the associated random metric spaces $\mathfrak{T}/\!\simeq$. Precisely,
the label $\Lambda_v$ of a point $v$  of $\mathfrak{T}$ is equal to the distance from (the equivalence class
of) $v$ to the set of (equivalence classes of) points of zero label in $\mathfrak{T}/\!\simeq$. The latter set is either
a single point in case 1, or a boundary homeomorphic to the circle in case 2, or a
boundary homeomorphic to the line in case 3. Amongst other applications, this interpretation of labels allows us to
prove the above-mentioned result about the complement of hulls in the Brownian plane.
Write $\BP_\infty$ for the Brownian plane, and 
recall that the hull $B^\bullet(r)$ is defined by saying that its complement $\check B^\bullet(r):=\BP_\infty\backslash B^\bullet(r)$ is the
unbounded component of the complement of the closed ball of radius $r$ centered at the
distinguished point (bottom of the spine) of $\BP_\infty$. Then Theorem \ref{infBdBp} below
states that (the closure of) $\check B^\bullet(r)$ equipped with its intrinsic metric
is an infinite-volume Brownian disk whose perimeter is the boundary size $|\partial B^\bullet(r)|$
(see \cite{CLG} for the definition of this boundary size). 

Much of the technical work in the present paper is devoted to
making sense of the conditioning $\z_0=z$ in case 2, which is not a trivial matter because $\z_0=\infty$ a.s.
Our approach is motivated by the previously mentioned result concerning the distribution of $\check B^\bullet(r)$. At first, it would seem that our construction of the Brownian plane
from an infinite tree $\mathfrak{T}$ made of a spine equipped with labels $(X_t)_{t\geq 0}$ (given by
a nine-dimensional Bessel process),
and labeled subtrees conditioned to have positive labels,
would be suited perfectly to analyse 
the distribution of a hull or of its complement. In fact, it is observed in \cite{CLG} that  the set $\check B^\bullet(r)$ exactly corresponds to
a subtree $\mathfrak{T}_{(r)}$ consisting of the part of the spine of $\mathfrak{T}$ above level
$\mathrm{L}_r:=\sup\{t\geq 0:X_t=r\}$ and of the subtrees branching off the spine of $\mathfrak{T}$ above level $\mathrm{L}_r$
and truncated at points where labels hit $r$ --- furthermore the boundary size $|\partial B^\bullet(r)|$
is just the sum of the exit measures at level $r$ of all these subtrees. However, this representation of $\check B^\bullet(r)$ seems 
to depend heavily on $r$, even if labels are shifted by $-r$: in particular, the distribution 
of the process $(X_{\mathrm{L}_r+t}-r)_{t\geq 0}$ depends on $r$. Nevertheless, and perhaps surprisingly,
it turns out that, if we condition the boundary size $|\partial B^\bullet(r)|$ to be equal to a fixed $z>0$, the 
conditional distribution of the labeled tree $\mathfrak{T}_{(r)}$ (with labels shifted by $-r$) does not depend on $r$,
and this leads to the probability measure $\Theta_z$ which is used in our construction 
of the infinite-volume Brownian disk. The precise construction of the measures $\Theta_z$, which involves
an appropriate truncation procedure, is given in Section \ref{cod-disk}, where we also explain in which sense these measures
correspond to the conditioning of case 2 above.

As the reader will have guessed from the preceding discussion, some of the technicalities in our proofs are made necessary by the problem of conditioning on events of zero probability. For instance, in order to define the free pointed Brownian disk with perimeter $z$ and a given height $r>0$, it is relevant to condition the Brownian snake excursion measure $\N_r$ (see Section \ref{sec:preli} for a definition)
on the event that the exit measure at $0$ is equal to $z$. It is not immediately obvious how to make a canonical choice of these conditional distributions, so that
they depend continuously on the pair $(r,z)$. We deal carefully with these questions in Section \ref{dis-rela}. 

Our proofs also rely on certain explicit distributions, which are 
of independent interest. In particular, we prove that, in a free pointed Brownian disk of perimeter $1$, the density of the distribution
of the distance from the distinguished point to the boundary is given by the function
$$p_1(r):=9\,r^{-6}\Big(r+\frac{2}{3}\,r^3- 
\Big(\frac{3}{2}\Big)^{1/2}\,\sqrt{\pi}\,(1+r^2)\,\exp\Big(\frac{3}{2r^2}\Big)\,\mathrm{erfc}\Big(\sqrt{\frac{3}{2r^2}}\Big)\Big),$$
with the usual notation $\mathrm{erfc}(\cdot)$ for the complementary error function. See Propositions \ref{densityuniform} and \ref{pointed-disk} below for a short proof, which uses the representation 
of Brownian disks found in \cite{Disks} (with some more work, the same formula could also be derived from the 
representation in \cite{Bet,BM}). 

The paper is organized as follows. Section \ref{sec:preli} contains a number of preliminaries, and in particular
we introduce the formalism of snake trajectories \cite{ALG}, and the associated
Brownian snake excursion measures, to code compact
continuous random trees equipped with real labels. We also introduce the notion of a ``coding triple'' for a non-compact continuous random tree.
Such a coding triple consists of a random process representing the labels on the spine,
and two random point measures on the space of all pairs $(t,\omega)$, where $t\geq 0$
and $\omega$ is a snake trajectory (the idea is that, for every such pair, the labeled tree coded by $\omega$ will be grafted
 to the left or to the right of the spine at level $t$). 
The main goal of Section \ref{dis-rela}
is to define the coding triple associated with the infinite-volume Brownian disk or, in other words, to
make sense of the conditioning appearing in case 2 above. Section \ref{coding-metric}
then gives the construction of the random metric spaces of interest from the corresponding 
coding triples, starting from the construction of the Brownian plane in
\cite{CLG}. As an important ingredient of our discussion, we consider the free pointed Brownian disk $\D_z^{(a)}$ with perimeter 
$z$ and height $a$ (recall that the height refers to the distance from the distinguished
point to the boundary). The infinite-volume Brownian disk with perimeter $z$ is then obtained as the
local limit of $\D_z^{(a)}$ in the Gromov-Hausdorff sense when $a\to\infty$.
In an analogous manner, we construct the Brownian half-plane and we verify that it is the tangent cone
in distribution of the free Brownian disk at a point chosen uniformly on its boundary.
Finally, Section \ref{appli} is devoted to our applications to the complement of hulls and
to horohulls in the Brownian plane, and Section \ref{sec:consist} shows that our definitions of the
infinite-volume Brownian disk and of the Brownian half-plane are consistent with previous work.

\section{Preliminaries}

\label{sec:preli}

\subsection{Snake trajectories}
\label{sna-tra}

Continuous random trees whose vertices are assigned real labels play a fundamental role in this work.
The formalism of snake trajectories, which has been introduced in \cite{ALG}, provides a convenient
framework to deal with such labeled trees.

A (one-dimensional) finite path $\w$ is a continuous mapping $\w:[0,\zeta]\la\R$, where the
number $\zeta=\zeta_{(\w)}\geq 0$ is called the lifetime of $\w$. We let 
$\W$ denote the space of all finite paths, which is a Polish space when equipped with the
distance
$$d_\W(\w,\w')=|\zeta_{(\w)}-\zeta_{(\w')}|+\sup_{t\geq 0}|\w(t\wedge
\zeta_{(\w)})-\w'(t\wedge\zeta_{(\w')})|.$$
The endpoint or tip of the path $\w$ is denoted by $\wh \w=\w(\zeta_{(\w)})$.
For $x\in\R$, we
set $\W_x=\{\w\in\W:\w(0)=x\}$. The trivial element of $\W_x$ 
with zero lifetime is identified with the point $x$ of $\R$. We also use the notation $\W^\infty$, resp. $\W^\infty_x$,
for the space of all continuous functions $\w:[0,\infty)\la\R$, resp. the set of all such functions
with $\w(0)=x$. 
%Occasionally we will use the notation $\underline\w=\min\{\w(t):0\leq t\leq \zeta_{(\w)}\}$.

\begin{definition}
\label{def:snakepaths}
Let $x\in\R$. 
A snake trajectory with initial point $x$ is a continuous mapping $s\mapsto \omega_s$
from $\R_+$ into $\W_x$ 
which satisfies the following two properties:
\begin{enumerate}
\item[\rm(i)] We have $\omega_0=x$ and the number $\sigma(\omega):=\sup\{s\geq 0: \omega_s\not =x\}$,
called the duration of the snake trajectory $\omega$,
is finite (by convention $\sigma(\omega)=0$ if $\omega_s=x$ for every $s\geq 0$). 
\item[\rm(ii)] {\rm (Snake property)} For every $0\leq s\leq s'$, we have
$\omega_s(t)=\omega_{s'}(t)$ for every $t\in[0,\displaystyle{\min_{s\leq r\leq s'}} \zeta_{(\omega_r)}]$.
\end{enumerate} 
\end{definition}

We will write $\S_x$ for the set of all snake trajectories with initial point $x$
and $\S=\bigcup_{x\in\R}\S_x$ for the set of all snake trajectories. If $\omega\in \S$, we often write $W_s(\omega)=\omega_s$ and $\zeta_s(\omega)=\zeta_{(\omega_s)}$
for every $s\geq 0$. The set $\S$ is a Polish space for the distance
$$d_{\S}(\omega,\omega')= |\sigma(\omega)-\sigma(\omega')|+ \sup_{s\geq 0} \,d_\W(W_s(\omega),W_{s}(\omega')).$$
A snake trajectory $\omega$ is completely determined 
by the knowledge of the lifetime function $s\mapsto \zeta_s(\omega)$ and of the tip function $s\mapsto \wh W_s(\omega)$: See \cite[Proposition 8]{ALG}.
If $\omega$ is a snake trajectory, its time reversal $\check\omega$ defined by $\check\omega_s=\omega_{(\sigma(\omega)-s)^+}$ is also a snake trajectory.

Let $\omega\in \S$ be a snake trajectory and $\sigma=\sigma(\omega)$. The lifetime function $s\mapsto \zeta_s(\omega)$ codes a
compact $\R$-tree, which will be denoted 
by $\t_{(\omega)}$ and called the {\it genealogical tree} of the snake trajectory. This $\R$-tree is the quotient space $\t_{(\omega)}:=[0,\sigma]/\!\sim$ 
of the interval $[0,\sigma]$
for the equivalence relation
$$s\sim s'\ \hbox{if and only if }\ \zeta_s(\omega)=\zeta_{s'}(\omega)= \min_{s\wedge s'\leq r\leq s\vee s'} \zeta_r(\omega),$$
and $\t_{(\omega)}$ is equipped with the distance induced by
$$d_{(\omega)}(s,s')= \zeta_s(\omega)+\zeta_{s'}(\omega)-2 \min_{s\wedge s'\leq r\leq s\vee s'} \zeta_r(\omega).$$
(notice that $d_{(\omega)}(s,s')=0$ if and only if $s\sim s'$, and see e.g.~\cite[Section 3]{LGM} for more information about the
coding of $\R$-trees by continuous functions).  We write $p_{(\omega)}:[0,\sigma]\la \t_{(\omega)}$
for the canonical projection. By convention, $\t_{(\omega)}$ is rooted at the point
$\rho_{(\omega)}:=p_{(\omega)}(0)$, and the volume measure on $\t_{(\omega)}$ is defined as the pushforward of
Lebesgue measure on $[0,\sigma]$ under $p_{(\omega)}$. If $u,v\in\t_{(\omega)}$, $\llbracket u,v\rrbracket$ denotes 
the geodesic segment between $u$ and $v$ in $\t_{(\omega)}$, and we also use the notation $\llbracket u,v\llbracket$
or $\rrbracket u,v\llbracket$ with an obvious meaning. 

It will be useful to define also intervals on the tree $\t_{(\omega)}$. For $s,s'\in[0,\sigma]$, 
we use the convention $[s,s']=[s,\sigma]\cup [0,s']$ if $s>s'$
(and of course, $[s,s']$ is the usual interval if $s\leq s'$). If $u,v\in \t_{(\omega)}$ are distinct, then we can find
$s,s'\in[0,\sigma]$ in a unique way  so that $p_{(\omega)}(s)=u$ and $p_{(\omega)}(s')=v$ and 
the interval $[s,s']$ is as small as possible, and we define $[u,v]:=p_{(\omega)}([s,s'])$. 
Informally, $[u,v]$ is the set of all points that are visited when going from $u$ to $v$
in ``clockwise order'' around the tree. We take $[u,u]=\{u\}$.

By property (ii) in the definition of  a snake trajectory, the condition $p_{(\omega)}(s)=p_{(\omega)}(s')$ implies that 
$W_s(\omega)=W_{s'}(\omega)$. So the mapping $s\mapsto W_s(\omega)$ can be viewed as defined on the quotient space $\t_{(\omega)}$.
For $u\in\t_{(\omega)}$, we set $\ell_u(\omega):=\wh W_s(\omega)$ whenever  $s\in[0,\sigma]$ is such that $u=p_{(\omega)}(s)$  (by the previous observation, this does not
depend on the choice of $s$). We can interpret $\ell_u(\omega)$ as a ``label'' assigned to the ``vertex'' $u$ of $\t_{(\omega)}$. 
Notice that the mapping $u\mapsto \ell_u(\omega)$ is continuous on $\t_{(\omega)}$, and that, for every $s\geq 0$, the path
$W_s(\omega)$ records the labels $\ell_u(\omega)$ along the ``ancestral line'' $\llbracket \rho_{(\omega)},p_{(\omega)}(s)\rrbracket$. 
We will use the notation $W_*(\omega):=\min\{\ell_u(\omega):u\in\t_{(\omega)}\}$.

We now introduce two important operations on snake trajectories in $\S$. The first one 
is the re-rooting operation (see \cite[Section 2.2]{ALG}). Let $\omega\in \S$ and
$r\in[0,\sigma(\omega)]$. Then $\omega^{[r]}$ is the snake trajectory in $\S_{\wh W_r(\omega)}$ such that
$\sigma(\omega^{[r]})=\sigma(\omega)$ and for every $s\in [0,\sigma(\omega)]$,
\begin{align*}
\zeta_s(\omega^{[r]})&= d_{(\omega)}(r,r\oplus s),\\
\wh W_s(\omega^{[r]})&= \wh W_{r\oplus s}(\omega),
\end{align*}
where we use the notation $r\oplus s=r+s$ if $r+s\leq \sigma(\omega)$, and $r\oplus s=r+s-\sigma(\omega)$ otherwise. 
By a remark following the definition of snake trajectories, these prescriptions completely determine $\omega^{[r]}$.
The genealogical tree $\t_{(\omega^{[r]})}$ may be interpreted as the tree $\t_{(\omega)}$ re-rooted at the vertex $p_{(\omega)}(r)$ (see \cite[Lemma 2.2]{DLG} for a precise statement) and vertices of the re-rooted tree
receive the same labels as in $\t_{(\omega)}$. We sometimes write $W^{[t]}(\omega)$ instead of $\omega^{[t]}$. 

The second operation is the truncation of snake trajectories. Let $x,y\in \R$ with $y<x$. For every $\w\in\W_x$, set
$$\tau_y(\w):=\inf\{t\in[0,\zeta_{(\w)}]: \w(t)=y\}$$
with the usual convention $\inf\varnothing =\infty$ (this convention will be in force throughout this work
unless otherwise indicated). Then, if 
$\omega\in \S_x$, we set, for every $s\geq 0$,
$$\eta_s(\omega)=\inf\Big\{t\geq 0:\int_0^t \mathrm{d}u\,\mathbf{1}_{\{\zeta_{(\omega_u)}\leq\tau_y(\omega_u)\}}>s\Big\}.$$
Note that the condition $\zeta_{(\omega_u)}\leq\tau_y(\omega_u)$ holds if and only if $\tau_y(\omega_u)=\infty$ or $\tau_y(\omega_u)=\zeta_{(\omega_u)}$.
Then, setting $\omega'_s=\omega_{\eta_s(\omega)}$ for every $s\geq 0$ defines an element $\omega'$ of $\S_x$,
which will be denoted by  $\tr_y(\omega)$ and called the truncation of $\omega$ at $y$
(see \cite[Proposition 10]{ALG}). The effect of the time 
change $\eta_s(\omega)$ is to ``eliminate'' those paths $\omega_s$ that hit $y$ and then survive for a positive
amount of time. We leave it as an exercise for the reader to check that the genealogical tree $\t_{(\tr_y(\omega))}$ is
canonically and isometrically identified to the closed set $\{v\in\t_{(\omega)}:\ell_u(\omega)>y\hbox{ for every }u\in\llbracket \rho_{(\omega)},v\llbracket\}$,
and this identification preserves labels. 
%The genealogical tree of 
%$\tr_y(\omega)$ is canonically and isometrically identified 
%with the closed subset of $\t_\zeta$ consisting of all $a$ such that
%$\ell_b(\omega)\not =y$ for every strict ancestor $b$ of $a$ (excluding the root when $y=0$). 
%By abuse of notation, we often write ${\rm tr}_y(W)$ instead of ${\rm tr}_y(\omega)$.

\subsection{The Brownian snake excursion 
measure on snake trajectories}
\label{sna-mea}

Let $x\in\R$. The Brownian snake excursion 
measure $\N_x$ is the $\sigma$-finite measure on $\S_x$ that satisfies the following two properties: Under $\N_x$,
\begin{enumerate}
\item[(i)] the distribution of the lifetime function $(\zeta_s)_{s\geq 0}$ is the It\^o 
measure of positive excursions of linear Brownian motion, normalized so that, for every $\ve>0$,
$$\N_x\Big(\sup_{s\geq 0} \zeta_s >\ve\Big)=\frac{1}{2\ve};$$
\item[(ii)] conditionally on $(\zeta_s)_{s\geq 0}$, the tip function $(\wh W_s)_{s\geq 0}$ is
a Gaussian process with mean $x$ and covariance function 
$$K(s,s'):= \min_{s\wedge s'\leq r\leq s\vee s'} \zeta_r.$$
\end{enumerate}
Informally, the lifetime process $(\zeta_s)_{s\geq 0}$ evolves under $\N_x$ like a Brownian excursion,
and conditionally on $(\zeta_s)_{s\geq 0}$, each path $W_s$ is a linear Brownian path started from $x$ with lifetime $\zeta_s$, which
is ``erased'' from its tip when $\zeta_s$ decreases and is ``extended'' when $\zeta_s$ increases.
The measure $\N_x$ can be interpreted as the excursion measure away from $x$ for the 
Markov process in $\W_x$ called the Brownian snake.
We refer to 
\cite{Zurich} for a detailed study of the Brownian snake. 
For every $y<x$, we have
\begin{equation}
\label{hittingpro}
\N_x(W_*\leq y)={\displaystyle \frac{3}{2(x-y)^2}}.
\end{equation}
See e.g. \cite[Section VI.1]{Zurich} for a proof. 

The measure $\N_x$ is invariant under the time-reversal operation $\omega\mapsto \check\omega$. Furthermore,
the following scaling property is often useful. For $\lambda>0$, for every 
$\omega\in \S_x$, we define $\theta_\lambda(\omega)\in \S_{x\sqrt{\lambda}}$
by $\theta_\lambda(\omega)=\omega'$, with
$$\omega'_s(t):= \sqrt{\lambda}\,\omega_{s/\lambda^2}(t/\lambda)\;,\quad
\hbox{for } s\geq 0\hbox{ and }0\leq t\leq \zeta'_s:=\lambda\zeta_{s/\lambda^2}.$$
Then it is a simple exercise to verify that $\theta_\lambda(\N_x)= \lambda\, \N_{x\sqrt{\lambda}}$. 

\medskip
\noindent{\it Exit measures.} Let $x,y\in\R$, with $y<x$. Under the measure $\N_x$, one can make sense of a quantity that ``measures the quantity'' of paths 
$W_s$ that hit level $y$. One shows \cite[Proposition 34]{Disks} that the limit
\begin{equation}
\label{formu-exit}
L^y_t:=\lim_{\ve \da 0} \frac{1}{\ve^2} \int_0^t \dd s\,\mathbf{1}_{\{\tau_y(W_s)=\infty,\, \wh W_s<y+\ve\}}
\end{equation}
exists uniformly for $t\geq 0$, $\N_x$ a.e., and defines a continuous nondecreasing function, which is 
obviously constant on $[\sigma,\infty)$. 
The process $(L^y_t)_{t\geq 0}$ is called the exit local time at level $y$, and the exit measure 
$\z_y$ is defined by $\z_y=L^y_\infty=L^y_\sigma$. Then, $\N_x$ a.e., the topological support of the measure 
$\dd L^y_t$ is exactly the set $\{s\in[0,\sigma]:\tau_y(W_s)=\zeta_s\}$, and, in particular, $\z_y>0$ if and only if one of the paths $W_s$ hits $y$. The definition of $\z_y$
is a special case of the theory of exit measures (see \cite[Chapter V]{Zurich} for this general theory). We will use the formula for
the Laplace transform of $\z_y$: For $\lambda> 0$,
\begin{equation}
\label{Lap-exit}
\N_x\Big(1-\exp (-\lambda \z_y)\Big)= \Big((x-y)\sqrt{2/3} + \lambda^{-1/2}\Big)^{-2}.
\end{equation}
See formula (6) in \cite{CLG} for a brief justification. 

It will be useful to observe that $\z_y$ can be defined in terms of the truncated snake $\tr_y(\omega)$. 
To this end, recall the time change $(\eta_s(\omega))_{s\geq 0}$ used to define 
$\tr_y(\omega)$ at the end of Section \ref{sna-tra}, and set $\wt L^y_t=L^y_{\eta_t}$
for every $t\geq 0$. Then $\wt L^y_\infty= L^y_\infty=\z_y$, whereas formula \eqref{formu-exit}
implies that
\begin{equation}
\label{formu-exit-bis}
\wt L^y_t=\lim_{\ve \da 0} \frac{1}{\ve^2} \int_0^t \dd s\,\mathbf{1}_{\{\wh W_s(\tr_y(\omega))<y+\ve\}}
\end{equation}
uniformly for $t\geq 0$, $\N_x$ a.e.

\subsection{The positive excursion measure}
\label{posi-mea}

We now introduce another $\sigma$-finite measure
on $\S_0$, which is supported on snake trajectories taking only nonnegative values. For $\delta\geq 0$, let $\S^{(\delta)}$
be the set of all $\omega\in \S$ such that 
$\sup_{s\geq 0}(\sup_{t\in[0,\zeta_s(\omega)]}|\omega_s(t)|)> \delta$. 
Also set
$$\S_0^+=\{\omega\in\S_0: \omega_s(t)\geq 0 \hbox{ for every }s\geq 0,t\in[0,\zeta_s(\omega)]\}\cap \S^{(0)}.$$
 By \cite[Theorem 23]{ALG}, there exists a $\sigma$-finite measure $\N^*$ on
$\S$, which is supported on  $\S_0^+$
and gives finite mass to the sets $\S^{(\delta)}$ for every $\delta>0$, such that
$$\N^*(G)=\lim_{\ve\da 0}\frac{1}{\ve}\,\N_\ve(G(\tr_{0}(\omega))),$$
for every bounded continuous function $G$ on $\S$ that vanishes 
on $\S\backslash\S^{(\delta)}$ for some $\delta>0$. 
Under $\N^*$, each of the paths $W_s$, $0<s<\sigma$, starts from $0$, then
stays positive during some time interval $(0,\alpha)$, and is stopped immediately
when it returns to $0$, if it does return to $0$. 

\medskip
\noindent{\it The re-rooting formula.} We can relate the measure $\N^*$ to the excursion measures $\N_x$ of the 
preceding section via a re-rooting formula which we now state \cite[Theorem 28]{ALG}. Recall the
notation $\omega^{[t]}$ for the snake trajectory $\omega$ re-rooted at $t$. For any 
nonnegative measurable function $G$ on $\S$, we have
\begin{equation}
\label{reroot-for}
\N^*\Big(\int_0^\sigma \dd t\,G(\omega^{[t]})\Big)= 2\int_0^\infty \dd x\,\N_x\Big(\z_0\,G(\tr_0(\omega))\Big).
\end{equation}

\medskip
\noindent{\it Conditioning on the exit measure at $0$.} In a way analogous to the definition of exit measures, one can make sense of the ``quantity'' of paths $W_s$ that return to $0$ under $\N^*$.
To this end, one observes that the limit
\begin{equation}
\label{approxz*0}
\z^*_0:=\lim_{\ve\da 0} \frac{1}{\ve^2} \int_0^\sigma \dd s\,\mathbf{1}_{\{\wh W_s<\ve\}}
\end{equation}
exists $\N^*$ a.e. (this indeed follows from \eqref{formu-exit-bis}, using \eqref{reroot-for} 
to relate $\N^*$ to the distribution of $\tr_0(\omega)$ under $\N_x$, $x>0$). 

According to
\cite[Proposition 33]{ALG}, there exists a unique collection $(\N^{*,z})_{z>0}$ of probability measures on $\S_0^+$
such that:
\begin{enumerate}
\item[\rm(i)] We have
$$\N^*=
 \sqrt{\frac{3}{2 \pi}} 
\int_0^\infty \mathrm{d}z\,z^{-5/2}\, \N^{*,z}.$$
\item[\rm(ii)] For every $z>0$, $\N^{*,z}$ is supported on $\{\z^*_0=z\}$.
\item[\rm(iii)] For every $z,z'>0$, $\N^{*,z'}=\theta_{z'/z}(\N^{*,z})$.
\end{enumerate}
Informally, $\N^{*,z}=\N^*(\cdot\mid \z^*_0=z)$.

It will be convenient to have a ``pointed version'' of the measures $\N^{*,z}$.
We note that $\N^{*,z}(\sigma)=z^2$ (see the remark after \cite[Proposition 15]{Disks}) and define a probability measure on $\S_0\times \R_+$ by setting
$$\ov{\N}^{*,z}(\dd\omega \dd t)= z^{-2}\,\N^{*,z}(\dd \omega)\,\mathbf{1}_{[0,\sigma(\omega)]}(t)\,\dd t.$$

\subsection{Coding finite or infinite labeled trees}
\label{coding-infinite}

Many of the random compact (resp.~non-compact) metric spaces that we discuss in the present work 
are coded by triples $(Z,\mm,\mm')$ where $Z=(Z_t)_{t\in [0,h]}$ (resp. $Z=(Z_t)_{t\in[0,\infty)}$) is a finite (resp.~infinite) random path, and $\mm$ and $\mm'$
are random point measures on $[0,h]\times \S$ (resp.~on $[0,\infty)\times\S$).  Such a triple is called  a coding triple, and we interpret it as coding a labeled tree, having a spine 
isometric to $[0,h]$ or to $[0,\infty)$, in such a way that the path $Z$ corresponds to labels along the spine, and,
for each atom $(t_i,\omega_i)$ of $\mm$ (resp. of $\mm'$), the genealogical tree of $\omega_i$
corresponds a subtree branching off the left side (resp. off the right side) of the spine at level $t_i$. 
The random metric spaces of interest are then obtained via some identification of vertices in the labeled trees, and equipped with a
metric which is determined from the labels.

In this section we explain how coding triples are used to construct labeled trees. We follow closely the presentation given in \cite{CLG} for a special case.

\medskip

\noindent{\it The infinite spine case.} We consider a (deterministic) triple $(\w,\pp,\pp')$ such that:
\begin{description}
\item[\rm(i)] $\w\in \W^\infty$;% and $\w(t)>0$ for every $t>0$;
\item[\rm(ii)] $\pp=\sum_{i\in I} \delta_{(t_i,\omega_i)}$ and $\pp'=\sum_{i\in J} \delta_{(t_i,\omega_i)}$ are point measures on $(0,\infty)\times \S$ (the indexing sets $I$ and $J$ are disjoint),
and, for every $i\in I\cup J$, $\omega_i\in \S_{\w(t_i)}$ and $\sigma(\omega_i)>0$;
\item[\rm(iii)] all numbers $t_i$, $i\in I\cup J$, are distinct;
\item[\rm(iv)] the functions 
$$u\mapsto \beta_u:= \sum_{i\in I} \mathbf{1}_{\{t_i\leq u\}}\,\sigma(\omega_i)\;,\quad u\mapsto \beta'_u:= \sum_{i\in J} \mathbf{1}_{\{t_i\leq u\}}\,\sigma(\omega_i)\;.$$
 take finite values and are monotone increasing
on $\R_+$, and tend to $\infty$ at $\infty$ (in particular, the sets $\{t_i:i\in I\}$ and $\{t_i:i\in J\}$ are dense in $(0,\infty)$);
\item[\rm(v)] for every $t>0$ and $\ve>0$,
\begin{equation}
\label{conti-spine}
\#\Big\{i\in I\cup J: t_i\leq t\hbox{ and }\sup_{0\leq s\leq \sigma(\omega_i)}|\wh W_s(\omega_i)-\w(t_i)|>\ve\Big\}<\infty.
\end{equation}
\end{description}

Such a triple $(\w,\pp,\pp')$ will be called an infinite spine coding triple.
Recall the notation $\t_{(\omega)}$ for the genealogical tree of the snake trajectory $\omega$ and $\rho_{(\omega)}$ for the root of $\t_{(\omega)}$. 
The tree $\t_\infty$ associated with the coding triple $(\w,\pp,\pp')$ is obtained from the disjoint union
$$[0,\infty) \cup\Bigg(\bigcup_{i\in I\cup J} \t_{(\omega_i)}\Bigg)$$
by identifying the point $t_i$ of $[0,\infty)$ with the root $\rho_{(\omega_i)}$ of $\t_{(\omega_i)}$, for every $i\in I\cup J$. The metric $\dd_{\t_\infty}$
on $\t_\infty$ is determined as follows. The restriction of $\dd_{\t_\infty}$ to each tree $\t_{(\omega_i)}$ is the metric $d_{(\omega_i)}$ on
$\t_{(\omega_i)}$, and the restriction of $\dd_{\t_\infty}$ to the spine $[0,\infty)$ is the usual 
Euclidean distance. If $x\in \t_{(\omega_i)}$ and $t\in[0,\infty)$, we take $\dd_{\t_\infty}(x,t)= d_{(\omega_i)}(x,\rho_{(\omega_i)})+ |t_i-t|$. If $x\in\t_{(\omega_i)}$ and $y\in\t_{(\omega_j)}$, with $i\not = j$,
we take $\dd_{\t_\infty}(x,y)= d_{(\omega_i)}(x,\rho_{(\omega_i)}) + |t_i-t_j| + d_{(\omega_j)}(\rho_{(\omega_j)},y)$. 
We note that $\t_\infty$ is a non-compact $\R$-tree. By convention, $\t_\infty$
is rooted at $0$. The tree $\t_\infty$ is equipped with a volume measure, which is defined as the sum
of the volume measures on the trees $\t_{(\omega_i)}$, $i\in I\cup J$.

We can also define labels on $\t_\infty$. The label $\Lambda_x$
of $x\in\t_\infty$ is defined by $\Lambda_x= \w(t)$ if $x=t$ belongs to the spine $[0,\infty)$, and 
$\Lambda_x=\ell_x(\omega_i)$ if $x$ belongs to $\t_{(\omega_i)}$, for some 
$i\in I\cup J$. Note that the mapping $x\mapsto \Lambda_x$ is continuous (use property \eqref{conti-spine}
to check continuity at points of the spine).

For our purposes, it is important to define an order structure on $\t_\infty$. To this end, we introduce a
``clockwise exploration'' of $\t_\infty$, which is defined as follows. Write $\beta_{u-}$ and $\beta'_{u-}$ for the respective left limits at $u$ of the
functions $u\mapsto \beta_u$ and $u\mapsto \beta'_u$ introduced in (iv) above, with the convention $\beta_{0-}=\beta'_{0-}=0$. 
Then, for every $s\geq 0$, there is a unique $u\geq 0$ such that $\beta_{u-}\leq s\leq \beta_{u}$, and:
\begin{description}
\item[$\bullet$] Either we have $u=t_i$ for some $i\in I$ (then $\sigma(\omega_i)=\beta_{t_i}-\beta_{t_i-}$), and we set
$\ee^+_s:= p_{(\omega_i)}(s-\beta_{t_i-})$.
\item[$\bullet$] Or there is no such $i$ and we set $\ee^+_s=u$.
\end{description}
We define similarly $(\ee^-_s)_{s\geq 0}$. For every $s\geq 0$, there is a unique $u\geq 0$ such that $\beta'_{u-}\leq s\leq \beta'_{u}$, and:
\begin{description}
\item[$\bullet$] Either we have $u=t_i$ for some $i\in J$, and we set
$\ee^-_s:= p_{(\omega_i)}(\beta'_{t_i}-s)(\,= p_{(\check\omega_i)}(s-\beta'_{t_i-}))$.
\item[$\bullet$] Or there is no such $i$ and we set $\ee^-_s=u$.
\end{description}
Informally, $(\ee^+_s)_{s\geq 0}$ and $(\ee^-_s)_{s\geq 0}$ correspond to the exploration of the left and right side
of the tree $\t_\infty$ respectively. 
Noting that $\ee^+_0=\ee^-_0=0$, we define $(\ee_s)_{s\in \R}$ 
by setting
$$\ee_s:=\left\{\begin{array}{ll}
\ee^+_s\quad&\hbox{if }s\geq 0,\\
\ee^-_{-s}\quad&\hbox{if }s\leq 0.
\end{array}
\right.
$$
It is straightforward to verify that the mapping $s\mapsto \ee_s$ from $\R$ onto $\t_\infty$ is continuous.
We also note that the volume measure on $\t_\infty$ is  the pushforward of Lebesgue
measure on $\R$ under the mapping $s\mapsto \ee_s$.

This exploration process allows us to define intervals on
$\t_\infty$, in a way similar to what we did in Section \ref{sna-tra}. 
Let us make the convention that, if $s>t$, the ``interval'' $[s,t]$ is defined by $[s,t]=[s,\infty)\cup (-\infty,t]$. 
Then, for every $x,y\in\t_\infty$, such that $x\not =y$, there is a smallest interval $[s,t]$, with $s,t\in\R$, such that
$\ee_s=x$ and $\ee_t=y$, and we define  
$$[x,y]:=\{\ee_r:r\in[s,t]\}.$$
Note that we have typically $[x,y]\not =[y,x]$. Of course, we take $[x,x]=\{x\}$. We sometimes also use the self-evident notation $]x,y[$. For $x\in\t_\infty$, we finally define $[x,\infty)=\{\ee_r:r\in[s,\infty)\}$, where $s$ is the largest real such that
$\ee_s=x$, and we define $(-\infty,x]$ in a similar manner. Note that $[x,\infty)\cap (-\infty,x]=\llbracket x,\infty\llbracket$
is the range of the geodesic ray starting from $x$ in $\t_\infty$.

\medskip
\noindent{\it The finite spine case.} 
It will also be useful to consider the case where $\w=(\w(t))_{0\leq t\leq\zeta}$ is a finite path in $\W$
with positive lifetime $\zeta$, and
$\pp$ and $\pp'$ are now point measures supported on $(0,\zeta]\times \S$. We then assume that the obvious adaptations
of properties (i)---(v) hold, and in particular (iv) is replaced by  
\begin{description}\item[\rm(iv)'] the functions 
$u\mapsto \beta_u:= \sum_{i\in I} \mathbf{1}_{\{t_i\leq u\}}\,\sigma(\omega_i)\;,\ u\mapsto \beta'_u:= \sum_{i\in J} \mathbf{1}_{\{t_i\leq u\}}\,\sigma(\omega_i)$
 take finite values and are monotone increasing
on $[0,\zeta]$.
\end{description}
The same construction yields a (labeled) compact $\R$-tree $(\t,(\Lambda_v)_{v\in\t})$, with a spine represented by the interval $[0,\zeta]$. 
The distance on $\t$ is denoted by $\dd_\t$ and the labels are defined in exactly the same way as in the infinite spine case. The tree
$\t$ has a cyclic order structure induced by a clockwise exploration function $\ee_s$,
which is conveniently defined on the interval $[0,\beta_\zeta+\beta'_\zeta]$: Informally
$(\ee_s,s\in [0,\beta_\zeta])$ is obtained by concatenating the mappings $p_{(\omega_i)}$
for all atoms $(t_i,\omega_i)$ of $\pp$, in the increasing order of the $t_i$'s, and
$(\ee_s,s\in [\beta_\zeta,\beta_\zeta+\beta'_\zeta])$ is obtained by concatenating the mappings $p_{(\omega_j)}$
for all atoms $(t_j,\omega_j)$ of $\pp'$, in the decreasing order of the $t_j$'s  (in particular, $\ee_0=\ee_{\beta_\zeta+\beta'_\zeta}$ is the root or bottom of the spine
and $\ee_{\beta_\zeta}$ is the top of the spine). In order to define ``intervals'' on the tree $\t$, we now make the 
convention that, if $s,t\in [0,\beta_\zeta+\beta'_\zeta]$ and $t<s$, $[s,t]=[s,\beta_\zeta+\beta'_\zeta]\cup  [0,t]$. 
In that setting, we again refer to $(\w,\pp,\pp')$ as a (finite spine) coding triple.

\smallskip
Here, in contrast with the infinite spine case, we can also represent the labeled tree $(\t,(\Lambda_v)_{v\in\t})$ by a snake trajectory $\omega\in\S_{\w(0)}$ such that $\t_{(\omega)}=\t$, which is defined as follows.
The duration $\sigma(\omega)$ is equal to $\beta_\zeta+\beta'_\zeta$, and, for every $s\in[0,\beta_\zeta+\beta'_\zeta]$, 
the finite path $\omega_s$ is such that $\zeta_{(\omega_s)}=\dd_\t(\ee_0,\ee_s)$ and $\wh\omega_s=\Lambda_{\ee_s}$ 
(by a remark in Section \ref{sna-tra}, this completely determines $\omega$). The snake trajectory $\omega$ obtained in this way will be denoted by $\omega=\Omega(\w,\pp,\pp')$. We note that the triple $(\w,\pp,\pp')$ contains more information than
$\omega$: Roughly speaking, in order to recover this triple from $\omega$, we need to know $s_0\in(0,\sigma(\omega))$,
such that the ancestral line of $p_{(\omega)}(s_0)$ in the genealogical tree of $\omega$
corresponds to the spine. 

It will be useful to consider the time-reversal operation on (finite spine) coding triples satisfying our assumptions,
which is defined by
\begin{equation}
\label{time-rev}
\mathbf{TR}\Big((\w(t))_{0\leq t\leq \zeta},\sum_{i\in I}\delta_{(t_i,\omega_i)}, \sum_{j\in J}\delta_{(t'_j,\omega'_j)}\Big)
:=\Big((\w(\zeta-t))_{0\leq t\leq \zeta},\sum_{j\in J}\delta_{(\zeta-t'_j,\omega'_j)},\sum_{i\in I}\delta_{(\zeta-t_i,\omega_i)}\Big).
\end{equation}
We note that the labeled trees
associated with the coding triples in the left and right sides of \eqref{time-rev} are identified via an isometry that preserves 
labels and
intervals, but the roles of the top and the bottom of the spine are interchanged. 
\medskip

In Section \ref{dis-rela} below we investigate relations between different distributions
on coding triples, and in Section \ref{coding-metric} we explain how to go from (random) coding triples
to the random metric spaces of interest in this work.

\smallskip
\noindent{\it Important remark.} Later, when we speak about the tree associated with a coding triple (as we just defined 
both in the finite and in the infinite spine case), it will always be understood that this includes the labeling on the
tree and the clockwise exploration, which is needed to make sense of intervals on the tree.

\medskip
\noindent{\it Spine decomposition under $\N_a$.}
Let $a>0$. We conclude this section with a result connecting the measure $\N_a$  with a (finite spine) coding triple. 
Arguing under $\N_a(\dd\omega)$, for every $r\in (0,\sigma)$, we can define two point measures $\pp_{(r)}$ and $\pp'_{(r)}$ that account for
the (labeled) subtrees branching off the ancestral line of $p_{(\omega)}(r)$ in the genealogical tree $\t_{(\omega)}$. Precisely, if $r$ is fixed, we consider all connected components
$(u_i,v_i)$, $i\in I$, of the open set $\{s\in[0,r]: \zeta_s(\omega)>\min_{t\in[s,r]}\zeta_t(\omega)\}$, and for each $i\in I$, we define
a snake trajectory $\omega^i$ by setting $\sigma(\omega^i)=v_i-u_i$ and,
for every $s\in[0,\sigma(\omega^i)]$,
$$\omega^i_s(t):=\omega_{u_i+s}(\zeta_{u_i}(\omega)+t)\;,\quad\hbox{for } 
0\leq t\leq \zeta_{(\omega^i_s)}:=\zeta_{u_i+s}(\omega)-\zeta_{u_i}(\omega).$$
Note that $\omega^i\in \S_{\wh \omega_{u_i}}$, and $\wh \omega_{u_i}=\omega_r(\zeta_{u_i})$ by the snake property.
We then set $\pp_{(r)}=\sum_{i\in I} \delta_{(\zeta_{u_i},\omega^i)}$. To define $\pp'_{(r)}$, we proceed in a very similar manner,
replacing the interval $[0,r]$ by $[r,\sigma]$.

Recall our notation $(L^0_s)_{s\in [0,\sigma]}$ for the exit local time at level $0$. Let $M_p(\R_+\times \S)$ stand for the set
of all point measures on $\R_+\times \S$.

\begin{proposition}
\label{Palm}
Let $a>0$ and let $Y=(Y_t)_{0\leq t\leq T^Y}$ stand for a linear Brownian motion started from $a$ and stopped at 
its first hitting time of $0$. Conditionally given $Y$, let $\mm$ and $\mm'$ be two independent Poisson point 
measures on $\R_+\times \S$ with intensity
$$2\,\mathbf{1}_{[0,T^Y]}(t)\,\dd t\,\N_{Y_t}(\dd \omega).$$
Then, for any nonnegative measurable function $F$ on $\W\times M_p(\R_+\times \S)^2$, we have
$$\N_a\Big(\int_0^\sigma \dd L^0_r\,F(W_r,\pp_{(r)},\pp'_{(r)})\Big)= \E\Big[F(Y,\mm,\mm')\Big].$$
\end{proposition}

It is straightforward to verify that $\N_a(\dd\omega)$ a.e., for every $r\in(0,\sigma)$, $(W_r,\pp_{(r)},\pp'_{(r)})$
is a coding triple in the sense of the previous discussion (finite spine case), and $\Omega(W_r,\pp_{(r)},\pp'_{(r)})=\omega$.

\proof We may assume that $F(W_r,\pp_{(r)},\pp'_{(r)})=F_1(W_r)F_2(\pp_{(r)})F_3(\pp'_{(r)})$ where 
$F_1$ is defined on $\W$ and $F_2$ and $F_3$ are defined on $M_p(\R_+\times \S)$.
Let $(\tau^0_r)_{r\geq 0}$ be the inverse local time defined by 
$\tau^0_r=\inf\{s\geq 0:L^0_s \geq r\}$. Then,
$$\N_a\Big(\int_0^\sigma \dd L^0_r\,F_1(W_r)F_2(\pp_{(r)})F_3(\pp'_{(r)})\Big)
=\int_0^\infty dr\,\N_a\Big(\mathbf{1}_{\{\tau^0_r<\infty\}}\,F_1(W_{\tau^0_r})F_2(\pp_{(\tau^0_r)})F_3(\pp'_{(\tau^0_r)})\Big).$$
We may now apply the strong Markov property of the Brownian snake \cite[Theorem IV.6]{Zurich}, noting that
$F_1(W_{\tau^0_r})F_2(\pp_{\tau^0_r})$ is measurable with respect to the past up to time $\tau^0_r$. Using
also \cite[Lemma V.5]{Zurich}, we get, for every $r>0$,
$$\N_a\Big(\mathbf{1}_{\{\tau^0_r<\infty\}}\,F_1(W_{\tau^0_r})F_2(\pp_{(\tau^0_r)})F_3(\pp'_{(\tau^0_r)})\Big)
=\N_a\Big(\mathbf{1}_{\{\tau^0_r<\infty\}}\,F_1(W_{\tau^0_r})F_2(\pp_{(\tau^0_r)})\,\P^{(W_{\tau^0_r})}(F_3)\Big),$$
where, for any finite path $\w$, we write $\P^{(\w)}$ for the distribution of a Poisson point measure
on $\R_+\times M_p(\S)$
with intensity $2\,\mathbf{1}_{[0,\zeta_{(\w)}]}(t)\,\dd t\,\N_{\w(t)}(\dd \omega)$. From the preceding two
displays, we arrive at
$$\N_a\Big(\int_0^\sigma \dd L^0_r\,F_1(W_r)F_2(\pp_{(r)})F_3(\pp'_{(r)})\Big)
= \N_a\Big(\int_0^\sigma \dd L^0_r\,F_1(W_r)F_2(\pp_{(r)})\,\P^{(W_{r})}(F_3)\Big).$$
By the invariance of the excursion measure $\N_a$ under time-reversal (this is an immediate 
consequence of the similar property for the It\^o measure of Brownian excursions), the right-hand side of
the last display is also equal to 
$$\N_a\Big(\int_0^\sigma \dd L^0_r\,F_1(W_r)F_2(\check\pp'_{(r)})\,\P^{(W_{r})}(F_3)\Big),$$
where we write $\check\pp'_{(r)}(\dd s\dd\omega)$ for the image of $\pp'_{(r)}(\dd s\dd\omega)$ under the mapping $(s,\omega)\mapsto(s,\check\omega)$.
Then the same application of the strong Markov property shows that this equals
$$\N_a\Big(\int_0^\sigma \dd L^0_r\,F_1(W_r)\,\P^{(W_r)}(F_2)\,\P^{(W_{r})}(F_3)\Big).$$
Finally, the first-moment formula in \cite[Proposition V.3]{Zurich} shows that this quantity is 
also equal to
$$\E\Big[ F_1(Y)\,\P^{(Y)}(F_2)\,\P^{(Y)}(F_3)\Big]$$
with the notation of the proposition. This completes the proof. \endproof

Proposition \ref{Palm} will allow us to relate properties valid under $\N_a$
to similar properties for the triple $(Y,\mm,\mm')$. Let us illustrate this on an example that will be useful later.
Recall from \eqref{formu-exit} that the exit measure $\z_0$ satisfies
$$\z_0=\lim_{\ve \da 0} \frac{1}{\ve^2} \int_0^\sigma \dd s\,\mathbf{1}_{\{\tau_0(W_s)=\infty,\, \wh W_s<\ve\}},\quad \N_a\hbox{ a.e.}$$
Replacing the limit by a liminf, we may assume that $\z_0(\omega)$ is defined for every $\omega\in\W$. 
It is a simple matter to verify that, $\N_a$ a.e. for every $r$ such that
$\zeta_r=\tau_0(W_r)$, we have
\begin{equation}
\label{tech-exit1}
\z_0=\int \pp_{(r)}(\dd t\dd \varpi)\,\z_0(\varpi) + \int \pp'_{(r)}(\dd t\dd \varpi)\,\z_0(\varpi).
\end{equation}
Then, if we define $\z^Y:=\int \mm(\dd t\dd \varpi)\,\z_0(\varpi) + \int \mm'(\dd t\dd \varpi)\,\z_0(\varpi)$, we deduce
from Proposition \ref{Palm} that, for any nonnegative measurable function $\varphi$ on $\R_+$, we have 
$\N_a(\z_0\,\varphi(\z_0))= \E[ \varphi(\z^Y)]$.
More generally, set $\tau^Y_u=\inf\{t\geq 0:Y_t=u\}$ and $\z^Y_u:=\int \mm(\dd t\dd \varpi)\,\mathbf{1}_{\{t<\tau^Y_u\}}\,\z_u(\varpi) + \int \mm'(\dd t\dd \varpi)\,\mathbf{1}_{\{t<\tau^Y_u\}}\,\z_0(\varpi)$, for every $u\in (0,a)$. Then  we have
\begin{equation}
\label{tech-exit2}
\N_a(\z_0\,\varphi(\z_0,\z_{u_1},\ldots\z_{u_p}))= \E[ \varphi(\z^Y,\z^Y_{u_1},\ldots,\z^Y_{u_p})],
\end{equation}
for every $0<u_1<\cdots<u_p<a$ and any nonnegative measurable function $\varphi$ on $\R_+^{p+1}$. 

\section{Distributional relations between coding processes}
\label{dis-rela}

\subsection{Some explicit distributions}

Let us introduce the function
$$\psi(x)=\frac{2}{\sqrt{\pi}}(x^{1/2} + x^{-1/2}) - 2(x+\frac{3}{2})\,e^x\,\hbox{erfc}(\sqrt{x}),\quad x>0.$$
Note that $\psi(x)=x^{-1}\chi_3(x)$ in the notation of the Appendix below, and thus, by formula (A.3) there,
\begin{equation}
\label{laplace1}
\int_0^{\infty} e^{-\lambda x}\,x\psi(x)\,\dd x = (1+\sqrt{\lambda})^{-3},\quad \lambda\geq 0.
\end{equation}
%Since $x\psi(x)=\gamma_1*\gamma_1*\gamma_1(x)$, with $\gamma_1$ defined in (A.0), it is clear that
%$\psi$ takes positive values. 
Furthermore one checks from the explicit formula for $\psi$ that
$\psi(x)= \frac{2}{\sqrt{\pi}}x^{-1/2} + O(1)$ as $x\to 0$, and $\psi(x)=\frac{3}{2\sqrt{\pi}}x^{-5/2} +O(x^{-7/2})$ as $x\to \infty$. 

\begin{proposition}
\label{densityuniform}
{\rm(i)} Let $a>0$. The density of $\z_0$ under $\N_a(\cdot\cap\{ \z_0\not =0\})$ is
\begin{equation}
\label{density-exit}
h_a(z):= \Big(\frac{3}{2a^2}\Big)^2 \, \psi(\frac{3z}{2a^2}),\quad z>0.
\end{equation}
{\rm(ii)} For every $z>0$ and $a>0$, set
$$p_z(a):=2\,\Big(\frac{3}{2}\Big)^{3/2}\,\sqrt{\pi}\,z^{3/2}a^{-4}\,\psi(\frac{3z}{2a^2}).$$
Then, $a\mapsto p_z(a)$ defines a probability density on $(0,\infty)$, and 
for every nonnegative measurable function $g$ on $[0,\infty)$,
$$z^{-2}\N^{*,z}\Big( \int_0^\sigma \dd s\,g(\wh W_s)\Big)
= \int_0^\infty \dd a\,p_z(a)\,g(a).$$
\end{proposition}

\rem The construction of Brownian disks in \cite{Disks} allows us to interpret
(ii) by saying that $p_z$ is the density of the distribution of the distance from the distinguished point to
the boundary in a free pointed Brownian disk with perimeter $z$. See Proposition \ref{pointed-disk} below.

\proof (i) From \eqref{Lap-exit}, we have for $\lambda\geq 0$,
$$\N_a\Big(1-\exp (-\lambda \z_0)\Big)= \Big(a\sqrt{2/3} + \lambda^{-1/2}\Big)^{-2}
=\frac{3}{2a^2}\,\Big(1+ \Big(\frac{2a^2\lambda}{3}\Big)^{-1/2}\Big)^{-2},$$
and in particular $\N_a(\z_0\not = 0)=\frac{3}{2a^2}$ in agreement with \eqref{hittingpro}. 
On the other hand, 
by formula (A.4) in the Appendix,
$$\int_0^{\infty} (1-e^{-\lambda x})\,\psi(x)\,\dd x
= (1+\lambda^{-1/2})^{-2}.$$
Part (i) follows by comparing the last two displays.

\smallskip
\noindent(ii) From \eqref{formu-exit} and \eqref{approxz*0}, we get the existence of a measurable function $\Gamma$ on $\S$ such that 
$\z_0=\Gamma(\tr_0(\omega))$, $\N_a(\dd \omega)$ a.e., for any $a>0$, 
and also $\z^*_0=\Gamma(\omega)=\Gamma(\omega^{[t]})$
for every $t\in[0,\sigma(\omega)]$, $\N^*(\dd\omega)$ a.e.
By applying the re-rooting formula \eqref{reroot-for} with $G(\omega)=f(\Gamma(\omega))g(\omega_0)$, where $f$ and $g$ are nonnegative real functions, we get
$$\N^*\Big(f(\z^*_0) \int_0^\sigma \dd s\,g(\wh W_s)\Big)=2\,\int_0^\infty \dd a\, g(a)\,\N_a\Big(\z_0f(\z_0)\Big).$$
The left-hand side can be written as
$$\sqrt{\frac{3}{2\pi}}\,\int_0^\infty \dd z\,z^{-5/2}\,f(z)\,\N^{*,z}\Big( \int_0^\sigma \dd s\,g(\wh W_s)\Big).$$
On the other hand, part (i) allows us to rewrite the right-hand side as 
$$2\,\int_0^\infty \dd a\, g(a)\,\int_0^\infty \dd z\,h_a(z)\,z\,f(z)
=2\,\int_0^\infty \dd z\, z\,f(z)\int_0^\infty \dd a\,h_a(z)\,g(a).$$
By comparing the last two displays, we get, $\dd z$ a.e.,
$$z^{-2}\N^{*,z}\Big( \int_0^\sigma \dd s\,g(\wh W_s)\Big)
= 2\sqrt{\frac{2\pi}{3}}\,z^{3/2}\,\int_0^\infty \dd a\,h_a(z)\,g(a)= \int_0^\infty \dd a\,p_z(a)\,g(a),$$
where $p_z(a)$ is as in the proposition. A scaling argument shows that this identity indeed holds for every $z>0$.
Since $z^{-2}\N^{*,z}(\sigma)=1$,
$p_z$ is a probability density, which may also be checked directly.
\endproof

\subsection{A distributional identity for coding triples}
\label{sec-ident}

As we already explained, coding triples will be used to construct the random metric spaces of interest 
in this work. 
The relevant case for the forthcoming construction of the
infinite-volume Brownian disk with perimeter $z>0$ may be described as follows: we let $R=(R_t)_{t\in[0,\infty)}$ be a three-dimensional Bessel process started from $0$,  we assume that,
conditionally on $R$, $\pp$ and $\pp'$ are independent Poisson measures with intensity 
$2\,\dd t\,\N_{R_t}(\dd \omega)$, and finally we condition on the event $\z=z$, where 
$\z$ denotes the total exit measure at $0$ of the atoms of $\pp$ and $\pp'$. In Section \ref{cod-disk},
we will give a precise meaning to this conditioning and obtain a conditional distribution 
$\Theta_z$ on coding triples, which plays an important role in the next sections. Before doing that, we need to develop
certain preliminary tools, and we first recall 
special cases of a well-known time-reversal property for Bessel processes. Let $R$ be as above and let $X$ be a Bessel process
of dimension $9$ started from $0$. Then, for every $a>0$,
\begin{itemize}
\item[\rm(a)] If $\mathbf{L}_a:=\sup\{t\geq 0: R_t=a\}$, the process $(R_{\mathbf{L}_a-t}, 0\leq t\leq \mathbf{L}_a)$ is distributed 
as a linear Brownian motion started from $a$ and stopped upon hitting $0$.
\item[\rm(b)] If $\mathrm{L}_a:=\sup\{t\geq 0: X_t=a\}$, the process $(X_{\mathrm{L}_a-t}, 0\leq t\leq \mathrm{L}_a)$ is distributed 
as a Bessel process of dimension $-5$  started from $a$ and stopped upon hitting $0$.
\end{itemize}
Both (a) and (b) are special cases of a more general result for Bessel processes, which is itself a consequence of 
Nagasawa's time-reversal theorem (see \cite[Theorem VII.4.5]{RY}, and \cite[Exercise XI.1.23]{RY} for the 
case of interest here, and note that part (a) is due to Williams \cite{Wil}). As a consequence of (a) and (b), one
gets that, for every $a>0$, the process $(R_{\mathbf{L}_a+t}, t\geq 0)$ is independent of $(R_t,0\leq t\leq \mathbf{L}_a)$, and 
similarly $(X_{\mathrm{L}_a+t}, t\geq 0)$ is independent of $(X_t,0\leq t\leq \mathrm{L}_a)$. This property is used implicitly in what follows.

Let us introduce some notation needed for the technical results that follow. 
We fix $a>0$ and consider 
a triple $(Y,\mm,\mm')$ distributed as in Proposition \ref{Palm}:
$Y=(Y_t,0\leq t\leq T^Y)$ is a linear Brownian motion started from $a$
and stopped at the first time it hits $0$, and,
conditionally on $Y$, $\mm$ and $\mm'$ are independent Poisson point measures on $\R_+\times \S$
with intensity
$2\,\mathbf{1}_{[0,T^Y]}(t)\,\dd t\,\N_{Y_t}(\dd \omega)$.
We also introduce the point measures $\wt\mm$ and $\wt\mm'$ obtained by
truncating the atoms of $\mm$ and $\mm'$ at level $0$. More precisely, for any nonnegative measurable
function $\Phi$ on $\R_+\times \S$, we set
$$\int \wt\mm(\dd t\dd \omega)\,\Phi(t,\omega):= \int \mm(\dd t\dd \omega)\,\Phi(t,\tr_0(\omega))$$
and $\wt\mm'$ is defined similarly from $\mm'$. We will be interested in the triple $(Y,\wt\mm,\wt\mm')$, which we may view as
a coding triple in the sense of Section \ref{coding-infinite}.

We define
$$\z^Y = \int \mm(\dd t\dd \omega) \z_0(\omega) + \int \mm'(\dd t\dd \omega)\z_0(\omega)$$
in agreement with the end of Section \ref{coding-infinite}.
We saw that,
for any nonnegative measurable function $\varphi$ on $\R_+$, we have $\E[ \varphi(\z^Y)]=\N_a(\z_0\,\varphi(\z_0))$.
It then follows from Proposition \ref{densityuniform}$\,$(i) that the distribution 
of $\z^Y$ has density $z\,h_a(z)$. 
%By construction,
%$$\E[e^{-\lambda \z^Y}]= \E\Big[\exp\Big(-4 \int_0^{T^Y} \dd t\,\N_{Y_t}(1- e^{-\lambda \z_0})\Big)\Big]
%=\E\Big[\exp\Big( -6\int_0^{T^Y} \dd t \Big((\frac{2\lambda}{3})^{-1/2} + Y_t \Big)^{-2}\Big)\Big]$$
%using formula \eqref{Lap-exit}. 
%The latter quantities are equal to
%$$\E_v\Big[\exp\Big(-6\int_0^{T_u} \frac{\dd t}{B_t^2}\Big)\Big]$$
%where $u=(2\lambda/3)^{-1/2}$, $v=a+u$, $B$ stands for a linear Brownian motion that starts from $v$ under the
%probability measure $\P_v$, and $T_u=\inf\{t\geq 0:B_t=u\}$. The quantity in the last display is equal to $(u/v)^3$ (a simple
%way to get this is to apply the optional stopping theorem to the martingale
%$B^{-3}_{t\wedge T_u}\exp(-6\int_0^{t\wedge T_u} B_s^{-2}\, \dd s)$). We have thus proved that
%\begin{equation}
%\label{Laplace-ex}
%\E[e^{-\lambda \z^Y}]= \Big(\frac{u}{v}\Big)^3= \Big(1+a\sqrt{2\lambda/3}\Big)^{-3}.
%\end{equation}
%Recalling \eqref{laplace1}, we get that the distribution of $\z^Y$ has density
%\begin{equation}
%\label{Laplace-ex2}
%\Big(\frac{3}{2a^2}\Big)^2 z\,\psi\Big(\frac{3z}{2a^2}\Big)=z\,h_{a}(z).
%\end{equation}
In particular, we have
\begin{equation}
\label{Laplace-ex}
\E[e^{-\lambda \z^Y}]=\int_0^\infty z\,h_a(z)\,e^{-\lambda z}\,\dd z=\Big(1+a\sqrt{2\lambda/3}\Big)^{-3},
\end{equation}
by \eqref{laplace1} and \eqref{density-exit}.

We write $(\check\Theta^{(a)}_z)_{z>0}$ for a regular version of the conditional distributions 
of the triple $(Y,\wt\mm,\wt\mm')$ knowing that $\z^Y=z$. The collection
$(\check\Theta^{(a)}_z)_{z>0}$ is well defined only up to a set of values of $z$ 
of zero Lebesgue measure, but we will see later how to make a canonical choice of
this collection.

Let us also fix $r>0$.
We next consider a triple $(V,\nn,\nn')$, where
\begin{itemize}
\item[$\bullet$] $V=(V_t,0\leq t\leq T^V)$ is distributed as a Bessel process of dimension $-5$ started from $r+a$
and stopped at the first time it hits $r$. 
\item[$\bullet$] Conditionally on $V$, $\nn$ and $\nn'$ are independent Poisson point measures on $\R_+\times \S$
with intensity
$$2\,\mathbf{1}_{[0,T^V]}(t)\,\dd t\,\N_{V_t}(\dd \omega\cap \{W_*>0\}),$$
where we recall the notation $W_*(\omega)=\min\{\wh W_s(\omega):0\leq s\leq \sigma(\omega)\}$.
\end{itemize}
We write  $\wt\nn$ and $\wt\nn'$ for the point measures obtained by
truncating the atoms of $\nn$ and $\nn'$ at level $r$, in the same way as
$\wt\mm$ was defined above from $\mm$ by truncation at level $0$. We also introduce the exit measure
$$\z^V= \int \nn(\dd t\dd \omega) \z_r(\omega) + \int \nn'(\dd t\dd \omega)\z_r(\omega).$$
As we will see in the next proof, the distributions of $\z^V$ and $\z^Y$ are related by the formula
\begin{equation}
\label{relation-density}
\E[h(\z^V)]= \Big(\frac{r+a}{r}\Big)^3\,\E[h(\z^Y)\,e^{-\frac{3}{2r^2}\z^Y}].
\end{equation}
In particular, the distribution of $\z^V$ also has a positive density on $(0,\infty)$.

We let $\vartheta_r$ stand for the obvious shift that maps snake trajectories with initial point $x$ to snake
trajectories with initial point $x-r$. If $\mu=\sum_{i\in I}\delta_{(t_i,\omega_i)}$ is a point
measure on $\R_+\times \S$, we also write 
$\vartheta_r\mu=\sum_{i\in I}\delta_{(t_i,\vartheta_r\omega_i)}$, by abuse of notation.

\begin{proposition}
\label{ident-cond-dis}
The collection $(\check\Theta^{(a)}_z)_{z>0}$ is a regular version of the  conditional distributions of $(V-r,\vartheta_r\wt\nn,\vartheta_r\wt\nn')$ knowing that $\z^V=z$.
\end{proposition}

\smallskip
In other words, the conditional distribution of $(V-r,\vartheta_r\wt\nn,\vartheta_r\wt\nn')$ knowing that $\z^V=z$
coincides with the conditional distribution of  $(Y,\wt\mm,\wt\mm')$ knowing that $\z^Y=z$.
In particular, the conditional distribution of $(V-r,\vartheta_r\wt\nn,\vartheta_r\wt\nn')$ knowing that $\z^V=z$ does not depend on $r$,
which is by no means an obvious fact.

\smallskip
\begin{proof} 
Recall our notation $M_p(\R_+\times \S)$ for the set
of all point measures on $\R_+\times \S$.
As in the proof of Proposition \ref{Palm}, if $\w$ is a finite path taking nonnegative values, we
write $\P^{(\w)}(\dd \mu)$ for the probability measure on $M_p(\R_+\times \S)$ which is 
the distribution of a Poisson point measure on $\R_+\times \S$
with intensity 
$2\,\mathbf{1}_{[0,\zeta_{(\w)}]}(t)\,\dd t\,\N_{\w(t)}(\dd \omega)$.
Denoting the generic element of $M_p(\R_+\times \S)\times M_p(\R_+\times \S)$ by $(\mu,\mu')$, we have
the formula
\begin{equation}
\label{condis-11}
\P^{(\w)}\otimes \P^{(\w)}(\mu(W_*\leq 0)\!=\!\mu'(W^*\leq 0)\!=\!0)=\exp\Big(-4\!\int_0^{\zeta_{(\w)}} \!\dd t\,\N_{\w(t)}(W_*\leq 0)\Big)
=\exp\Big(-6\!\int_0^{\zeta_{(\w)}}\! \frac{\dd t}{\w(t)^2}\Big),
\end{equation}
where in the left-hand side we abuse notation by writing $\mu(W_*\leq 0)$ instead of  $\mu(\{(t,\omega)\in\R_+\times \S: W_*(\omega)\leq 0\})$.

We then introduce a random finite path $(U_t)_{0\leq t\leq T^U}$, which is distributed as a linear Brownian motion
started from $r+a$ and stopped when hitting $r$ (so $U-r$ has the same distribution as $Y$).
Let $\pp$ and $\pp'$ be random elements of $M_p(\R_+\times \S)$ such that the conditional distribution of the pair $(\pp,\pp')$ given 
$U$ is $\P^{(U)}\otimes \P^{(U)}(\dd \mu\dd\mu')$. Define
$$\z^U= \int \pp(\dd t\dd \omega) \z_r(\omega) + \int \pp'(\dd t\dd \omega)\z_r(\omega),$$
and also write $\wt\pp$, resp.~$\wt\pp'$, for the point measure $\pp$, resp.~$\pp'$, truncated at level $r$. 
Then, the statement of the proposition reduces to showing that the conditional distribution of $(V,\wt\nn,\wt\nn')$
knowing that $\z^V=z$ coincides $\dd z$ a.e.~with the conditional distribution of $(U,\wt\pp,\wt\pp')$
knowing that $\z^U=z$
(an obvious translation argument yields that $(U-r,\vartheta_r\wt\pp,\vartheta_r\wt\pp',\z^U)$ has the same law as $(Y,\wt\mm,\wt\mm',\z^Y)$ above).

The first step of the proof is to verify that, for any nonnegative measurable functions 
$F$ and $G$ defined on $\W$ and on $M_p(\R_+\times \S)^2$ respectively,
\begin{align}
\label{keystep}
&\E\Big[ F(U)G(\pp,\pp')\,\Big|\, \pp(W_*\leq 0)=\pp'(W_*\leq 0)=0\Big]\nonumber\\
&\quad=\E\Big[ F(V)\,\P^{(V)}\otimes\P^{(V)}[G(\mu,\mu')\mid \mu(W_*\leq 0)=\mu'(W_*\leq 0)=0]\Big].
\end{align}

To prove \eqref{keystep}, we first apply \eqref{condis-11} to get
$$
\P\Big(\pp(W_*\leq 0)=\pp'(W_*\leq 0)=0\Big) 
=  \E\Big[\exp\Big(-6\!\int_0^{T^U}\!\!\frac{\dd t}{U_t^2}\Big)\Big] = \Big(\frac{r}{r+a}\Big)^3
$$
where the last equality is easily derived by using It\^o's formula to verify that 
$U^{-3}_{t\wedge T^U}\exp(-6\int_0^{t\wedge T^U} U_s^{-2}\, \dd s)$ is a martingale.
So we have
\begin{align*}
 &\Big(\frac{r}{r+a}\Big)^3\E\Big[ F(U)G(\pp,\pp')\,\Big|\, \pp(W_*\leq 0)=\pp'(W_*\leq 0)=0\Big]\cr
 &\qquad=\E\Big[ F(U)G(\pp,\pp')\,\mathbf{1}_{\{\pp(W_*\leq 0)=\pp'(W_*\leq 0)=0\}}\Big]\cr
&\qquad=\E\Big[ F(U)\,\P^{(U)}\otimes\P^{(U)}[G(\mu,\mu')\,\mathbf{1}_{\{\mu(W_*\leq 0)=\mu'(W_*\leq 0)=0\}}]\Big]\cr
&\qquad=\E\Big[ F(U)\,\exp\Big(-6\int_0^{T^U} \frac{\dd t}{U_t^2}\Big)\;
\P^{(U)}\otimes\P^{(U)}[G(\mu,\mu')\mid \mu(W_*\leq 0)=\mu'(W_*\leq 0)=0]\Big]
\end{align*}
using \eqref{condis-11} in the last equality. 
%to write
%$$\P^{(U)}\otimes\P^{(U)}[\mu(W_*\leq 0)=\mu'(W_*\leq 0)=0]
%=\exp\Big(-6\int_0^{T^U} \frac{\dd t}{U_t^2}\Big).$$
To complete the proof of \eqref{keystep}, we just observe that, by classical absolute continuity relations 
between Brownian motion and Bessel processes, the law of $U$ under the probability measure
$$ \Big(\frac{r+a}{r}\Big)^3 \exp\Big(-6\int_0^{T^U} \frac{\dd t}{U_t^2}\Big)\cdot \P$$
coincides with the law of $V$ under $\P$ (see \cite[Lemma 1]{Bessel} for a short proof). 

Let us complete the proof of the proposition. By a standard property of Poisson measures and the
definition of the pair $(\nn,\nn')$, we have
$$\E[G(\nn,\nn')\mid V]= \P^{(V)}\otimes \P^{(V)}[G(\mu,\mu')\mid \mu(W_*\leq 0)=\mu'(W_*\leq 0)=0].$$
It thus follows from \eqref{keystep} that
$$\E\Big[ F(U)G(\pp,\pp')\,\Big|\, \pp(W_*\leq 0)=\pp'(W_*\leq 0)=0\Big]
= \E[F(V)\,\E[G(\nn,\nn')\!\mid\! V]]= \E[ F(V)\,G(\nn,\nn')].$$
In particular, for any nonnegative measurable function $h$ on $[0,\infty)$, we have
$$\E\Big[ F(U)G(\wt\pp,\wt\pp')h(\z^U)\,\Big|\,  \pp(W_*\leq 0)=\pp'(W_*\leq 0)=0\Big]
= \E[ F(V)\,G(\wt\nn,\wt\nn')\,h(\z^V)].$$
The left-hand side of the last display is equal to
$$\Big(\frac{r+a}{r}\Big)^3\, \E\Big[ F(U)G(\wt\pp,\wt\pp')h(\z^U)\,\exp(-\frac{3}{2r^2}\,\z^U)\Big]$$
because, on one hand, we saw that $\P(\pp(W_*\leq 0)=\pp'(W_*\leq 0)=0)=(\frac{r}{r+a})^3$ and, on the 
other hand, the special Markov property (see e.g. the appendix of \cite{subor}) and \eqref{hittingpro} show that
$$\P(\pp(W_*\leq 0)=\pp'(W_*\leq 0)=0\mid U,\wt\pp,\wt\pp')= \exp(-\frac{3}{2r^2}\,\z^U).$$
We can find nonnegative measurable functions $\varphi_1$ and $\varphi_2$ on $[0,\infty)$ such that
$$\E\Big[ F(U)G(\wt\pp,\wt\pp')\mid \z^U]=\varphi_1(\z^U)\;,\quad  \E[ F(V)\,G(\wt\nn,\wt\nn')\mid \z^V]
=\varphi_2(\z^V)\;,$$
and it follows from the preceding considerations that, for any function $h$,
\begin{equation}
\label{tech-iden}
\Big(\frac{r+a}{r}\Big)^3\, \E[\varphi_1(\z^U)\,h(\z^U)\,\exp(-\frac{3}{2r^2}\,\z^U)] =
\E[\varphi_2(\z^V)h(\z^V)].
\end{equation}
By specializing this identity to the case $F=1$, $G=1$, we get the relation \eqref{relation-density} between the
distributions of $\z^V$ and $\z^Y$ (recall that $\z^U$ has the same distribution as $\z^Y$). But then it also
follows from \eqref{tech-iden} that (for arbitrary $F$ and $G$) we have  
$\E[\varphi_1(\z^V)h(\z^V)]=\E[\varphi_2(\z^V)h(\z^V)]$ for any test function $h$, so that $\varphi_1(\z^V)=\varphi_2(\z^V)$ a.s. and
$\varphi_1(z)=\varphi_2(z)$,
$\dd z$ a.e., which completes the proof.
\end{proof}

For every $z>0$, we let $\Theta^{(a)}_z$ denote the image of $\check\Theta^{(a)}_z$ under the time-reversal transformation
$\mathbf{TR}$ in \eqref{time-rev}.
Since the process $(Y_t,0\leq t\leq T^Y)$
is mapped by time-reversal 
to a three-dimensional Bessel process started from $0$ and stopped at its last passage at $a$
(by property (a) stated at the beginning of the section), we could have defined  $\Theta^{(a)}_z$ 
directly in terms of conditioning a coding triple whose first component is a three-dimensional Bessel process up to a last
passage time. The connection with the discussion at the beginning of this section should then be clear: $\Theta^{(a)}_z$ is the analog
of the probability measure $\Theta_z$ we are aiming at, when the three-dimensional Bessel process is truncated at 
a last passage time.% (in terms of trees, we have a ``finite spine'' instead  of an infinite one).
%Still our 
%presentation turns out to be slightly more convenient for our purposes. 

\subsection{The coding triple of the infinite-volume Brownian disk}
\label{cod-disk}

In this section, we define the probability measures $\Theta_z$, $z>0$, which were introduced informally at
the beginning of the preceding section. Roughly speaking, the idea is to get $\Theta_z$ as the limit of 
$\Theta^{(a)}_z$ as $a\to\infty$. Proposition \ref{condi-dis2} below will also show that, for every $r>0$, the collection $(\Theta_z)_{z>0}$ corresponds to conditional 
distributions of a coding triple whose
first component is a nine-dimensional Bessel process considered after its last passage time at $r$ (compare with Proposition \ref{ident-cond-dis}). The latter
fact is the key to the identification as infinite-volume Brownian disks of the complement of hulls in the Brownian plane. 

We consider a triple  $(X,\ll,\rr)$ such that
$X=(X_t)_{t\geq 0}$ is a nine-dimensional Bessel process started from $0$
and, conditionally on $X$, $\ll$ and $\rr$ are two independent Poisson measures
on $\R_+\times \S$ with intensity
$$2\,\dd t\,\N_{X_t}(\dd \omega \cap\{W_*>0\}).$$
We also set, for every $r>0$,
\begin{equation}
\label{last-pass}
\mathrm{L}_r:=\sup\{t\geq 0:X_t=r\}.
\end{equation}

In what follows, we fix $r>0$, and we shall be interested in atoms $(t,\omega)$
of $\ll$ or $\rr$ such that $t>\mathrm{L}_r$. More precisely, we introduce 
a point measure $\ll^{(r)}$ as the image of 
$$\mathbf{1}_{(\mathrm{L}_r,\infty)}(t)\,\ll(\dd t\,\dd\omega)$$
under the mapping $(t,\omega)\mapsto (t-\mathrm{L}_r,\vartheta_r\omega)$ (where $\vartheta_r$
is the shift operator already used in Proposition \ref{ident-cond-dis}). In a way similar to the previous section, we
define $\wt\ll^{(r)}$ by truncating the atoms of $\ll^{(r)}$ at level $0$ (more precisely,
$\wt\ll^{(r)}$ is the image of $\ll^{(r)}$ under the mapping
$(t,\omega)\mapsto (t,\tr_0(\omega))$). We define similarly $\rr^{(r)}$ and
$\wt\rr^{(r)}$ from the point measure $\rr$. Finally, we set
\begin{align}
\label{peri-hull}
\z^{(r)}&= \int \ll^{(r)}(\dd t\dd\omega) \,\z_0(\omega)
+\int \rr^{(r)}(\dd t\dd\omega) \,\z_0(\omega)\nonumber\\
&= \int \ll(\dd t\dd\omega)\,\mathbf{1}_{(\mathrm{L}_r,\infty)}(t)\,\z_r(\omega)
+ \int \rr(\dd t\dd\omega)\,\mathbf{1}_{(\mathrm{L}_r,\infty)}(t)\,\z_r(\omega)
\end{align}
and we also consider the process $(X^{(r)}_t)_{t\geq 0}$ defined by
$$X^{(r)}_t=X_{\mathrm{L}_r+t}-r.$$
By \cite[Proposition 1.2]{CLG}, the distribution of $\z^{(r)}$ has a density given by
\begin{equation}
\label{density-hull-peri}
k_r(z):= \frac{1}{\sqrt{\pi}}\,3^{3/2} 2^{-1/2}\,r^{-3}\,z^{1/2}\,e^{-\frac{3z}{2r^2}}.
\end{equation}

Our first goal is to verify that the conditional distribution of the triple 
$(X^{(r)},\wt\ll^{(r)},\wt\rr^{(r)})$ knowing that $\z^{(r)}=z$
does not depend on $r$. Note that, for instance, the unconditional 
distribution of $X^{(r)}$ depends on $r$. 

We will deduce the preceding assertion from Proposition \ref{ident-cond-dis},
but to this end a truncation argument is needed. So we consider $a>0$, and we
set
$$\mathrm{L}^{(r)}_a:=\mathrm{L}_{r+a}-\mathrm{L}_r = \sup\{t\geq 0: X^{(r)}_t=a\}.$$
We then set\footnote{Our notation is somewhat misleading since $\ll^{(r+a,\infty)}$ and $\rr^{(r+a,\infty)}$
both depend on $r$ and not only on $r+a$. Since $r$ is fixed in most of this section, this should not
be confusing.}
\begin{align*}
&\ll^{(r,r+a)}(\dd t\dd\omega)=\mathbf{1}_{[0,\mathrm{L}^{(r)}_a]}(t)\,\ll^{(r)}(\dd t\dd\omega),\\
&\ll^{(r+a,\infty)}(\dd t\dd\omega)=\mathbf{1}_{(\mathrm{L}^{(r)}_a,\infty)}(t)\,\ll^{(r)}(\dd t\dd\omega),
\end{align*}
and we define $\rr^{(r,r+a)}$ and $\rr^{(r+a,\infty)}$ in a similar way from $ \rr^{(r)}$.
As previously, we let $\wt\ll^{(r,r+a)}$, $\wt\ll^{(r+a,\infty)}$,  $\wt\rr^{(r,r+a)}$, $\wt\rr^{(r+a,\infty)}$
stand for these point measures truncated at level $0$. We finally set
\begin{align*}
&\z^{(r,r+a)}= \int \ll^{(r,r+a)}(\dd t\dd\omega) \,\z_0(\omega)
+\int \rr^{(r,r+a)}(\dd t\dd\omega) \,\z_0(\omega)\\
&\z^{(r+a,\infty)}= \int \ll^{(r+a,\infty)}(\dd t\dd\omega) \,\z_0(\omega)
+\int \rr^{(r+a,\infty)}(\dd t\dd\omega) \,\z_0(\omega).
\end{align*}
Obviously $\z^{(r)}= \z^{(r,r+a)}+\z^{(r+a,\infty)}$. Also, the random variables $\z^{(r,r+a)}$ and $\z^{(r+a,\infty)}$
are independent, as a consequence
of the independence properties stated at the beginning of Section \ref{sec-ident} after properties (a) and (b).

\begin{lemma}
\label{lem-tech1}
The collection $(\Theta^{(a)}_z)_{z>0}$ is a regular version of the conditional distributions of the triple
$$\Big( (X^{(r)}_t)_{0\leq t\leq \mathrm{L}^{(r)}_a}, \wt\ll^{(r,r+a)}, \wt\rr^{(r,r+a)}\Big)$$
knowing that $\z^{(r,r+a)}=z$.
\end{lemma}

This lemma is merely a reformulation of Proposition \ref{ident-cond-dis}. The point is
that the time-reversed process $ (X_{\mathrm{L}_{r+a}-t})_{0\leq t\leq \mathrm{L}^{(r)}_a}$ is distributed as 
a Bessel process of dimension $-5$ started from $r+a$ and stopped upon hitting $r$
(by property (b) stated at the beginning of Section \ref{sec-ident}). Recalling our notation $\mathbf{TR}$
for the time-reversal operation defined in \eqref{time-rev}, it follows that 
$$\Big(\mathbf{TR}\Big((X^{(r)}_t)_{0\leq t\leq \mathrm{L}^{(r)}_a}, \wt\ll^{(r,r+a)}, \wt\rr^{(r,r+a)}\Big),\z^{(r,r+a)}\Big)$$
has the same distribution as $((V-r,\vartheta_r\wt\nn,\vartheta_r\wt\nn'),\z^V)$, with the notation introduced
before Proposition \ref{ident-cond-dis}.
The result of the lemma now follows from Proposition \ref{ident-cond-dis}.

\smallskip
Since the distribution of $\z^{(r,r+a)}$
is the same as the distribution of $\z^V$ in the preceding section, it has a positive density with respect to 
Lebesgue measure, which we denote by $g_{r,a}(z)$. Recalling that $\z^Y$
has density $zh_a(z)$, \eqref{relation-density} gives the explicit expression
\begin{equation}
\label{densityZX}
g_{r,a}(z)=\Big(\frac{r+a}{r}\Big)^3\,e^{-\frac{3z}{2r^2}}\,z\,h_{a}(z)
\end{equation}
where $h_a$ is defined in \eqref{density-exit}. 

On the other hand, the distribution of $\z^{(r+a,\infty)}$ may be written in the form
$$(1-\ve_{r,a})\,\delta_0(\dd z) + \Upsilon_{r,a}(\dd z)$$
where $\ve_{r,a}\in[0,1]$ and the measure $\Upsilon_{r,a}$ is supported on $(0,\infty)$. Note that 
$$\ve_{r,a}=\Upsilon_{r,a}((0,\infty))=\P(\z^{(r+a,\infty)}>0)=1-\Big(\frac{a}{r+a}\Big)^3$$
where the last equality follows from Lemma 4.2 in \cite{CLG}, using the fact that
$\N_x(0<W_*\leq r)= \frac{3}{2}((x-r)^{-2}-x^{-2})$ for $x>r$. In particular, $\ve_{r,a}\la 0$
as $a\to\infty$.

Recall that $k_r(z)$ denotes the density of $\z^{(r)}$ (cf.~\eqref{density-hull-peri}). Since $\z^{(r)}=\z^{(r,r+a)}+\z^{(r+a,\infty)}$, and $\z^{(r,r+a)}$ and $\z^{(r+a,\infty)}$ are independent, 
the conditional distributions of $\z^{(r+a,\infty)}$ knowing that $\z^{(r)}=z$ 
are defined in a canonical manner by
$$\nu_{r,a}(\dd z'\,|\, z)= \frac{1}{k_r(z)}\Big((1-\ve_{r,a})g_{r,a}(z)\,\delta_0(\dd z') + g_{r,a}(z-z')\,\Upsilon_{r,a}(\dd z')\Big).$$
In particular, we have for every $z>0$,
$$\nu_{r,a}(\{0\}\,|\, z)= \frac{(1-\ve_{r,a})g_{r,a}(z)}{k_r(z)},$$
and the explicit expression \eqref{densityZX} can be used to verify that
$g_{r,a}(z)\la k_r(z)$ as $a\to\infty$. It follows that
\begin{equation}
\label{conv-mass0}
\nu_{r,a}(\{0\}\,|\, z) \build{\la}_{a\to\infty}^{} 1.
\end{equation}

Recall the scaling transformations $\theta_\lambda$ on snake trajectories defined in 
Section \ref{sna-mea}. It will also be useful to consider restriction operators which are
defined as follows. For every $a>0$, $\mathfrak{R}_a$ acts both on $\W^\infty_0\times M_p(\R_+\times\S)^2$
and on $\W_0\times M_p(\R_+\times\S)^2$ by 
\begin{equation}
\label{restric-opera}
\mathfrak{R}_a:\Big(\w, \sum_{i\in I} \delta_{(t_i,\omega_i)}, \sum_{j\in J}\delta_{(t'_j,\omega'_j)}\Big) \mapsto \Big((\w(t))_{t\leq \lambda^{(a)}(\w)},
\sum_{i\in I,t_i\leq \lambda^{(a)}(\w)}\delta_{(t_i,\omega_i)}, \sum_{j\in J,t'_j\leq \lambda^{(a)}(\w)}\delta_{( t_j,\omega'_j)}\Big),
\end{equation}
where $\lambda^{(a)}(\w)=\sup\{t\geq 0:\w(t)\leq a\}$ for $\w\in\mathcal{W}^\infty_0$ or $\w\in\W_0$. 
%For $\lambda >0$ 
%and $\omega\in\S$, $\theta_\lambda(\omega)$ is the element of $\S$ defined by
%$$\theta_\lambda(\omega)_s(t)=\sqrt{\lambda}\,\omega_{s/\lambda^2}(\frac{t}{\lambda})\;, \quad 0\leq t\leq \zeta_{\theta_\lambda(\omega)_s}=\lambda \zeta_{\omega_{s/\lambda^2}}.$$

\begin{proposition}
\label{condi-dis2}
We can find a collection $(\Theta_z)_{z>0}$ of probability measures on $\W^\infty\times M_p(\R_+\times \S)^2$ that does not depend on $r$ and is such that, for every $r>0$, 
$(\Theta_z)_{z>0}$ is a regular version of the 
conditional distributions of the triple
$$(X^{(r)},\wt\ll^{(r)},\wt\rr^{(r)})$$ 
knowing that $\z^{(r)}=z$. This collection is unique if we impose the additional
scaling invariance property: for every $\lambda >0$ and $z>0$, $\Theta_{\lambda z}$ is the
image of $\Theta_z$ under the scaling transformation
$$\Sigma_\lambda:\Big(\w, \sum_{i\in I} \delta_{(t_i,\omega_i)}, \sum_{j\in J}\delta_{(t'_j,\omega'_j)}\Big) \mapsto \Big(\sqrt{\lambda}\w(\cdot/\lambda),
\sum_{i\in I} \delta_{(\lambda t_i,\theta_\lambda(\omega_i))}, \sum_{j\in J}\delta_{(\lambda t'_j,\theta_\lambda(\omega'_j))}\Big).$$
\end{proposition}

\proof Let $r>0$, and let $(\Theta_{z,r})_{z>0}$ be a regular version of the
conditional distributions of $(X^{(r)},\wt\ll^{(r)},\wt\rr^{(r)})$ knowing $\z^{(r)}=z$. Our first goal 
is to verify that $(\Theta_{z,r})_{z>0}$ does not depend on $r$, except possibly on a values of 
$z$ of zero Lebesgue measure. To this end, let $c>0$ and let $G$ be a measurable function 
on $\W\times M_p(\R_+\times \S)^2$ such that $0\leq G\leq 1$. Then, for every $a\geq c$, and
for every nonnegative measurable function $f$ on $(0,\infty)$,
\begin{align*}
\int \dd z\,k_r(z)\,f(z)\Theta_{z,r}(G\circ \mathfrak{R}_c)
&=\E[G((X^{(r)}_t)_{0\leq t\leq \mathrm{L}^{(r)}_c}, \wt\ll^{(r,r+c)},\wt\rr^{(r,r+c)})\,f(\z^{(r)})]\\
&=\E\Big[ \E[G((X^{(r)}_t)_{0\leq t\leq \mathrm{L}^{(r)}_c}, \wt\ll^{(r,r+c)},\wt\rr^{(r,r+c)})\,|\, \z^{(r,r+a)}]\,f(\z^{(r,r+a)}+\z^{(r+a,\infty)})\Big]
\end{align*}
where we use the fact that $\z^{(r+a,\infty)}$ is independent of 
$((X^{(r)}_t)_{0\leq t\leq \mathrm{L}^{(r)}_c}, \wt\ll^{(r,r+c)},\wt\rr^{(r,r+c)}), \z^{(r,r+a)})$ to write the last equality. 
By Lemma \ref{lem-tech1}, we have 
$$\E[G((X^{(r)}_t)_{0\leq t\leq \mathrm{L}^{(r)}_c}, \wt\ll^{(r,r+c)},\wt\rr^{(r,r+c)})\,|\, \z^{(r,r+a)}]
=\Phi(\z^{(r,r+a)})$$
where $\Phi(z)=\Theta^{(a)}_z(G\circ\mathfrak{R}_c)$. Using the explicit distribution of $\z^{(r,r+a)}$
and $\z^{(r+a,\infty)}$, we thus get
\begin{align*}
\int \dd z\,k_r(z)\,f(z)\Theta_{z,r}(G\circ \mathfrak{R}_c)
&=\int\dd y\, g_{r,a}(y)\int ((1-\ve_{r,a})\delta_0+\Upsilon_{r,a})(\dd y')\,f(y+y')\,\Theta^{(a)}_y(G\circ\mathfrak{R}_c)\\
&=\int \dd z\,k_r(z)\,f(z)\,\int \nu_{r,a}(\dd z'\,|\, z)\,\Theta^{(a)}_{z-z'}(G\circ\mathfrak{R}_c).
\end{align*}
It follows that we have, $\dd z$ a.e.,
$$\Theta_{z,r}(G\circ \mathfrak{R}_c)=\int \nu_{r,a}(\dd z'\,|\, z)\,\Theta^{(a)}_{z-z'}(G\circ\mathfrak{R}_c)
= \nu_{r,a}(\{0\}\,|\,z)\,\Theta^{(a)}_z(F) + \kappa_{r,a}(z)$$
where the ``remainder'' $\kappa_{r,a}(z)$ is nonnegative and bounded above by $1-\nu_{r,a}(\{0\}\,|\, z)$. Specializing to integer 
values of $a$ and using \eqref{conv-mass0}, we get, $\dd z$ a.e.,
$$\lim_{\N\ni k\to\infty} \Theta^{(k)}_z(G\circ\mathfrak{R}_c) = \Theta_{z,r}(G\circ\mathfrak{R}_c).$$
Since the left-hand side does not depend on $r$, we conclude that, for every $r,r'>0$, we must have 
$\Theta_{z,r}(G\circ\mathfrak{R}_c)=\Theta_{z,r'}(G\circ\mathfrak{R}_c)$, $\dd z$ a.e., and since this holds
for any $c>0$ and any function $G$, we conclude that $\Theta_{z,r}=\Theta_{z,r'}$, $\dd z$ a.e. So, if we take
$\bar\Theta_z=\Theta_{z,1}$, the collection $(\bar\Theta_z)_{z>0}$ satisfies the first part of the statement. 

It remains to obtain
the scaling invariance property. To this end, we first observe that the process
$$X^{\{\lambda\}}_t:= \sqrt{\lambda} \,X_{t/\lambda}$$
remains a nine-dimensional Bessel process started from $0$. 
Furthermore, with an obvious notation, we have  $\mathrm{L}^{\{\lambda\}}_{r\sqrt{\lambda}}=\lambda \mathrm{L}_{r}$ for every $r>0$. Then,
it is straightforward to verify that the image of $\ll$ under the transformation
\begin{equation}
\label{scaling-trans}
\sum_{i\in I} \delta_{(t_i,\omega_i)} \mapsto \sum_{i\in I} \delta_{(\lambda t_i,\theta_\lambda(\omega_i))}
\end{equation}
is, conditionally on $X^{\{\lambda\}}$, a Poisson point measure with intensity
$$2\,\dd t\,\N_{X^{\{\lambda\}}_t}(\dd \omega \cap\{W_*>0\}).$$
It follows that
the image of $\ll^{(r)}$ under the scaling transformation \eqref{scaling-trans} has the same distribution as $\ll^{(r\sqrt{\lambda})}$.

We also note that, for every $x>0$, we have $\z_0(\theta_\lambda(\omega))=\lambda\,\z_0(\omega)$, $\N_x(\dd \omega)$ a.e. By combining
the preceding observations, we get that, for every $r>0$, the image of the triple
$(X^{(r)},\wt\ll^{(r)},\wt\rr^{(r)})$
under the scaling transformation $\Sigma_\lambda$ has the same distribution as 
$(X^{(r\sqrt{\lambda})},\wt\ll^{(r\sqrt{\lambda})},\wt\rr^{(r\sqrt{\lambda})})$,
and moreover the exit measure at $0$ associated with $\Sigma_\lambda(X^{(r)},\wt\ll^{(r)},\wt\rr^{(r)})$ is $\lambda\z^{(r)}$.
By considering conditional distributions with respect to $\z^{(r)}$ and using the first part of the proof, we obtain that
$\Sigma_\lambda(\bar\Theta_z)=\bar\Theta_{\lambda z}$ for a.e. $z>0$. A Fubini type argument allows us to
single out $z_0>0$ such that the equality $\bar\Theta_{\lambda z_0}=\Sigma_\lambda(\bar\Theta_{z_0})$ holds for a.e. $\lambda>0$. 
We then define, for every $z>0$,
$$\Theta_z=\Sigma_{z/z_0}(\bar\Theta_{z_0}).$$
Clearly the collection $(\Theta_z)_{z>0}$ is also 
a regular version of the 
conditional distributions of the triple
$(X^{(r)},\wt\ll^{(r)},\wt\rr^{(r)})$
knowing that $\z^{(r)}=z$ (for any $r>0$). Furthermore, by construction, the equality
$\Sigma_\lambda(\Theta_z)=\Theta_{\lambda z}$ holds 
for every $z>0$ and $\lambda>0$. This completes the proof, except for the uniqueness
statement, which is easy and
left to the reader. \endproof

From now on, $(\Theta_z)_{z>0}$ is the unique collection satisfying the properties stated in Proposition \ref{condi-dis2}. 
Thanks to the scaling invariance property, we can in fact define this collection without appealing to any conditioning.
We consider the triple $(X^{(r)},\wt\ll^{(r)},\wt\rr^{(r)})$ as defined at the beginning of the section, and recall
the notation $\z^{(r)}$ and the scaling operators $\Sigma_\lambda$ in Proposition \ref{condi-dis2}.

\begin{proposition}
\label{scaling-def}
Let $r>0$ and $z>0$. Then $\Theta_z$ is the distribution of $\Sigma_{z/\z^{(r)}}(X^{(r)},\wt\ll^{(r)},\wt\rr^{(r)})$.
\end{proposition}

\proof Let $F$ be a nonnegative measurable function on $\W^\infty\times M_p(\R_+\times \S)^2$, and recall that
the distribution of $\z^{(r)}$ has density $k_r(z)$. Then, using Proposition \ref{condi-dis2},
$$\E\Big[ F\Big(\Sigma_{z/\z^{(r)}}(X^{(r)},\wt\ll^{(r)},\wt\rr^{(r)})\Big)\Big]=\int_0^\infty \dd y\,k_r(y)\,\Theta_y(F\circ \Sigma_{z/y})= \Theta_z(F),$$
since the image of $\Theta_y$ under $\Sigma_{z/y}$ is $\Theta_z$. \endproof

Proposition \ref{scaling-def} is useful to derive almost sure properties of coding triples distributed 
according to $\Theta_z$. We give an important example.

\begin{corollary}
\label{transi-labels}
Let $z>0$, and let $\t$ be the labeled tree associated with a coding triple distributed according 
to $\Theta_z$. Write $(\Lambda_v)_{v\in\t}$ for the labels on $\t$ and $(\ee_s)_{s\in\R}$
for the clockwise exploration of $\t$. Then,
$$\lim_{|s|\to\infty}\Lambda_{\ee_s}=\infty\,,\quad\hbox{a.s.}$$
\end{corollary}

\proof By Proposition \ref{scaling-def}, it suffices to prove the similar statement for 
the labeled tree associated with $(X^{(r)},\wt\ll^{(r)},\wt\rr^{(r)})$, or even for the 
labeled tree associated with $(X,\ll,\rr)$. In the latter case this follows from \cite[Lemma 3.3]{CLG}. \endproof

\subsection{The coding triple of the Brownian disk with a given height}
\label{cod-height}

The fact that the collection $(\Theta_z)_{z>0}$ has been uniquely defined will now allow us to
make a canonical choice for the conditional distributions $(\Theta^{(a)}_z)_{z>0}$ (until now, these
conditional distributions were only defined up to a set of values of $z$ of zero Lesgue measure). This will be
important later as we use $\Theta^{(a)}_z$ to construct the free pointed Brownian disk with perimeter $z$ 
and height $a$.
% (as already mentioned in the introduction, the height is defined as the distance from the distinguished point to the boundary). 

Recall the restriction operator $\mathfrak{R}_a$ introduced in \eqref{restric-opera}, and the notation
$\lambda^{(a)}(\w)=\sup\{t\geq 0:\w(t)\leq a\}$ for $\w\in\mathcal{W}^\infty_0$ or $\w\in\W_0$.

\begin{proposition}
\label{finite-height}
Let $a>0$, and define a function $W_{*,a}:\W^\infty_0\times M_p(\R_+\times\S)^2\la \R_+\cup\{\infty\}$ by
$$W_{*,(a)}\Big(\w, \sum_{i\in I} \delta_{(t_i,\omega_i)}, \sum_{j\in J}\delta_{(t'_j,\omega'_j)}\Big)
= \min\Big(\inf_{i\in I, t_i>\lambda^{(a)}(\w)} W_*(\omega_i), \inf_{j\in J,t'_j>\lambda^{(a)}(\w)} W_*(\omega'_j)\Big).$$
Then, we have $\Theta_z(W_{*,(a)}>0)= \sqrt{\pi}\,2^{1/2}3^{-3/2}\,a^3\,z^{1/2}\,h_a(z)$ and 
$\Theta_z(W_{*,(a)}>0)\la 1$ as $a\to\infty$.
Furthermore, we can choose the collection $(\Theta^{(a)}_z)_{z>0}$
so that, for every $z>0$, 
$\Theta^{(a)}_z$ is the pushforward  of $\Theta_z(\cdot\mid W_{*,(a)}>0)$ under $\mathfrak{R}_a$.
\end{proposition}

%\rem The proof will show that $\Theta_z(W_{*,(a)}>0)>0$ for every $z>0$ and $a>0$, so that the definition
%of $\Theta_z(\cdot\mid W_{*,(a)}>0)$ makes sense. 

\proof 
Let $r>0$ and $a>0$. For test functions $f$ and $F$ defined on $\R_+$ and on $\W\times M_p(\R_+\times\W)^2$ respectively, we have from Proposition \ref{condi-dis2},
\begin{align*}
&\int \dd z\,k_r(z)\,f(z)\, \Theta_z\Big(F\circ\mathfrak{R}_a\,\mathbf{1}_{\{W_{*,(a)}>0\}}\Big)\\
&\quad=\E\Big[f(\z^{(r)})\,F\Big((X^{(r)}_t)_{0\leq t\leq \mathrm{L}^{(r)}_a}, \wt\ll^{(r,r+a)}, \wt\rr^{(r,r+a)}\Big)\,\mathbf{1}_{\{\z^{(r+a,\infty)}=0\}}\Big]\\
&\quad=\E\Big[f(\z^{(r,r+a)})\,F\Big((X^{(r)}_t)_{0\leq t\leq \mathrm{L}^{(r)}_a}, \wt\ll^{(r,r+a)}, \wt\rr^{(r,r+a)}\Big)\Big] \times \P(\z^{(r+a,\infty)}=0)\\
&\quad= \int \dd z\,g_{r,a}(z)\,f(z)\,\Theta^{(a)}_z(F)\,\times \P(\z^{(r+a,\infty)}=0),
\end{align*}
using Lemma \ref{lem-tech1} in the last equality.
It follows that we have $\dd z$ a.e.,
$$k_r(z)\,\Theta_z\Big(F\circ\mathfrak{R}_a\,\mathbf{1}_{\{W_{*,(a)}>0\}}\Big) = 
g_{r,a}(z)\, \P(\z^{(r+a,\infty)}=0)\,\,\Theta^{(a)}_z(F).$$
For $F=1$, we get that the equality $k_r(z)\,\Theta_z(W_{*,(a)}>0)= g_{r,a}(z)\,\P(\z^{(r+a,\infty)}=0)=g_{r,a}(z)\,(\frac{a}{r+a})^3$ holds $\dd z$ a.e., but then, by a scaling argument using also the monotonicity of
$\Theta_z(W_{*,(a)}>0)$ in the variable $a$, it must hold for every $z>0$ and $a>0$. It follows that  $\Theta_z(W_{*,(a)}>0)=(k_r(z))^{-1}g_{r,a}(z)\,(\frac{a}{r+a})^3  $, and the explicit 
formulas for $k_r(z)$ and $g_{r,a}(z)$ give the first assertion of the proposition. Furthermore, the previous display gives
$$\Theta^{(a)}_z(F)= \frac{\Theta_z\Big(F\circ\mathfrak{R}_a\,\mathbf{1}_{\{W_{*,(a)}>0\}}\Big)}{\Theta_z(W_{*,(a)}>0)},$$
$\dd z$ a.e. The second assertion follows. \endproof

In what follows, we assume that, for every $a>0$, the collection $(\Theta^{(a)}_z)_{z>0}$ is chosen as in the preceding proposition, and that the collection
 $(\check\Theta^{(a)}_z)_{z>0}$ is then derived from $(\Theta^{(a)}_z)_{z>0}$ via the time-reversal operation.
 From the scaling properties of $(\Theta_z)_{z>0}$, one checks that, for every $\lambda>0$,
 the pushforward of $\Theta^{(a)}_z$ under the 
 scaling operator $\Sigma_\lambda$ is $\Theta_{\lambda z}^{(\sqrt{\lambda} a)}$. 
 
 The following corollary, which  relates the measures $\Theta^{(a)}_z$ when $a$ varies (and $z$ is fixed)
 is an immediate consequence of Proposition \ref{finite-height}. Before stating this corollary, we note
 that both $\mathfrak{R}_a$ and $W_{*,(a)}$ still make sense as mappings defined on
 $\W_0\times M_p(\R_+\times\S)^2$.
 
 \begin{corollary}
 \label{finite-height-cor}
 Let $0<a<a'$. Then we have
 $$\Theta^{(a')}_z (W_{*,(a)}>0)= \frac{a^3 h_a(z)}{a'^3h_{a'}(z)},$$
 and $\Theta^{(a)}_z$ is the pushforward of $\Theta^{(a')}_z (\cdot\mid W_{*,(a)}>0)$ under $\mathfrak{R}_a$.
  \end{corollary}
% Let
% For every 
% fixed $z>0$, the probability measure $\Theta^{(a)}_z$ also depends continuously on $a$ in the following
% sense. Given $\delta>0$, for every $a>a'\geq \delta$, we can couple a triple distributed according to 
% $\Theta^{(a)}_z$ with a triple distributed according to $\Theta^{(a')}_z$ in such a way that the second 
% one coincides with the image of the first one under the restriction operator $\mathfrak{R}_{a'}$,
% with high probability when $a-a'$ is small. This property is a straightforward consequence of the 
% definition of $\Theta^{(a)}_z$ in Proposition \ref{finite-height}.
% 
% \medskip
We now use the collection $(\check\Theta^{(a)}_z)_{z>0}$ to construct 
a regular version of the conditional 
distributions of $\tr_0(\omega)$ under $\N_a$ knowing $\z_0=z$, for every $a>0$ and $z>0$. This regular version 
is a priori unique up to sets of values of $z$ of zero Lebesgue measure, but for our purposes
it is important that the conditional distribution is defined for {\it every} $z>0$. 

We fix $a>0$ and consider a triple $(Y^{(z)},\wt\mm^{(z)},\wt\mm'^{(z)})$ distributed according to $\check\Theta^{(a)}_z$.
As explained at the end of Section \ref{coding-infinite} (finite spine case), we can use this triple 
to construct a snake trajectory, which belongs to $\S_a$ and is denoted by $\Omega(Y^{(z)},\wt\mm^{(z)},\wt\mm'^{(z)})$.
We write $\N_a^{(z)}$ for the distribution of the snake trajectory $\Omega(Y^{(z)},\wt\mm^{(z)},\wt\mm'^{(z)})$.
%From scaling properties and Corollary \ref{finite-height-cor}, one checks that $\N_a^{(z)}$ depends continuously on the pair $(a,z)$. 

\begin{proposition}
\label{snake-condi-exit}
The collection $(\N_a^{(z)})_{z>0}$ forms a regular version of the conditional 
distributions of $\tr_0(\omega)$ under $\N_a$ knowing that $\z_0=z$.
\end{proposition}

\proof Recall the notation introduced before Proposition
\ref{Palm}: Under the measure $\N_a(\dd\omega)$, we can consider, for every $s\in(0,\sigma(\omega))$, the point measure $\pp_{(s)}$ (resp. $\pp'_{(s)}$) 
that gives the snake trajectories
associated with the subtrees branching off the left side (resp. off the right side) of the ancestral line of the vertex 
$p_{(\omega)}(s)$ in the genealogical tree of $\omega$. 
Also use the notation $\wt\pp_{(s)}$ (resp. $\wt\pp'_{(s)}$) for the point measure
$\pp_{(s)}$ (resp. $\pp'_{(s)}$) ``truncated at level $0$''. This makes sense if
$s$ is such that $W_s(t)>0$ for $0\leq t<\zeta_s$, which is the case we will consider.
From Proposition \ref{Palm}, and using also \eqref{tech-exit1}, we have, for every nonnegative measurable functions $f$ and $F$ defined on $\R_+$ and on $\W\times M_p(\R_+\times \S)^2$ respectively,
%$$\N_a\Big(\int_0^\sigma \dd L^0_s\,F(W_s,\wt\pp_{(s)},\wt\pp'_{(s)})\Big) = \E\Big[F(Y,\wt\mm,\wt\mm')\Big],$$
%where $(Y,\wt\mm,\wt\mm')$ is as in Section \ref{sec-ident}. 
%From the approximation formula \eqref{formu-exit-bis}, we see that $\z_0$ (under the measure $\N_a(\dd\omega)\dd L^0_s$) and $\z^Y$ (under $\P$) are given by the same measurable function of the triples $(W_s,\wt\pp_{(s)},\wt\pp'_{(s)})$ and 
%$(Y,\wt\mm,\wt\mm')$ respectively. Hence, we have also, for every nonnegative measurable function $f$ on $\R_+$,
\begin{equation}
\label{tech-exit11}
\N_a\Big(\int_0^\sigma \dd L^0_s\,f(\z_0)\,F(W_s,\wt\pp_{(s)},\wt\pp'_{(s)})\Big) = \E\Big[f(\z^Y)\,F(Y,\wt\mm,\wt\mm')\Big],
\end{equation}
where $(Y,\wt\mm,\wt\mm')$ and $\z^Y$ are as in Section \ref{sec-ident}.
Notice that, $\dd L^0_s$ a.e., we have $\tr_0(\omega)=\Omega(W_s,\wt\pp_{(s)},\wt\pp'_{(s)})$. Hence the previous identity also gives, for 
every nonnegative measurable function $H$ on $\S$,
\begin{equation}
\label{tech-exit12}
\N_a( \z_0 f(\z_0)\,H(\tr_0(\omega)))= \E[f(\z^Y)\,H(\Omega(Y,\wt\mm,\wt\mm'))].
\end{equation}
Since the density of $\z^Y$ is $zh_a(z)$ and $\check\Theta^{(a)}_z$ is the conditional distribution 
of $(Y,\mm,\mm')$ given $\z^Y=z$, the right-hand side can 
be written as
$$\int \dd z\,zh_a(z)\,f(z)\,\check\Theta^{(a)}_z(H\circ\Omega)= \int \dd z\,zh_a(z)\,f(z)\,\N_a^{(z)}(H),$$
by the very definition of $\N_a^{(z)}$. 
The statement of the proposition follows.\endproof

From the scaling properties of the measures $\Theta^{(a)}_z$, we immediately get that,
for every $\lambda>0$,
 the pushforward of $\N^{(z)}_a$ under the 
 scaling transformation $\theta_\lambda$ is $\N^{(\lambda z)}_{a\sqrt{\lambda}}$.
 
In view of further applications, we also note that the definition of the exit local time at $0$
makes sense under $\N^{(z)}_a$. Precisely, one gets that, $\N^{(z)}_a(\dd \omega)$ a.e.,
the limit
\begin{equation}
\label{formu-exit-ter}
\wt L^0_t:=\lim_{\ve \da 0} \frac{1}{\ve^2} \int_0^t \dd s\,\mathbf{1}_{\{\wh W_s(\omega)<\ve\}}
\end{equation}
exists uniformly for $t\geq 0$, and $\wt L^0_\infty= \wt L^0_\sigma= z$.
If $\N^{(z)}_a$ is replaced by $\N_a$ (and $\wh W_s(\omega)$ by $\wh W_s(\tr_0(\omega))$) this is just
formula \eqref{formu-exit-bis} in Section \ref{sna-mea}. So \eqref{formu-exit-ter} is a conditional version
of \eqref{formu-exit-bis}, which must therefore hold $\N^{(z)}_a$ a.e., at least for a.e.~value of $z$. But then a scaling argument,
using also the way we have defined the conditional 
distributions $\N^{(z)}_a$ and Corollary \ref{finite-height-cor}, shows that \eqref{formu-exit-ter} indeed holds for every $z>0$. We omit the details.

\section{From coding triples to random metric spaces}
\label{coding-metric}

\subsection{The pseudo-metric functions associated with a coding triple}
\label{sec:pseudo}

Let $(\w,\pp,\pp')$ be a coding triple satisfying the assumptions of Section \ref{coding-infinite}
 in the infinite spine case,
and let $(\t,(\Lambda_v)_{v\in\t})$ be the associated labeled tree. We suppose here that
labels take nonnegative values, $\Lambda_v\geq 0$ for every $v\in\t$, and we set $\t^\circ:=\{v\in\t:\Lambda_v>0\}$
and $\partial\t=\t\backslash\t^\circ$. We assume that $\partial\t$ is not empty and that all points of 
$\partial\t$ are leaves (points whose removal does not
disconnect $\t$). In particular, $\t^\circ$ is dense in $\t$.
We denote the clockwise exploration of $\t$ by $(\ee_t)_{t\in\R}$, and we assume that
either $\Lambda_{\ee_t}\la \infty$ as $|t|\to\infty$, or the set $\{t\in\R:\Lambda_{\ee_t}=0\}$ intersects both intervals $[K,\infty)$ and $(-\infty,-K]$, for every $K>0$. This ensures that
$\inf_{w\in[u,v]}\Lambda_w$ is attained for every ``interval'' $[u,v]$ of $\t$.

We define, for every $u,v\in \t^\circ$,
\begin{equation}
\label{delta0}\Delta^\circ(u,v):=\left\{ \begin{array}{ll}
{\displaystyle \Lambda_u + \Lambda_v -2 \max\Big( \inf_{w\in[u,v]} \Lambda_w,\inf_{w\in[v,u]} \Lambda_w\Big) \quad}&\hbox{if } {\displaystyle\max\Big( \inf_{w\in[u,v]} \Lambda_w,\inf_{w\in[v,u]} \Lambda_w\Big)>0},\\
+\infty&\hbox{otherwise.}
\end{array}
\right.
\end{equation}
We then let $\Delta(u,v)$, $u,v\in \t^\circ$ be the maximal symmetric function on $\t^\circ\times\t^\circ$
that is bounded above by $\Delta^\circ$ and satisfies the triangle inequality:
\begin{equation}
\label{delta-infty}
\Delta(u,v) = \inf_{u_0=u,u_1,\ldots,u_p=v} \sum_{i=1}^p \Delta^{\circ}(u_{i-1},u_i)
\end{equation}
where the infimum is over all choices of the integer $p\geq 1$ and of the
finite sequence $u_0,u_1,\ldots,u_p$ in $\t$ such that $u_0=u$ and
$u_p=v$. Then $\Delta(u,v)<\infty$ for every $u,v\in \t^\circ$. Indeed, a compactness argument shows
that we can find finitely many points $u_0=u,u_1,\ldots,u_{p-1},u_p=v$ belonging to
the geodesic segment $\llbracket u,v\rrbracket$ of $\t$ and such that
$\Delta^\circ(u_{i-1},u_{i})<\infty$ for every $1\leq i\leq p$.

Furthermore, the mapping $(u,v)\mapsto \Delta(u,v)$ is continuous on $\t^\circ\times\t^\circ$
(observe that $\Delta^\circ(u_n,u)\la 0$ if $u_n\to u$ in $\t^\circ$, and use the triangle
inequality). We note the trivial bound $\Delta^\circ(u,v)\geq |\Lambda_u-\Lambda_v|$, which also implies
\begin{equation}
\label{easybound}
\Delta(u,v)\geq |\Lambda_u-\Lambda_v|.
\end{equation}

We will call $\Delta^\circ(u,v)$ and $\Delta(u,v)$ the {\it pseudo-metric functions associated with the triple} $(\w,\pp,\pp')$.
From now on, let us assume that the function $(u,v)\mapsto \Delta(u,v)$ has a continuous extension
to $\t\times\t$, which is therefore a pseudo-metric on $\t$. 
We will be interested in the resulting quotient metric space $\t/\approx$ where the equivalence 
relation $\approx$ is defined by saying that $u\approx v$ if and only if $\Delta(u,v)=0$. By abuse of 
notation, we will write $\t/\{\Delta=0\}$ instead of $\t/\approx$. We write $\Pi$ for the canonical projection from $\t$ onto $\t/\{\Delta=0\}$. We also write $\Lambda_x=\Lambda_u$ when $x\in \t/\{\Delta=0\}$ and $u\in\t$
are such that $x=\Pi(u)$ (this is unambiguous by \eqref{easybound}). 

If $x\in \t/\{\Delta=0\}$ is such that $\Lambda_x>0$, we can define 
a geodesic path starting from $x$ in the following way. We pick $u\in\t$ such that $\Pi(u)=x$ and then
$s\in\R$ such that $\ee_s=u$.
We then define $\gamma^{(s)}=(\gamma^{(s)}_r)_{0\leq r\leq \Lambda_x}$ by setting $\gamma^{(s)}_r=\Pi(\ee_{\eta^{(s)}_r})$, with
$$\eta^{(s)}_r:=\left\{\begin{array}{ll}
\inf\{t\geq s: \Lambda_{\ee_t}=\Lambda_x-r\}\quad&\hbox{if }\ \inf\{\Lambda_{\ee_t}:t\geq s\}\leq \Lambda_x-r,\\
\inf\{t\leq s: \Lambda_{\ee_t}=\Lambda_x-r\}\quad&\hbox{if }\ \inf\{\Lambda_{\ee_t}:t\geq s\}> \Lambda_x-r.
\end{array}
\right.$$
It is then a simple matter to verify that $\gamma^{(s)}$ is a geodesic path in $(\t/\{\Delta=0\},\Delta)$, which starts from $x$
and ends at a point belonging to $\Pi(\partial \t)$. On the other hand, the bound \eqref{easybound}
shows that $\Delta(x,y)\geq \Lambda_x$ if $y\in\Pi(\partial \t)$. It follows that $\Delta(x,\Pi(\partial \t))=\Lambda_x$
for every $x\in \t/\{\Delta=0\}$. The path $\gamma^{(s)}$ is called a simple geodesic (see e.g.~\cite[Section 2.6]{Uniqueness} for 
the analogous definition
in the Brownian map). 

We finally note that
$\t/\{\Delta=0\}$ is a length space, meaning that the distance between two points is equal
to the infimum of the lengths of paths connecting these two points. To get this property, just notice that,
if $u,v\in\t^\circ$ and $\Delta^\circ(u,v)<\infty$, then $\Delta^\circ(u,v)$ coincides 
with the length of a path from $\Pi(u)$ to $\Pi(v)$, that is obtained by concatenating two simple geodesics starting from $\Pi(u)$ and $\Pi(v)$ respectively, up to the 
time when they merge. More explicitly, if $\Delta^\circ(u,v)= \Lambda_u+\Lambda_v-\inf_{w\in[u,v]} \Lambda_w$,
and if the reals $s'$ and $s''$ are such that $\ee_{s'}=u$, $\ee_{s''}=v$ and $[u,v]=\{\ee_r:r\in[s',s'']\}$, then 
the concatenation of $(\Pi(\gamma^{(s')}_r),0\leq r\leq \Lambda_u-\inf_{w\in[u,v]} \Lambda_w)$ and 
$(\Pi(\gamma^{(s'')}_r),0\leq r\leq \Lambda_v-\inf_{w\in[u,v]} \Lambda_w)$ gives a continuous
path from $\Pi(u)$ to $\Pi(v)$ with length $\Delta^\circ(u,v)$, which furthermore is contained in $\Pi([u,v])$. 

\subsection{The Brownian plane}
\label{Br-plane}

As an illustration of the procedure described in the previous section, and in view of further developments,
we briefly recall the construction of the Brownian plane given in \cite{CLG}. We consider a (random) coding
triple $(X,\ll,\rr)$ distributed as in Section \ref{cod-disk}:
\begin{itemize}
\item[$\bullet$]
$X=(X_t)_{t\geq 0}$ is a nine-dimensional Bessel process started from $0$.
\item[$\bullet$]
Conditionally on $X$, $\ll$ and $\rr$ are independent Poisson point measures on $\R_+\times \S$
with intensity
$$2\,\dd t\,\N_{X_t}(\dd \omega\cap \{W_*>0\}).$$
\end{itemize}
%where we recall the notation $W_*(\omega)=\min\{\wh W_s(\omega):0\leq s\leq \sigma(\omega)\}$.

It is easy to verify that the assumptions of Section \ref{coding-infinite} hold a.s. for  $(X,\ll,\rr)$,
and thus
we can associate an infinite labeled tree $(\t^p_\infty, (\Lambda_v)_{v\in\t^p_\infty})$ with this coding triple.
The assumptions of the beginning of Section \ref{sec:pseudo} also hold
(notice that the condition $\lim_{|s|\to\infty}\Lambda_{\ee_s}=\infty$ holds by \cite[Lemma 3.3]{CLG}), and we introduce
the two pseudo-metric functions $\Delta^{p,\circ}(u,v)$ and  $\Delta^{p}(u,v)$ defined for $u,v\in \t^{p,\circ}_\infty:=\{v\in\t^p_\infty:\Lambda_v>0\}$ via formulas \eqref{delta0}
and \eqref{delta-infty}. In that case, since the root of $\t^{p}_\infty$ is the only point with zero label, 
it is easy to see that at least one of the two infima $\inf_{w\in[u,v]} \Lambda_w$ and $\inf_{w\in[v,u]} \Lambda_w$ 
is positive, for any $u,v\in\t^{p,\circ}_\infty$ . Furthermore, it is immediate to obtain that $\Delta^{p,\circ}(u,v)$
and $\Delta^{p}(u,v)$ can be extended continuously to $\t^p_\infty$ --- in fact in that case we can define 
$\Delta^{p,\circ}(u,v)$ for every $u,v\in\t^p_\infty$ by the quantity in the first line of \eqref{delta0}, and use formula
\eqref{delta-infty} to define $\Delta^{p}(u,v)$ for every $u,v\in\t^p_\infty$.
One can prove \cite[Section 3.2]{CLG} that, for any $u,v\in \t^p_\infty$,  $\Delta^p(u,v)=0$ if and only if $\Delta^{p,\circ}(u,v)=0$.

The Brownian plane 
$\BP_\infty$
is defined as the quotient space $\t^p_\infty/\{\Delta^p=0\}$ equipped with the distance induced 
by $\Delta^p$ (for which we keep the same notation $\Delta^p$) and with the volume measure which is the pushforward of the volume
measure on $\t^p_\infty$ under the canonical projection. We note that $\BP_\infty$ comes with a distinguished point $\rho$,
which is the image of the root of $\t^p_\infty$ under the canonical projection.
Furthermore, we have $\Delta^p(\rho,x)=\Lambda_x$ for every $x\in\BP_\infty$.
% (labels $\Lambda_x$ make sense
%for $x\in\BP_\infty$, since \eqref{easybound} shows that
%$\Delta^p(u,v)=0$ implies $\Lambda_u=\Lambda_v$). 

The Brownian plane is scale invariant in the following sense. If $E$ is a pointed measure metric space and $\lambda>0$, we write 
$\lambda\cdot E$ for the same space $E$ with the metric multiplied by the factor $\lambda$
and the volume measure multiplied by the factor $\lambda^4$ (and the same distinguished point). Then, for every
$\lambda>0$, $\lambda\cdot\BP_\infty$ has the same distribution as $\BP_\infty$.

\subsection{The pointed Brownian disk with given perimeter and height}
\label{sec:point-disk}

In this section, we explain how a free pointed Brownian disk with perimeter $z$ and height $a$
is constructed from the measure $\N^{(z)}_a$ defined in Section \ref{cod-height}. This is basically an adaptation of
\cite{Disks}, but we provide some details in view of further developments.

We start with a preliminary result. Recall the notation $h_a(z)$ and $p_z(a)$ in Proposition \ref{densityuniform}.

\begin{proposition}
\label{condi-re-root}
For any nonnegative measurable functions $G$ and $f$ defined respectively on $\S$ and on $\R_+$, for every
$z>0$, we have
$$z^{-2}\N^{*,z}\Big(\int_0^\sigma \dd t \,G(W^{[t]})\,f(\wh W_t)\Big)
= \int_0^\infty \dd a\,p_z(a)\,f(a)\,\N^{(z)}_a(G).$$
\end{proposition}

\proof We may assume that both $G$ and $f$ are bounded and 
continuous. Then the argument is very similar to the proof of Proposition \ref{densityuniform}\,(ii) (which we recover when $G=1$). 
Let $g$ be a nonnegative measurable function on $\R_+$. We 
use the re-rooting formula \eqref{reroot-for}, and then Proposition \ref{densityuniform}\,(i), to get
\begin{align*}
\N^{*}\Big(\int_0^\sigma \dd t \,G(W^{[t]})\,f(\wh W_t)g(\z_0^*)\Big)
&=2\int_{0}^\infty \dd a\,\N_a\Big(\z_0 \,G(\tr_0(\omega))\,f(a)\,g(\z_0)\Big)\\
&=2\int_{0}^\infty \dd a\,f(a)\int_0^\infty \dd z\,h_{a}(z)\,zg(z)\,\N_a(G(\tr_0(\omega))\mid \z_0=z)\\
&=2\int_0^\infty \dd z\,zg(z)\int_0^\infty \dd a\,h_a(z)\,f(a)\,\N^{(z)}_a(G).
\end{align*}
On the other hand, the left-hand side is also equal to
$$\sqrt{\frac{3}{2\pi}}\,\int_0^\infty \dd z\,z^{-5/2}g(z)\,\N^{*,z}\Big(\int_0^\sigma \dd t\,G(W^{[t]})\,f(\wh W_t)\Big).$$
Since this holds for any function $g$, we must have, $\dd z$ a.e.,
$$z^{-2}\N^{*,z}\Big(\int_0^\sigma \dd t \,G(W^{[t]})\,f(\wh W_t)\Big)= \sqrt{\frac{8\pi}{3}} z^{3/2}
\int_0^\infty \dd a \,h_a(z)\,f(a)\,\N^{(z)}_a(G).$$
This is the identity of the proposition, except that we get it only $\dd z$ a.e. However, a scaling argument shows that
both sides of the preceding display are continuous functions of $z$, which gives the desired result for
every $z>0$. \endproof

Recall from Section \ref{posi-mea} the definition of the probability measure $\ov{\N}^{*,z}(\dd\omega \dd t)$ on $\S\times \R_+$,
and, for $(\omega,t)\in \S\times\R_+$, write $U(\omega,t)=t$. We can then rewrite the identity
of Proposition \ref{condi-re-root} in the form
\begin{equation}
\label{condition-height}
\ov{\N}^{*,z}(f(\wh \omega_U)G(\omega^{[U]})) = \int_0^\infty \dd a\,p_z(a)\,f(a)\,\N^{(z)}_a(G).
\end{equation}
%and we know (Proposition \ref{densityuniform}) that $p_z(\cdot)$ is the density 
%of the distribution of $\wh W_U$ under $\ov{\N}^{*,z}$. In other words, we have obtained that
%$$\ov{\N}^{*,z}(G(W^{[U]})\mid \wh W_U=a)= \N^{(z)}_a(G).$$

Let us now come to Brownian disks. We write $\mathbb{K}$, resp.  $\mathbb{K}^\bullet$, for the space of all 
compact measure metric spaces, resp. pointed compact measure metric spaces, equipped 
with the Gromov-Hausdorff-Prokhorov topology. Theorem 1 in \cite{Disks} provides a measurable mapping 
$\Xi:\S\la \mathbb{K}$ such that the distribution of $\Xi(\omega)$ under $\N^{*,z}(\dd \omega)$
is the law of the free Brownian disk with perimeter $z$. Let us briefly recall the construction of this mapping, which is
essentially an adaptation of the procedure of Section \ref{sec:pseudo} to the finite spine case.
Under 
$\N^{*,z}(\dd \omega)$, each vertex $u$ of the the genealogical tree $\t_{(\omega)}$ receives a nonnegative label $\ell_u(\omega)$, and
we define
the  ``boundary'' $\partial\t_{(\omega)}:=\{u\in\t_{(\omega)}: \ell_u(\omega)=0\}$. 
We also set $\t_{(\omega)}^\circ:=\t_{(\omega)}\backslash\partial\t_{(\omega)}$, and for every
$u,v\in\t_{(\omega)}^\circ$, we define $\Delta^{d,\circ}(u,v)$ and $\Delta^d(u,v)$ by the exact analogs 
of formulas \eqref{delta0} and \eqref{delta-infty}, where $\t$ and $\t^\circ$ are replaced by $\t_{(\omega)}$
and $\t^\circ_{(\omega)}$ respectively,
and the labels $(\Lambda_u)_{u\in\t}$ are replaced by $(\ell_u)_{u\in\t_{(\omega)}}$
(recall the definition of intervals on $\t_{(\omega)}$ in Section \ref{sna-tra}). 

A key technical point (Proposition 31 in \cite{Disks}) is to verify that 
the mapping $(u,v)\mapsto\Delta^{d}(u,v)$ can be extended continuously (in a unique way) to $\t_{(\omega)}\times\t_{(\omega)}$, $\N^{*,z}(\dd \omega)$ a.s. We then define
$\Xi(\omega)$
as the quotient space $\t_{(\omega)}/\{\Delta^d=0\}$, which is equipped with the metric induced by $\Delta^d$ and with a volume measure
which is the pushforward of the volume measure on $\t_{(\omega)}$ under the canonical projection. 
In the next definition, we use the notation $\bp_\omega$ for the composition of the 
canonical projection $p_{(\omega)}$ from $[0,\sigma(\omega)]$ onto $\t_{(\omega)}$ with the canonical projection
from $\t_{(\omega)}$ onto $\Xi(\omega)$.

\begin{definition}
\label{def-disk}
The distribution of the free pointed Brownian disk with perimeter $z$ is the law
of the random measure metric space $(\Xi(\omega),\Delta^d)$ pointed at the point $\bp_\omega(t)$,
under
the probability
measure $\ov{\N}^{*,z}(\dd \omega\,\dd t)$.
\end{definition}

The consistency of this definition of the free pointed Brownian disk
with the one in \cite{BM} follows from the results in \cite{Disks} (to be specific, \cite[Theorem 1]{Disks} identifies
the law of $\Xi(\omega)$ under $\N^{*,z}$ as the distribution of the free Brownian disk of
perimeter $z$, and one can then use formula (42) in \cite{Disk-Bdry} to see that Definition \ref{def-disk}
is consistent with the definition of the free pointed Brownian disk in \cite{BM}). 
From \cite{Bet,BM} one knows that, $\N^{*,z}(\dd \omega)$ a.s., $\Xi(\omega)$ is homeomorphic to the unit disk. This makes it possible to
define the boundary $\partial \Xi(\omega)$ of $\Xi(\omega)$, and this boundary is identified in \cite{Disks} with the image of $\partial\t_{(\omega)}$ under the
canonical projection. Furthermore, a.s. for every $u\in\t_{(\omega)}$, the label $\ell_u(\omega)$ is equal to the 
distance in $\Xi(\omega)$ from (the equivalence class of) $u$ to
the boundary $\partial \Xi(\omega)$.

We note that the random measure metric space $\Xi(\omega)$ pointed at $\bp_\omega(t)$ is in fact a function of
the re-rooted snake trajectory $\omega^{[t]}$, 
because $\Xi(\omega)$ is canonically identified to $\Xi(\omega^{[t]})$, and $\bp_\omega(t)$ is mapped to $\bp_{\omega^{[t]}}(0)$ in this identification. To simplify notation, we 
may thus write $\Xi^\bullet(\omega^{[t]})$ for the random metric space $\Xi(\omega)$ pointed at $\bp_\omega(t)$, or equivalently for $\Xi(\omega^{[t]})$ pointed at
$\bp_{\omega^{[t]}}(0)$.

In the following developments, it will be convenient to write $\D^\bullet_z$ for a free pointed Brownian disk with 
perimeter $z$ and $\partial\D^\bullet_z$ for the boundary of $\D^\bullet_z$. With a slight abuse of notation, we keep the notation $\Delta^d$
for the distance on $\D^\bullet_z$.  By definition, the height $H_z$ of $\D^\bullet_z$ is
the
distance from the distinguished point to the boundary. From the preceding interpretation of $\Xi(\omega)$ pointed at $\bp_\omega(t)$ as a function of the re-rooted snake trajectory $\omega^{[t]}$, we get
that, for any nonnegative measurable function $\Phi$ on $\mathbb{K}^\bullet$, for any nonnegative measurable function $h$ on $\R_+$,
$$\E[\Phi(\D^\bullet_z)\,h(H_z)]=\ov{\N}^{*,z}\Big(\Phi(\Xi^\bullet(\omega^{[U]}))\,h(\wh\omega_U)\Big)$$
using the same notation as in \eqref{condition-height} and noting that $\ell_{p_{(\omega)}(U)}(\omega)=\wh\omega_U$
by definition. 
By \eqref{condition-height}, we can rewrite this as
\begin{equation}
\label{desinteg-height}
\E[\Phi(\D^\bullet_z)\,h(H_z)]=\int_0^\infty \dd a\,p_z(a)\,h(a)\,\N^{(z)}_a(\Phi\circ \Xi^\bullet).
\end{equation}
The next proposition readily follows from \eqref{desinteg-height}.
\begin{proposition}
\label{pointed-disk}
The height $H_z$ of $\D^\bullet_z$ is distributed according to the density $p_z(a)$. 
Furthermore, the conditional distribution of $\D^\bullet_z$ knowing that $H_z=a$ is the law of $\Xi^\bullet(\omega)$ under $\N^{(z)}_a$.
By definition, this is the distribution of the free pointed Brownian disk with perimeter $z$ and
height $a$.
\end{proposition}

At this point, we should note that the definition of 
$\Xi^\bullet(\omega)$ requires the continuous extension of the
mapping $(u,v)\mapsto \Delta^d(u,v)$ from $\t^\circ_{(\omega)}\times\t^\circ_{(\omega)}$
to $\t_{(\omega)}\times\t_{(\omega)}$. Proposition 31 in \cite{Disks} and Proposition \ref{snake-condi-exit} above give the existence
of this continuous extension $\N^{(z)}_a(\dd \omega)$ a.s. for a.a. $z>0$, for every fixed $a>0$, and then one can use
scaling arguments and Corollary \ref{finite-height-cor} to get the same result for {\it every} $z>0$ and $a>0$ (this is used in Proposition \ref{pointed-disk}).
Similar considerations allow us to deduce the following property from Proposition 32 (iii) in \cite{Disks}: $\N^{(z)}_a(\dd\omega)$ a.s.,
for every $u,v\in \partial\t_{(\omega)}$, we have $\Delta^d(u,v)=0$ if and only if $\ell_w(\omega)>0$ for every $w\in ]u,v[$, or for every $w\in]v,u[$.
Finally, from \cite[Proposition 30 (iv)]{Disks}, we also get that, $\N^{(z)}_a(\dd\omega)$ a.s. for every $u,v\in\t^\circ_{(\omega)}$, we have
$\Delta^d(u,v)=0$ if and only if $\Delta^{d,\circ}(u,v)=0$.

It will be useful to introduce the uniform measure on the boundary $\partial \D^\bullet_z$. % (see \cite[Section 10]{Disks}).
There exists a measure $\mu_z$ on $\partial \D^\bullet_z$ with total mass equal to $z$, such that, a.s. for any continuous function
$\varphi$ on $\D^\bullet_z$, we have
$$\langle \mu_z,\varphi\rangle =\lim_{\ve\to 0} \frac{1}{\ve^2}\int_{\D^\bullet_z} \mathrm{Vol}(\dd x)\,\mathbf{1}_{\{\Delta^d(x,
\partial \D^\bullet_z)<\ve\}}\,\varphi(x)$$
where $\mathrm{Vol}(\dd x)$ denotes the volume measure on $\partial \D^\bullet_z$ 
(see \cite[Corollary 37]{Disks}). The preceding approximation 
and the definition of $\mu_z$ are also valid for the free pointed Brownian disk with perimeter $z$ and
height $a$: In fact, if this Brownian disk is constructed as $\Xi^\bullet(\omega)$ under $\N^{(z)}_a$ (as in Proposition \ref{pointed-disk}), 
we define $\mu_z$ by setting
\begin{equation}
\label{def-uniform}
\langle \mu_z,\varphi\rangle = \int_0^\infty \dd \wt L^0_s \,\varphi(\bp_\omega(s)),
\end{equation}
where the exit local time $\wt L^0_s$ is defined under $\N^{(z)}_a$ as explained at the end
of Section \ref{cod-height}. Note that the approximation formula for $\mu_z$
then reduces to formula \eqref{formu-exit-ter}, thanks to the interpretation of
labels as distances to the boundary. 

For our purposes in the next sections, it will be important to consider a free Brownian disk
(with perimeter $z$ and height $a$) equipped with a distinguished point chosen uniformly
on the boundary. To this end, we proceed as follows. 
We start from a triple $(Y^{(z)},\wt\mm^{(z)},\wt\mm'^{(z)})$
distributed according to $\check\Theta^{(a)}_z$. 
As explained before Proposition \ref{snake-condi-exit},
the random snake trajectory $\Omega(Y^{(z)},\wt\mm^{(z)},\wt\mm'^{(z)})$ is then distributed according to 
$\N^{(z)}_a$, and (Proposition \ref{pointed-disk}) we obtain a free pointed Brownian disk 
$\D^{\bullet,a}_z$ with perimeter $z$ and
height $a$ by setting
\begin{equation}
\label{const-disk-point}
\D^{\bullet,a}_z:=\Xi^\bullet\Big(\Omega(Y^{(z)},\wt\mm^{(z)},\wt\mm'^{(z)})\Big).
\end{equation}
In this construction, 
$\D^{\bullet,a}_z$ comes with a distinguished vertex of its boundary, namely the one
corresponding to the top of the spine of the tree coded by $(Y^{(z)},\wt\mm^{(z)},\wt\mm'^{(z)})$. We denote this
special point by $\alpha$. 

In the next proposition, we verify that $\alpha$ is (in a certain sense) uniformly distributed
over $\partial\D^{\bullet,a}_z$. To give a precise statement of this property, it is convenient to introduce the doubly pointed 
measure metric space $\D^{\bullet\bullet,a}_z$ which is obtained by viewing $\alpha$ as a second
distinguished point of $\D^{\bullet,a}_z$. 

\begin{proposition}
\label{doubly-point}
Let $F$ be a nonnegative measurable function on the space of all doubly pointed compact
measure metric spaces, which is equipped with the Gromov-Hausdorff-Prokhorov
topology. Then,
$$\E[F(\D^{\bullet\bullet,a}_z)]=\frac{1}{z}\,\E\Big[ \int \mu_z(\dd x)\,F\Big([\D^{\bullet,a}_z,x]\Big)\Big],$$
where $\mu_z$ is the uniform measure on $\partial \D^{\bullet,a}_z$, and we use the notation $[\D^{\bullet,a}_z,x]$ for the doubly pointed space obtained by
equipping $\D^{\bullet,a}_z$ with the second distinguished point $x$.
\end{proposition}

\proof 
Let $(Y,\mm,\mm')$ be distributed as in Proposition \ref{Palm} (or in Section \ref{sec-ident}). 
As previously, let $\wt \mm$ and $\wt \mm'$ be obtained by truncating the atoms of 
$\mm$ and $\mm'$ at level $0$. Let us introduce the notation $\|\wt\mm\|$ for the 
sum of the quantities $\sigma(\omega_i)$ over all atoms $(t_i,\omega_i)$
of $\wt\mm$. Note that $\|\wt\mm\|$ 
is the time at which the clockwise exploration of the tree associated with the triple
$(Y,\wt\mm,\wt\mm')$ (which coincides with the genealogical tree of $\Omega(Y,\wt\mm,\wt\mm')$)
visits the top of the spine. We start by observing that, for any nonnegative measurable
function $G$ on $\S\times \R_+$, for any nonnegative measurable function $h$ on $\R_+$,
\begin{align*}
\E\Big[ G\Big(\Omega(Y,\wt\mm,\wt\mm'),\|\wt\mm\|\Big)\,h(\z^Y)\Big]
&= \N_a\Big(\int_0^\infty \dd L^0_r\, G\Big(\tr_0(\omega),
{\int_0^r \dd s\mathbf{1}_{\{\zeta_{(\omega_s)}\leq \tau_0(\omega_s)\}}}
\Big)\, h(\z_0)\Big)\\
&= \N_a\Big(\int_0^\infty \dd \wt L^0_r\, G(\tr_0(\omega),r)\, h(\z_0)\Big)
\end{align*}
where $\wt L^0_r$ is defined under $\N_a$ as in formula \eqref{formu-exit-bis}. 
The first equality follows from Proposition \ref{Palm} as in the derivation of \eqref{tech-exit11} and \eqref{tech-exit12} above,
and the second equality is just a time change formula.
By conditioning on $\z^Y=z$ in the left-hand side and on $\z_0=z$ in the right-hand side, we get
$$\E\Big[ G(\Omega(Y^{(z)},\wt\mm^{(z)},\wt\mm'^{(z)}),\|\wt\mm^{(z)}\|)\Big]
=\frac{1}{z}\, \N^{(z)}_a\Big(\int_0^\infty \dd \wt L^0_r\, G(\omega,r)\Big).$$
Now recall that $\D^{\bullet,a}_z=\Xi^\bullet(\Omega(Y^{(z)},\wt\mm^{(z)},\wt\mm'^{(z)}))$
and that $\D^{\bullet\bullet,a}_z$ is obtained by assigning to $\D^{\bullet,a}_z$ a second
distinguished point equal to $\alpha=\bp_\Omega(\|\wt\mm^{(z)}\|)$, where
$\bp_\Omega$ denotes the composition of the canonical projection onto the genealogical
tree of $\Omega(Y,\wt\mm,\wt\mm')$ with the projection from this genealogical tree
onto $\Xi^\bullet(\Omega(Y^{(z)},\wt\mm^{(z)},\wt\mm'^{(z)}))$. Thanks to these observations,
we obtain that
$$\E[F(\D^{\bullet\bullet,a}_z)]=\frac{1}{z}\,\N^{(z)}_a\Big(\int_0^\infty \dd \wt L^0_r\, F\Big([\Xi^\bullet(\omega),
\bp_\omega(r)]\Big)\Big),$$
and the desired result follows since 
\eqref{def-uniform} shows that $\mu_z$ is the pushforward of the measure 
$\dd \wt L^0_r$ under the mapping $r\mapsto \bp_\omega(r)$.
\endproof

In view of forthcoming limit results where the distinguished point of $\D^{\bullet,a}_z$ 
is ``sent to infinity'', it will be convenient to introduce the random pointed measure
metric space $\bar\D^{\bullet,a}_z$ defined from the doubly pointed space $\D^{\bullet\bullet,a}_z$
by forgetting the first distinguished point. So $\bar\D^{\bullet,a}_z$ is pointed at 
a point which is uniformly distributed over its  boundary.

\subsection{Infinite-volume Brownian disks}
\label{inf-Bdisk}

For every $z>0$ and $a>0$, we keep the notation $\D^{\bullet,a}_z$ for 
a pointed Brownian disk with perimeter $z$ and
height $a$.
We may and will assume that $\D^{\bullet,a}_z$ is constructed from a
coding triple distributed according to $\check\Theta^{(a)}_z$ as in formula \eqref{const-disk-point}, or equivalently,
by time-reversal, from a coding triple distributed according to $\Theta^{(a)}_z$. The idea is now to let $a\to\infty$
and to use the ``convergence'' of $\Theta^{(a)}_z$ to $\Theta_z$ in order to get the convergence
of $\D^{\bullet,a}_z$ as $a\to\infty$. As we already mentioned, for a precise statement of this convergence, it will be
more convenient to replace $\D^{\bullet,a}_z$ by $\bar\D^{\bullet,a}_z$. The limit, which will be denoted by $\D^\infty_z$, is a 
random pointed locally compact measure metric space, which we call the infinite-volume Brownian disk
with perimeter $z$.

Let us start by explaining the construction of $\D^\infty_z$, which follows the general pattern of Section \ref{sec:pseudo}. We consider
a coding triple $((\rho_t)_{t\geq 0},\qq,\qq')$ 
with distribution $\Theta_z$. From $((\rho_t)_{t\geq 0},\qq,\qq')$, Section \ref{coding-infinite} allows us to construct an infinite tree $\t^i_\infty$ equipped with nonnegative labels $(\Lambda_v)_{v\in\t^i_\infty}$, 
such that labels on the spine are given by the process $(\rho_t)_{t\geq 0}$. Note that the assumptions in 
Section \ref{sec:pseudo} hold in
particular thanks to Corollary \ref{transi-labels}. 

We set $\t^{i,\circ}_\infty:=\{v\in \t^i_\infty: \Lambda_v>0\}$ and $\partial  \t^i_\infty:=\t^i_\infty\backslash
\t^{i,\circ}_\infty=\{v\in \t^i_\infty: \Lambda_v=0\}$. We define the 
pseudo-metric functions $\Delta^{i,\circ}(u,v)$ and $\Delta^i(u,v)$ on $\t^{i,\circ}_\infty\times \t^{i,\circ}_\infty$
as explained in Section \ref{sec:pseudo}. 

\begin{lemma}
\label{cons-inf-tech}
{\rm (i)} The mapping $(u,v)\mapsto \Delta^i(u,v)$ has a.s. a continuous extension to $\t^i_\infty\times \t^i_\infty$.

\noindent{\rm (ii)} A.s., for every $u,v\in \t^{i,\circ}_\infty$, the property $\Delta^i(u,v)=0$ holds if and only if 
$\Delta^{i,\circ}(u,v)=0$.

\noindent{\rm(iii)}  A.s., for every $u,v\in\partial \t^{i}_\infty$, the property $\Delta^i(u,v)=0$ holds if and only if
$\Lambda_w>0$ for every $w\in]u,v[$, or for every $w\in]v,u[$.
\end{lemma}

\proof Thanks to Proposition \ref{finite-height} and Corollary \ref{transi-labels}, it is enough to verify 
that properties analogous to (i),(ii),(iii) hold  when $\t^i_\infty$ is replaced by the labeled tree associated with a coding triple distributed according
to $\Theta^{(a)}_z$, for some fixed $a>0$. But then this is a 
consequence of the similar results in \cite{Disks}, as it was explained in the discussion after the statement of
Proposition \ref{pointed-disk}. \endproof

We let 
$\D^\infty_z$ denote the quotient space $\t^i_\infty/\{\Delta^i=0\}$, which is equipped
with the metric induced by $\Delta^i$. The volume measure on $\D^\infty_z$ 
is (as usual) the pushforward of the volume measure on $\t^i_\infty$. We also distinguish a special point $\alpha_\infty$
of $\D^\infty_z$, which is the equivalence class of the root of $\t^i_\infty$. 

\begin{definition}
\label{def-infinite-disk}
The random pointed measure metric space $(\D^\infty_z,\Delta^i)$ is the {\it infinite-volume Brownian disk with perimeter $z$}. 
\end{definition}

As in Section \ref{sec:pseudo}, labels $\Lambda_x$ make sense for 
$x\in \D^\infty_z$, and $\Lambda_x$ is equal to the distance from $x$
to the ``boundary'' $\partial\D^\infty_z$, which  is defined as the set of all points of
$\D^\infty_z$ with zero label. 
From the scaling properties of the collection $(\Theta_z)_{z>0}$, one also gets that
$\lambda\cdot \D^\infty_z$ is distributed as $\D^\infty_{\lambda^2z}$, for every $\lambda>0$.

In Section \ref{sec:consist} below, we will verify that this definition of the infinite-volume Brownian disk
is consistent with \cite{BMR}. It then follows from \cite[Corollary 3.13]{BMR} that 
$\D^\infty_z$ is homeomorphic to the complement of the open unit disk of the plane, so that the boundary $\partial\D^\infty_z$ can
be understood in a topological sense. We will not use this result, which may also be derived from our interpretation of complements of
hulls of the Brownian plane as infinite-volume Brownian disks in Section \ref{iBdBp}. 

\smallskip
\rem We could have defined the infinite-volume Brownian disk without distinguishing a special point of the boundary. 
The reason for distinguishing $\alpha_\infty$ comes from the use of the local Gromov-Hausdorff convergence 
in Theorem \ref{local-G-H} below, which requires dealing with pointed spaces. 

\smallskip

%We write $\dd_{\D^{\infty}_z}$ for the metric on $\D^{\infty}_z$, and similarly $\dd_{\D^{\bullet,a}_z}$ for the metric
%on $\D^{\bullet,a}_z$. 
For the convergence result to follow, it is convenient to deal
with the following definition of ``balls'': for every $h>0$,
$$\mathcal{B}_h(\D^{\infty}_z)=\{v\in \D^{\infty}_z : \Delta^i(v,\partial \D^{\infty}_z)\leq h\},$$
and
$$\mathcal{B}_h(\bar\D^{\bullet,a}_z)=\{v\in \bar\D^{\bullet,a}_z : \Delta^{(a)}(v,\partial \bar\D^{\bullet,a}_z)\leq h\},$$
where we use the notation $\Delta^{(a)}$ for the metric on $\D^{\bullet,a}_z$.
We view both $\mathcal{B}_h(\D^{\infty}_z)$ and $\mathcal{B}_h(\bar\D^{\bullet,a}_z)$ as compact measure metric spaces, which are pointed
at $\alpha_\infty$ and $\alpha$ respectively.
The compactness of $\mathcal{B}_h(\D^{\infty}_z)$ is a consequence of the fact that the set $\{u\in\t^i_\infty:\Lambda_u\leq h\}$
is compact, by Corollary \ref{transi-labels}. 

\begin{proposition}
\label{cons-infi}
Let $z>0$ and $h>0$. There exists a function $(\ve(a),a>0)$ with $\ve(a)\la 0$ as $a\to\infty$, such that,
for every $a>0$,
we can define on the same probability space both the infinite-volume Brownian disk $\D^\infty_z$
and the pointed Brownian disk $\bar\D^{\bullet,a}_z$, in such a way that
$$\P(\mathcal{B}_{h'}(\bar\D^{\bullet,a}_z)=\mathcal{B}_{h'}(\D^\infty_z),\hbox{ for every }0\leq h'\leq h) \geq 1 - \ve(a).$$
\end{proposition}

In other words, when $a$ large, one can couple the spaces $\bar\D^{\bullet,a}_z$ and $\D^\infty_z$
so that their tubular neighborhoods of the boundary (of any fixed radius $h$) are isometric except
on a set of small probability.

\proof
Let $a>0$, and consider a triple $((\rho^{(a)}_t)_{0\leq t\leq \mathrm{L}^{(a)}},\qq^{(a)},\qq'^{(a)})$ with distribution $\Theta^{(a)}_z$. As it was explained before Proposition \ref{doubly-point}, this triple allows us to construct
the Brownian disk $\bar\D^{\bullet,a}_z$ of perimeter $z$ and height $a$ pointed at a boundary point.
To be specific, the triple $((\rho^{(a)}_t)_{0\leq t\leq \mathrm{L}^{(a)}},\qq^{(a)},\qq'^{(a)})$ codes 
a random compact
tree $\t^{(a)}$ equipped with labels $(\Lambda^{(a)}_v)_{v\in \t^{(a)}}$. The pseudo-metric functions $\Delta^{(a),\circ}(u,v)$
and $\Delta^{(a)}(u,v)$ are then defined as in Section \ref{sec:pseudo} for $u,v\in \t^{(a)}$ such that $\Lambda^{(a)}_u>0$ and $\Lambda^{(a)}_v>0$. The function
$(u,v)\mapsto\Delta^{(a)}$ is extended by continuity to $\t^{(a)}\times\t^{(a)}$, and the resulting quotient metric space 
pointed at the root of $\t^{(a)}$
is the pointed Brownian disk $\bar\D^{\bullet,a}_z$ --- here we observe that the distinguished point $\alpha$
corresponds to the root and not to the top of the spine of $\t^{(a)}$, because the effect of
dealing with the triple $((\rho^{(a)}_t)_{0\leq t\leq \mathrm{L}^{(a)}},\qq^{(a)},\qq'^{(a)})$ instead of its image under the time-reversal transformation \eqref{time-rev}
interchanges the roles of the root and the top of the spine (see the comments after \eqref{time-rev}).

Recall the restriction operator $\mathfrak{R}_a$ in \eqref{restric-opera}.
%, and set
%$$\Big((\rho_t)_{0\leq t\leq \lambda^{(a)}_\infty},\qq_\infty^{(a)},\qq'^{(a)}_\infty\Big)=\mathfrak{R}_a\Big(
%(\rho^{(a)}_t)_{0\leq t\leq \mathrm{L}^{(a)}},\qq^{(a)},\qq'^{(a)}\Big),$$
%so that 
%$\lambda_\infty^{(a)}=\sup\{t\geq 0: \rho_t=a\}$ and $\qq_\infty^{(a)}$, resp. $\qq'^{(a)}_\infty$, denotes the restriction
%of $\qq$, resp. of $\qq'$, to $[0,\lambda_\infty^{(a)}]\times \S$.
%and write $\t_\infty^{i,(a)}$ for the compact tree derived from $\t^i_\infty$ by removing the part of the spine above height $\lambda_\infty^{(a)}$
%(and of course the subtrees branching off this part of the spine). Also let $\qq_\infty^{(a)}$, resp. $\qq'^{(a)}_\infty$ denote the restriction
%of $\qq$, resp. of $\qq'$, to $[0,\lambda_\infty^{(a)}]\times \S$. 

\begin{lemma}
\label{couple-tree}
We can couple the triple $((\rho^{(a)}_t)_{0\leq t\leq \mathrm{L}^{(a)}},\qq^{(a)},\qq'^{(a)})$ distributed according to $\Theta^{(a)}_z$
 and the triple $((\rho_t)_{t\geq 0},\qq,\qq')$ distributed according to $\Theta_z$ so that the property
%$$ \P\Big(\Big(\t_\infty^{i,(a)},(\Lambda_v)_{v\in\t_\infty^{i,(a)}}\Big)= \Big(\t^{(a)},(\Lambda^{(a)}_v)_{v\in\t^{(a)}}\Big)\Big) \geq1-\ve_a,$$
%where $\ve_a\la 0$ as $a\to\infty$.
\begin{equation}
\label{couple-tech}
\Big((\rho^{(a)}_t)_{0\leq t\leq \mathrm{L}^{(a)}},\qq^{(a)},\qq'^{(a)}\Big)
=\mathfrak{R}_a\Big((\rho_t)_{t\geq 0},\qq,\qq'\Big)
\end{equation}
holds with probability tending to $1$ as $a\to\infty$. 
\end{lemma}

%\Big((\rho_t)_{0\leq t\leq \lambda^{(a)}_\infty},\qq_\infty^{(a)},\qq'^{(a)}_\infty\Big)=
\proof 
It suffices to verify that the variation distance between $\Theta^{(a)}_z$ and $\mathfrak{R}_a(\Theta_z)$ tends 
to $0$ as $a\to\infty$. This is an easy consequence of Proposition \ref{finite-height}. Indeed, let
$A$ be a measurable subset of $\W\times M_p(\R_+\times\S)^2$. Then,
\begin{align*}
\Theta_z(\mathfrak{R}_a^{-1}(A))&= \Theta_z(\mathfrak{R}_a^{-1}(A)\cap\{W_{*,(a)}>0\}) + 
 \Theta_z(\mathfrak{R}_a^{-1}(A)\cap\{W_{*,(a)}=0\})\\
 &=\Theta_z(W_{*,(a)}>0)\,\Theta_z^{(a)}(A) + \Theta_z(\mathfrak{R}_a^{-1}(A)\cap\{W_{*,(a)}=0\}),
 \end{align*}
 by Proposition \ref{finite-height}.  It follows that the variation distance between $\Theta^{(a)}_z$ and $\mathfrak{R}_a(\Theta_z)$ is bounded above by $1-\Theta_z(W_{*,(a)}>0)$, which tends to $0$
 as $a\to\infty$.%, by the explicit formula for $\Theta_z(W_{*,(a)}>0)$ found in Proposition \ref{finite-height}.
%Let $G$ be a measurable function on $\W\times M_p(\R_+\times\S)^2$ such that $0\leq G\leq 1$. 
%Fix $u>0$. It follows from formula \eqref{tech-import} (applied with $c=a$ and $r$ replaced by $u$) that
%\begin{equation}
%\label{couple-tech}
%\E\Big[G\Big((\rho_t)_{0\leq t\leq \lambda^{(a)}_\infty},\qq_\infty^{(a)},\qq'^{(a)}_\infty\Big)\Big]
%=\nu_{u,a}(\{0\}\mid z)\,\E\Big[G\Big((\rho^{(a)}_t)_{0\leq t\leq \mathrm{L}^{(a)}},\qq^{(a)},\qq'^{(a)}\Big)\Big] + \kappa_{u,a}(z),
%\end{equation}
%where the remainder $\kappa_{u,a}(z)$
%is nonnegative and bounded above by $1-\nu_{u,a}(\{0\}\mid z)$. It follows that we can couple the two triples
%$((\rho_t)_{t\geq 0},\qq,\qq')$ and $((\rho^{(a)}_t)_{0\leq t\leq \mathrm{L}^{(a)}},\qq^{(a)},\qq'^{(a)})$ in such a way that the event
%$$\Big((\rho_t)_{0\leq t\leq \lambda^{(a)}_\infty},\qq_\infty^{(a)},\qq'^{(a)}_\infty\Big)=\Big((\rho^{(a)}_t)_{0\leq t\leq \mathrm{L}^{(a)}},\qq^{(a)},\qq'^{(a)}\Big)$$
%occurs at least with probability $\nu_{u,a}(\{0\}\mid z)$. Now recall from \eqref{conv-mass0} that $\nu_{u,a}(\{0\}\mid z)\la 1$ as $a\to\infty$. The statement
%of the lemma follows since the labeled trees $(\t_\infty^{i,(a)},(\Lambda_v)_{v\in\t_\infty^{i,(a)}})$ 
%and $(\t^{(a)},(\Lambda^{(a)}_v)_{v\in\t^{(a)}})$ are given by the same function of the triples
%$((\rho_t)_{0\leq t\leq \lambda^{(a)}_\infty},\qq_\infty^{(a)},\qq'^{(a)}_\infty)$ and $((\rho^{(a)}_t)_{0\leq t\leq \mathrm{L}^{(a)}},\qq^{(a)},\qq'^{(a)})$
%respectively. 
\endproof

It will be convenient to write $\t_\infty^{i,(a)}$ for the (labeled) compact tree derived from $\t^i_\infty$ by removing the part of 
the spine above height $\lambda_\infty^{(a)}:=\sup\{t\geq 0:\rho_t=a\}$
(and of course the subtrees branching off this part of the spine). We can also view $\t_\infty^{i,(a)}$
as the labeled tree coded by $\mathfrak{R}_a((\rho_t)_{t\geq 0},\qq,\qq')$. On the event where
\eqref{couple-tech} holds, we can therefore also identify $\t_\infty^{i,(a)}$ with the labeled tree $\t^{(a)}$
coded by $((\rho^{(a)}_t)_{0\leq t\leq \mathrm{L}^{(a)}},\qq^{(a)},\qq'^{(a)})$, and this identification is used
in the next lemma and its proof. 

From now on, 
we assume that the triples $((\rho^{(a)}_t)_{0\leq t\leq \mathrm{L}^{(a)}},\qq^{(a)},\qq'^{(a)})$ and 
$((\rho_t)_{t\geq 0},\qq,\qq')$ are coupled as in Lemma \ref{couple-tree}, and that
$\D^{\bullet,a}_z$ and $\D^\infty_z$ are constructed from these triples as explained above. 

\begin{lemma}
\label{key-cond-infi}
Let $h>0$. Set 
$$A=\max\{\Delta^i(x,y):x,y\in \D^\infty_z, \Lambda_x\leq h,\Lambda_y\leq h\}.$$
On the intersection of the event where \eqref{couple-tech} holds 
%\begin{equation}
%\label{tree-equal}
%\Big(\t^{i,(a)}_\infty, (\Lambda_v)_{v\in \t_\infty^{i,(a)}}\Big)
%= \Big(\t^{(a)}, (\Lambda^{(a)}_v)_{v\in \t^{(a)}}\Big)
%\end{equation}
with the event where 
\begin{equation}
\label{label-large}
\inf\{ \Lambda_v: v\in \t^i_\infty \backslash\t^{i,(a)}_\infty\} \geq A+h+1
\end{equation}
we have 
$$\Delta^i(v,w)=\Delta^{(a)}(v,w)$$
for every $v,w\in\t^{(a)}$ such that $\Lambda_v\leq h$ and $\Lambda_w\leq h$.
\end{lemma}

\rem The statement of the lemma makes sense because on the event 
where \eqref{couple-tech} holds, the trees $\t^{(a)}$ and $\t^{i,(a)}_\infty$
are identified (as explained before the statement of the lemma), and so 
$\Delta^{(a)}(v,w)$ and $\Delta^i(v,w)$ both make sense when 
$v,w\in\t^{(a)}=\t^{i,(a)}_\infty$. 

\smallskip
The statement of the proposition follows from Lemma \ref{key-cond-infi}.
Indeed, by Corollary \ref{transi-labels},
$$\inf\{ \Lambda_u: u\in \t^i_\infty \backslash\t^{i,(a)}_\infty\}\build{\la}_{a\to\infty}^{} +\infty\;,\hbox{ a.s.}$$
and so, when $a$ is large, the property \eqref{label-large} will hold except
on a set of small probability. Also, by Lemma \ref{couple-tree}, we know that 
the property  \eqref{couple-tech} holds outside an event of small probability. 
Note that, when \eqref{label-large} holds, 
labels do not vanish on $\t^i_\infty \backslash\t^{i,(a)}_\infty$, and so the 
``boundary'' (set
of points with zero label) of $\t^i_\infty$ is identified with the boundary
of $\t^{(a)}$. 
When \eqref{couple-tech} and \eqref{label-large} both hold, the conclusion of
Lemma \ref{key-cond-infi} shows that $\mathcal{B}_{h'}(\D^{\infty}_z)$ and $\mathcal{B}_{h'}(\bar\D^{\bullet,a}_z)$ are isometric, for every $0\leq h'\leq h$. 

\proof[Proof of Lemma \ref{key-cond-infi}]
Throughout the proof we assume that both \eqref{couple-tech} and \eqref{label-large} hold, so that $\t^{(a)}$ and $\t^{i,(a)}_\infty$
are identified. If $u,v\in\t^{i,(a)}_\infty=\t^{(a)}$, we use the notation $[u,v]_{\t^i_\infty}$ for the interval from $u$ to $v$ in $\t^i_\infty$,
and similarly $[u,v]_{\t^{(a)}}$ for the same interval in $\t^{(a)}$. We note that either $[u,v]_{\t^i_\infty}\subset \t^{i,(a)}_\infty$, 
and then $[u,v]_{\t^i_\infty}=[u,v]_{\t^{(a)}}$, or $[u,v]_{\t^i_\infty}\not\subset \t^{i,(a)}_\infty$, and then $[u,v]_{\t^i_\infty}$
is the union of $[u,v]_{\t^{(a)}}$ and $\t^i_\infty\backslash\t^{i,(a)}_\infty$.

We use the notation $\t_\infty^{i,(a),\circ}=\{u\in\t_\infty^{i,(a)}: \Lambda_u>0\}$.
We first observe that, if $u,v\in \t_\infty^{i,(a),\circ}$, we have
\begin{equation}
\label{easybd}
\Delta^{(a),\circ}(u,v)\leq \Delta^{i,\circ}(u,v).
\end{equation}
Let us explain this bound. Since labels do not vanish on $\t^i_\infty\backslash\t_\infty^{i,(a)}$, it is immediate that 
$\Delta^{(a),\circ}(u,v)=\infty$ if and only if $\Delta^{i,\circ}(u,v)=\infty$. So we may assume that
both are finite and  then \eqref{easybd} directely follows from the definition of these quantities and the fact that we have always $[u,v]_{\t^{(a)}}\subset [u,v]_{\t^i_\infty}$
(and labels are the same on $\t^{(a)}$ and $\t^{i,(a)}_\infty$).
%\begin{align}
%\label{formula-Del}
%\Delta^{i,\circ}(u,v)&=\Lambda_u+\Lambda_v -2\max\Big(\inf_{w\in [u,v]_{\t^i_\infty}}\Lambda_w,
%\inf_{w\in [v,u]_{\t^i_\infty}}\Lambda_w\Big)\nonumber\\
%\Delta^{(a),\circ}(u,v)&=\Lambda_u+\Lambda_v -2\max\Big(\inf_{w\in [u,v]_{\t^{(a)}}}\Lambda_w,
%\inf_{w\in [v,u]_{\t^{(a)}}}\Lambda_w\Big)
%\end{align}
%with the obvious notation $[u,v]_{\t^i_\infty}$ for the interval from $u$ to $v$ in $\t^i_\infty$, and
%similarly for $[u,v]_{\t^{(a)}}$. Since $[u,v]_{\t^{(a)}}\subset [u,v]_{\t^i_\infty}$, \eqref{easybd}
%follows at once. 

Next suppose that $u,v\in \t_\infty^{i,(a),\circ}$ satisfy also $\Lambda_u\leq h$ and $\Lambda_v\leq h$. 
Recall the definition \eqref{delta-infty} of $\Delta^i(u,v)$ as an infimum involving all possible choices of 
$u_0,u_1,\ldots,u_p$ in $\t^{i,\circ}_\infty$ such that $u_0=u$ and $u_p=v$. We claim that in this definition
we can restrict our attention to the case when $u_0,u_1,\ldots,u_p$ satisfy 
$\Lambda_{u_j}< A+h+1$ for every $1\leq j\leq p$,
and therefore $u_0,u_1,\ldots,u_p\in\t^{i,(a)}_\infty$, by \eqref{label-large}.
To see this, suppose that $\Lambda_{u_k}\geq A+h+1$, for some 
$k\in\{1,\ldots,p-1\}$. Then, from the definition of $A$, we have also
$\Lambda_{u_k}\geq \Delta^i(u,v)+h+1$. 
Hence, using the bound \eqref{easybound}, we get
$$\sum_{j=1}^p \Delta^{i,\circ}(u_{j-1},u_j)\geq |\Lambda_{u_k}-\Lambda_u|\geq (\Delta^i(u,v)+h+1)-h=\Delta^i(u,v)+1,$$
so that we may disregard the sequence $u_0,u_1,\ldots,u_p$ in the infimum defining $\Delta^i(u,v)$.

By the previous considerations and \eqref{easybd}, we get, for $u,v\in \t_\infty^{i,(a),\circ}$ 
such that $\Lambda_u\leq h$ and $\Lambda_v\leq h$,
\begin{equation}
\label{cons-infi-1}
\Delta^i(u,v)=\build{\inf_{u_0=u,u_1,\ldots,u_p=v}}_{u_1,\ldots,u_{p-1}\in\t_\infty^{i,(a)}}^{}
\;\sum_{j=1}^p \Delta^{i,\circ}(u_{j-1},u_j)
\geq \build{\inf_{u_0=u,u_1,\ldots,u_p=v}}_{u_1,\ldots,u_{p-1}\in\t^{(a)}}^{}
\;\sum_{j=1}^p \Delta^{(a),\circ}(u_{j-1},u_j)= \Delta^{(a)}(u,v).
\end{equation}
We now want to argue that we have indeed the equality $\Delta^{(a)}(u,v)=\Delta^i(u,v)$. 
To this end it is enough to show that, for any sequence $u_0=u,u_1,\ldots,u_p=v$ in
$\t^{(a)}$ such that
\begin{equation}
\label{cons-infi-2}
\sum_{j=1}^p \Delta^{(a),\circ}(u_{j-1},u_j)<\Delta^{(a)}(u,v)+1,
\end{equation}
we have in fact $\Delta^{(a),\circ}(u_{j-1},u_j)=\Delta^{i,\circ}(u_{j-1},u_j)$ for every $j\in\{1,\ldots,p\}$
(this will entail that the two infima in \eqref{cons-infi-1} are equal). We argue by contradiction and suppose that
$\Delta^{(a),\circ}(u_{j-1},u_j)<\Delta^{i,\circ}(u_{j-1},u_j)$ for some $j\in\{1,\ldots,p\}$. This means that we have
$$\inf_{w\in[u_{j-1},u_j]_{\t^i_\infty}}\Lambda_w < \inf_{w\in [u_{j-1},u_j]_{\t^{(a)}}}\Lambda_w,$$
or the same with $[u_{j-1},u_j]$ replaced by $[u_j,u_{j-1}]$. However,
$[u_{j-1},u_j]_{\t^i_\infty}$ can be different from $[u_{j-1},u_j]_{\t^{(a)}}$ only if 
$[u_{j-1},u_j]_{\t^i_\infty}$ is the union of  $[u_{j-1},u_j]_{\t^{(a)}}$ and $\t^i_\infty\backslash \t_\infty^{i,(a)}$, and we get
$$\inf_{w\in [u_{j-1},u_j]_{\t^{(a)}}}\Lambda_w \geq \inf_{w\in\t^i_\infty\backslash \t_\infty^{i,(a)}}\Lambda_w
\geq A + h +1. $$
This implies in particular that $\Lambda_{u_j}\geq A+h+1$, and by the same argument as above this
gives a contradiction with \eqref{cons-infi-2}. This completes the proof of the lemma. \endproof

In the next statement, we use the local Gromov-Hausdorff-Prokhorov convergence  for pointed locally compact measure length spaces, as defined in \cite{ADH}.

\begin{theorem}
\label{local-G-H}
We have 
$$\bar\D^{\bullet,a}_z \build{\la}_{a\to\infty}^{\rm(d)} \D^\infty_z$$
in distribution in the sense of the local Gromov-Hausdorff-Prokhorov convergence. 
\end{theorem}

\proof 
The statement of the theorem is an immediate consequence of Proposition \ref{cons-infi}.
In fact, it is enough to verify that, for every $h>0$, the closed ball of
radius $h$ centered at the distinguished point $\alpha$ of $\bar\D^{\bullet,a}_z$ converges in
distribution  to the corresponding ball in $\D^\infty_z$ as $a\to\infty$, in the sense
of the Gromov-Hausdorff-Prokhorov convergence for pointed compact measure metric spaces.
However, this readily follows from the coupling obtained in Proposition \ref{cons-infi}, since
the closed ball of radius $h$ centered at $\alpha$ is obviously contained in $\mathcal{B}_h(\bar\D^{\bullet,a}_z)$
and similarly in the limiting space. \endproof

We conclude this section with a couple of almost sure properties of the infinite-volume
Brownian disk $\D^\infty_z$ that can be derived from our approach. First, from the analogous result
for the disk $\D^{\bullet,a}_z$ (see the remarks after Proposition \ref{pointed-disk})
and the coupling in Proposition \ref{cons-infi}, one easily obtains the existence 
of the uniform measure $\mu^\infty_z$ on $\partial \D^\infty_z$, which is a measure
of total mass $z$ satisfying
$$\langle \mu^\infty_z,\varphi\rangle =\lim_{\ve\to 0} \frac{1}{\ve^2}\int_{\D^\infty_z} \mathrm{Vol}(\dd x)\,\mathbf{1}_{\{\Delta^i(x,
\partial \D^\infty_z)<\ve\}}\,\varphi(x),$$
for any continuous function
$\varphi$ on $\D^\infty_z$, a.s. In particular the volume of the tubular neighborhhood of radius $\ve$ of $\partial \D^\infty_z$
behaves like $z\ve^2$ when $\ve\to 0$.

Our construction of $\D^\infty_z$ is well suited to the analysis of geodesics to
the boundary. Write $(\ee^i_s)_{s\in\R}$ for the clockwise exploration of the tree $\t^i_\infty$, and
set 
$$s_0:=\min\{s\in \R: \Lambda_{\ee^i_s}=0\}\,,\ s_1:=\max\{s\in \R: \Lambda_{\ee^i_s}=0\}.$$
Also set $x_0:=\Pi^i(\ee_{s_0})=\Pi^i(\ee_{s_1})$.

\begin{proposition}
\label{geode}
Almost surely, there exists $h_0>0$ such that, for 
every $x\in \D^\infty_z$ with $\Delta^i(x,\partial\D^\infty_z)> h_0$,
any geodesic from $x$ to $\partial\D^\infty_z$ hits $\partial\D^\infty_z$
at $x_0$.
\end{proposition}

The proof shows more precisely that all geodesics to $\partial\D^\infty_z$ starting outside a sufficiently large ball
coalesce before hitting the boundary. 

\proof Recall the notation $\t^{i,(a)}_\infty$ in the proof of
Proposition \ref{cons-infi}. By Corollary \ref{transi-labels}, we may choose $a$ large enough so that labels $\Lambda_v$
do not vanish on $\t^i_\infty\backslash \t^{i,(a)}_\infty$. Then we may take
$h_0=\max\{\Lambda_v:v\in \t^{i,(a)}_\infty\}$. To verify this, fix $x\in \D^\infty_z$ 
such that $\Delta^i(x,\partial\D^\infty_z)> h_0$, then we may write
$x=\Pi^i(v)$ with $v\in \t^i_\infty\backslash \t^{i,(a)}_\infty$, and we have
$\Delta^i(x,\partial\D^\infty_z)=\Lambda_v$. Consider a simple geodesic
$\gamma=(\gamma_r)_{0\leq r\leq \Lambda_v}$ from $x$ to $\partial\D^\infty_z$
constructed as in Section \ref{sec:pseudo}. Then it is straightforward to
verify that $\gamma_{\Lambda_v}=x_0$. To complete
the proof, we just need the fact that any geodesic from $x$ to $\partial\D^\infty_z$ is a simple geodesic.
This follows via Theorem \ref{infBdBp} below from the analogous result in the Brownian plane, which is itself a 
consequence of the study of geodesics in the Brownian map \cite{Acta}. We omit the details. \endproof

One may also consider geodesic rays in the infinite Brownian disk (a geodesic ray $\gamma=(\gamma_t)_{t\geq 0}$
is an infinite geodesic path). In a way analogous to the case of the Brownian plane (see \cite[Theorem 18]{Plane})
one obtains that any two geodesic rays in $\D^\infty_z$ coalesce in finite time. Again this can be deduced from the 
Brownian plane result via Theorem \ref{infBdBp}, but this also follows, with some more work, from the 
alternative construction of the infinite Brownian disk presented in Section \ref{sec:consist} below.

\subsection{The Brownian half-plane}
\label{half-plane}

In this section, we define the Brownian half-plane and show that it is the tangent cone 
in distribution of the free
pointed Brownian disk at a point chosen uniformly on its boundary. Let us start with the definition.
We consider a coding triple $(R,\pp,\pp')$, where $R=(R_t)_{t\in[0,\infty)}$ is a three-dimensional Bessel process started from $0$,  and,
conditionally on $R$, $\pp$ and $\pp'$ are independent Poisson point measures on $\R_+\times \S$ with intensity 
$2\,\dd t\,\N_{R_t}(\dd \omega)$. For every $r>0$, we set 
 $\mathbf{L}_r:=\sup\{t\geq 0: R_t=r\}$ (as in Section \ref{sec-ident}), and we let $\wt\pp$ and $\wt\pp'$ stand for the point measures $\pp$ and $\pp'$
truncated at level $0$. 

Following Section \ref{coding-infinite}, we can use the coding triple $(R,\wt\pp,\wt\pp')$
to construct a tree $\t^{hp}_\infty$ equipped with nonnegative labels $(\Lambda_v)_{v\in\t^{hp}_\infty}$.
In contrast with the measures $\Theta_z$ used to define the infinite-volume Brownian disk, 
there is no conditioning on the total exit measure at $0$, which is here infinite a.s., as it can be seen 
from a scaling argument. There are subtrees carrying zero labels that branch off the right side or the left side of the spine at arbitrary high levels, so that labels along the clockwise exploration of $\t^{hp}_\infty$
vanish in both intervals $(-\infty,-K]$ and $[K,\infty)$, for any $K>0$. 

We then follow the general procedure of Section \ref{sec:pseudo}. We set $\t_\infty^{hp,\circ}:=\{v\in \t^{hp}_\infty: \Lambda_v>0\}$ 
and $\partial \t^{hp}_\infty:=\t^{hp}_\infty\backslash \t^{hp,\circ}_\infty$, and we let $\Delta^{hp,\circ}(u,v)$ and $\Delta^{hp}(u,v)$, for $u,v\in \t_\infty^{hp,\circ}$, be the pseudo-metric functions associated with
the triple $(R,\wt\pp,\wt\pp')$ as in Section \ref{sec:pseudo}. 

\begin{lemma}
\label{cons-hp-tech}
{\rm (i)} The mapping $(u,v)\mapsto \Delta^{hp}(u,v)$ has a.s. a continuous extension to $\t^{hp}_\infty\times \t^{hp}_\infty$.

\noindent{\rm (ii)} A.s., for every $u,v\in \t^{hp,\circ}_\infty$, the property $\Delta^{hp}(u,v)=0$ holds if and only if 
$\Delta^{hp,\circ}(u,v)=0$.

\noindent{\rm(iii)}  A.s., for every $u,v\in\partial \t^{hp}_\infty$, the property $\Delta^{hp}(u,v)=0$ holds if and only if
$\Lambda_w>0$ for every $w\in]u,v[$, or for every $w\in]v,u[$.
\end{lemma}

\proof Property (i) can be derived by minor modifications of the proof of \cite[Proposition 31]{Disks}, noting that we
may restrict our attention to
the bounded subtree obtained by truncating $\t^{hp}_\infty$ at height $\mathbf{L}_r$ for some $r>0$. We
omit the details. As for (ii) and (iii), there is an additional complication due to the fact that
it is not immediately clear why we can restrict our attention to a bounded subtree.  Let us explain
the argument for (iii), which is the property we use below. The fact that 
$\Lambda_w>0$ for every $w\in]u,v[$ implies $\Delta^{hp}(u,v)=0$ is easy and left to the reader.
Suppose then that 
$u,v$ are distinct points of $\partial \t^{hp}_\infty$ are such that $\Delta^{hp}(u,v)=0$. Without loss
of generality we can assume that $[u,v]$ is compact, and we then have to check that 
$\Lambda_w>0$ for every $w\in]u,v[$. Recall the notation $\llbracket u,\infty\llbracket$ for
the unique geodesic ray from $u$ in the tree $\t^{hp}_\infty$, and 
$\rrbracket u,\infty\llbracket=\llbracket u,\infty\llbracket\backslash\{u\}$. We claim that, for every
$\delta>0$, we can find points $u'\in\rrbracket u,\infty\llbracket$ and $v'\in\rrbracket v,\infty\llbracket$
such that $\Lambda_{u'}<\delta$ and $\Lambda_{v'}<\delta$ and there exist
$w_0=u',w_1,\ldots,w_p=v'\in [u,v]$ such that
\begin{equation}
\label{cons-hp-tech1}
\sum_{i=1}^p \Delta^{hp,\circ}(w_{i-1},w_i) <\delta.
\end{equation}
If the claim holds, we can use \cite[Proposition 32$\,$(ii)]{Disks} to see that necessarily 
$\Lambda_w>0$ for every $w\in]u,v[$ (the point is the fact that all $w_i$'s belong to 
$[u,v]$, and thus we are dealing with a compact subtree of $\t^{hp}_\infty$). So it remains 
to prove our claim. First note that we can find $u'\in\rrbracket u,\infty\llbracket$ and $v'\in\rrbracket v,\infty\llbracket$
such that $\Lambda_{u'}<\delta/2$, $\Lambda_{v'}<\delta/2$ and $\Delta^{hp}(u',v')<\delta/2$ and 
in particular there exist $w_0,w_1,\ldots,w_p\in \t^{hp,\circ}_\infty$ 
such that \eqref{cons-hp-tech1} holds with $\delta$ replaced by $\delta/2$. It may happen that
some of the $w_i$'s do not belong to $[u,\infty)$, but then we can replace $u'$ by $w_{j+1}$, where 
$j=\max\{i:w_i\notin[u,\infty)\}$, noting that necessarily $w_{j+1}\in \rrbracket u,\infty\llbracket$ 
(otherwise $\Delta^{hp,\circ}(w_j,w_{j+1})$ would be $\infty$) and $\Lambda_{w_{j+1}}<\delta$
by the bound $|\Lambda_{w_i}-\Lambda_{w_{i-1}}|\leq \Delta^{hp,\circ}(w_{i-1},w_i)$ for $1\leq i\leq j+1$. Therefore we can assume that
all $w_i$'s belong to $[u,\infty)$, and then a symmetric argument shows that we can assume 
that they all belong to $[u,v]$ as desired. \endproof

%
%
%Again one verifies that the mapping $(v,w)\mapsto \Delta^{hp}(v,w)$ has a continuous extension to $\t^{hp}_\infty\times\t^{hp}_\infty$. 
%This is derived from the analogous result for the Brownian disk in \cite[Proposition 31]{Disks},
%using Proposition \ref{Palm} above to make the connection between the setting of \cite{Disks} and 
%the present one. Furthermore, if $u,v$ are two points of $\t^{hp}_\infty$ such that
%$\Lambda_u=\Lambda_v=0$, the property $\Delta^{hp}(u,v)=0$ holds if and only if
%we have $\Lambda_w>0$ for every $w\in]u,v[$ or for every $w\in]v,u[$. This property is the analog of
%Proposition 32$\,$(ii) in \cite{Disks}, and it can be derived from the latter result again by using Proposition \ref{Palm}.

We set  $\H_\infty=\t^{hp}_\infty/\{\Delta^{hp}=0\}$, and we let $\Pi^{hp}$ denote the canonical projection 
from $\t^{hp}_\infty$ onto $\H_\infty$. We equip
$\H_\infty$ with the distance induced by $\Delta^{hp}$ and the volume measure
which is the pushforward of the volume measure on $\t^{hp}_\infty$ under the canonical projection.
We observe that $\H_\infty$ has a distinguished vertex, namely the root $\rho$ of $\t^{hp}_\infty$ (or bottom of the spine).
By Lemma \ref{cons-hp-tech} (iii), the equivalence class of $\rho$ 
in the quotient $\t^{hp}_\infty/\{\Delta^{hp}=0\}$ must be a singleton,
since there are points of $\t^{hp}_\infty$ with zero label arbitrarily close to $\rho$, both on the 
left side and on the right side of the spine.

\begin{definition}
\label{def-half-plane}
The random pointed locally compact measure metric space $\H_\infty$ is called the Brownian half-plane.
\end{definition}

At the end of Section \ref{sec:consist}, we will explain why this definition is consistent with the one found in 
\cite{BMR} or in \cite{GM0}. The Brownian half-plane enjoys the same scale invariance property as the Brownian plane: Recalling  the notation
$\lambda\cdot E$ introduced in Section \ref{Br-plane},
$\lambda\cdot \H_\infty$ has the same distribution as $\H_\infty$, for every $\lambda>0$.
The boundary $\partial \H_\infty$
 is defined by $\partial\H_\infty:=\Pi^{hp}(\partial\t^{hp}_\infty)$ (one can prove
 that $\H_\infty$ is homeomorphic to the usual half-plane and then $\partial\H_\infty$ is also 
 the set of all points of $\H_\infty$ that have no neighborhood homeomorphic to an open disk, but we do not
 need these facts here). As noted in Section \ref{sec:pseudo},
 for any $v\in\t^{hp}_\infty$, $\Lambda_v$ is equal to the distance from
 $\Pi^{hp}(v)$ to the boundary $\partial\H_\infty$.

Let $r>0$ and let
  $\t^{hp}_{\infty,r}$ be the closed subset of $\t^{hp}_\infty$ consisting of the
 part $[0,\mathbf{L}_r]$ of the spine and the subtrees branching off $[0,\mathbf{L}_r]$. 
 
 \begin{lemma}
 \label{ball-root}
 We have 
 $$\inf_{v\notin \t^{hp}_{\infty,r}} \Delta^{hp}(\rho,v)>0\;,\quad\hbox{a.s.}$$
 \end{lemma}
 
 \proof  We argue by
 contradiction and assume that there is a sequence $(u_n)_{n\geq 1}$ in the complement of 
 $\t^{hp}_{\infty,r}$ such that $\Delta^{hp}(\rho,u_n)\la 0$ as $n\to\infty$. Suppose that infinitely many
 points of this sequence belong to $[\rho,\infty)$. Let $v_{(r)}$ be the last point of $\t^{hp}_{\infty,r}\cap\partial \t^{hp}_\infty$ visited by
 the exploration of $\t^{hp}_\infty$, and note that $\Pi^{hp}(\rho)\not =\Pi^{hp}(v_{(r)})$, by Lemma \ref{cons-hp-tech} (iii). Then an argument very similar to
 the proof of Lemma \ref{cons-hp-tech} (iii) shows that we can find another sequence
 $(v_n)_{n\geq 1}$ with $v_n\in\llbracket v_{(r)},\infty\llbracket$ and such that we still
 have $\Delta^{hp}(\rho,v_n)\la 0$ as $n\to\infty$. In particular, $\Lambda_{v_n}\la 0$,
 and this implies that $v_n\la v_{(r)}$ in $\t^{hp}_\infty$, and thus $\Delta^{hp}(v_{(r)},v_n)\la 0$ as $n\to\infty$. Finally we get $\Delta^{hp}(\rho,v_{(r)})=0$,
 which is a contradiction. The case when 
  infinitely many
 points of this sequence belong to $(\infty,\rho]$ is treated in a symmetric manner. \endproof
 
 It follows from Lemma \ref{ball-root}  that $\Pi^{hp}(\t^{hp}_{\infty,r})$ contains 
a ball of positive radius centered at $\rho$ in $\H_\infty$. Then, by scale invariance, we have a.s.
$$\lim_{r\to\infty}\Big(\inf_{v\notin \t^{hp}_{\infty,r}} \Delta^{hp}(\rho,v)\Big) =+\infty.$$
This implies in particular that $\H_\infty$ is boundedly compact (any ball centered at $\rho$
is contained in the image of a compact subtree of $\t^{hp}_\infty$ under $\Pi^{hp}$). 

Our next goal is to prove that $\H_\infty$ is the tangent cone in distribution of the
pointed Brownian disk at a point chosen uniformly on its boundary --- this will eventually allow us to make the connection
with previous definitions of the Brownian half-plane. 
Recall from the end of Section \ref{sec:point-disk} the notation $\bar\D^{\bullet,a}_z$ for the pointed measure metric space obtained from $\D^{\bullet\bullet,a}_z$ by
 ``forgetting'' the first distinguished point (so $\bar\D^{\bullet,a}_z$ is pointed at a point chosen
 uniformly on its boundary). 
% 
% For $\lambda>0$, we use the notation $\lambda \cdot \bar\D^{\bullet,a}_z$ for the
%same space $\bar\D^{\bullet,a}_z$ with the metric multiplied by the factor $\lambda>0$.

\begin{theorem}
\label{conv-hp}
Let $z>0$ and $a>0$. We have
$$\lambda \cdot \bar \D^{\bullet,a}_z \build{\la}_{\lambda\to\infty}^{\rm(d)} \H_\infty,$$
in distribution in the sense of the local Gromov-Hausdorff-Prokhorov convergence.
\end{theorem}

We give below the proof of Theorem \ref{conv-hp} for $a=1$, but a scaling argument yields the general case.
Before we proceed to the proof of Theorem \ref{conv-hp}, we start with some preliminary estimates.
We consider again a triple $(Y,\wt\mm,\wt\mm')$ distributed as
explained at the beginning of Section \ref{sec-ident} with $a=1$.
%(so $(Y^{(z)},\wt\mm^{(z)},\wt\mm'^{(z)})$ is distributed as $(Y,\wt\mm,\wt\mm')$ conditioned
%on $\{\z^Y=z\}$). 
Recall that the random path 
$Y$ is defined on the interval $[0,T^Y]$, that $Y_0=1$ and $T^Y=\inf\{t\geq 0: Y_t=0\}$. For every
$\ve\in(0,1)$, we also set
$$T_\ve:=\inf\{t\geq 0:Y_t=\ve\}.$$
We let  $\wt\mm^\ve(\dd t\dd \omega)$, resp. $\wt\mm'^\ve(\dd t\dd \omega)$, be the image of $\mathbf{1}_{[T_\ve,T^Y]}(t)\,\wt\mm(\dd t\dd\omega)$,
resp. of $\mathbf{1}_{[T_\ve,T^Y]}(t)\,\wt\mm'(\dd t\dd\omega)$, under the mapping $(t,\omega)\mapsto (t-T_\ve,\omega)$.
We also set $Y^\ve_t:=Y_{T_\ve +t}$ for $0\leq t\leq T^Y-T_\ve$. Recall that
$$\z^Y = \int \wt\mm(\dd t\dd \omega) \z_0(\omega) + \int \wt\mm'(\dd t\dd \omega)\z_0(\omega),$$
and also set
$$\z^{Y,\ve}:= \int \wt\mm^\ve(\dd t\dd \omega) \z_0(\omega) + \int \wt\mm'^\ve(\dd t\dd \omega)\z_0(\omega).$$
Set $\Gamma^\ve:=(Y^\ve,\wt\mm^\ve,\wt\mm'^\ve)$, and observe that $\Gamma^\ve$
is a coding triple in the sense of Section \ref{coding-infinite}. Moreover, the conditional distribution of
$\Gamma^\ve$ knowing $\z^{Y,\ve}=z$ is $\check\Theta^{(\ve)}_z$. Our first goal is to show that the conditional distribution
of $\Gamma^\ve$ given $\z^Y=z$ is close to its unconditional distribution when $\ve\to 0$.

From \eqref{Laplace-ex}, we have
$$\E[e^{-\lambda \z^Y}]=\Big(1+\sqrt{2\lambda/3}\,\Big)^{-3}\;,\quad 
\E[e^{-\lambda \z^{Y,\ve}}]=\Big(1+\ve\sqrt{2\lambda/3}\,\Big)^{-3}.$$
Furthermore, we may write $\z^Y=\z^{Y,\ve}+ \wh\z^{Y,\ve}$, where $\z^{Y,\ve}$ and $\wh\z^{Y,\ve}$ are independent
(more precisely, $\wh\z^{Y,\ve}$ is independent of $\Gamma^\ve$).
Hence
\begin{equation}
\label{LapT}
\E[e^{-\lambda \wh\z^{Y,\ve}}]= \Bigg(\frac{1+\ve\sqrt{2\lambda/3}}{1+\sqrt{2\lambda/3}}\Bigg)^3.
\end{equation}
The distribution of $\wh \z^{Y,\ve}$ can be written in the form
$$\ve^3\,\delta_0(\dd y)+ \wh\Upsilon_\ve(\dd y),$$
where the measure $\wh\Upsilon_\ve(\dd y)$ is supported on $(0,\infty)$.
To simplify notation, we also write $\varphi(y)=y\,h_1(y)$ for the density of 
$\z^Y$ and $\varphi_\ve(y)=y\,h_\ve(y)$ for the density of $\z^{Y,\ve}$. 

\begin{lemma}
\label{var-dis0}
We have $\wh\Upsilon_\ve(\dd y)=\wh\varphi_\ve(y)\,\dd y$, where the functions $\wh\varphi_\ve(y)$ satisfy
$$\lim_{\ve\to 0} \wh\varphi_\ve(y)= \varphi(y)$$ 
uniformly
on every interval of the form $[\delta,\infty)$, $\delta>0$.
\end{lemma}

\proof %We use the explicit form of $\wh\varphi_\ve$ that can be derived from the Laplace transform\eqref{LapT}. 
From \eqref{LapT}, we have, for every $\lambda>0$,
$$
\int_0^\infty \wh\Upsilon_\ve(\dd y)\,e^{-\lambda y}=  \Bigg(\frac{1+\ve\sqrt{2\lambda/3}}{1+\sqrt{2\lambda/3}}\Bigg)^3 
-\ve^3.$$
%Hence we get $\wh\varphi_\ve(z)=\frac{3}{2}\gamma_{(\ve)}(\frac{3z}{2})$, where the function $\gamma_{(\ve)}$ satisfies
Now observe that
\begin{equation}
\label{formu-gamma}
\Bigg(\frac{1+\ve\sqrt{\lambda}}{1+\sqrt{\lambda}}\Bigg)^3 -\ve^3
= \frac{(1-\ve)^3+3\ve(1-\ve)^2(1+\sqrt{\lambda}) + 3\ve^2(1-\ve)(1+\sqrt{\lambda})^2}{(1+\sqrt{\lambda})^3}
\end{equation}
where we have expanded $(1+\ve\sqrt{\lambda})^3=((1-\ve)+\ve(1+\sqrt{\lambda}))^3$. 
It follows from formulas (A.1), (A.2), (A.3) in the Appendix that the Laplace transform of the function $\chi_{(\ve)}$
defined by
$$\chi_{(\ve)}(y) = (1-\ve)^3 \chi_3(y) + 3\ve(1-\ve)^2 \chi_2(y) + 3\ve^2(1-\ve) \chi_1(y).
$$
is the quantity in \eqref{formu-gamma}. Consequently, we have $\wh\Upsilon_\ve(\dd y)=\wh\varphi_\ve(y)\,\dd y$ 
with $\wh\varphi_\ve(y)=\frac{3}{2}\chi_{(\ve)}(\frac{3y}{2})$. Furthermore,
the explicit formulas for $\chi_1,\chi_2,\chi_3$ show that $\chi_{(\ve)}(y)$ converge 
to $\chi_3(y)$ as $\ve\to0$, uniformly
on every interval of the form $[\delta,\infty)$, $\delta>0$. The result of the proposition follows 
since
$\varphi(y)=\frac{3}{2}\chi_3(\frac{3y}{2})$ by definition.
\endproof

\begin{lemma}
\label{var-dis}
Let $z>0$. The total variation distance between the conditional distribution of $\Gamma^\ve$ knowing that
$\z^Y=z$ and the unconditional distribution of $\Gamma^\ve$ converges to $0$
as $\ve\to 0$.
\end{lemma}

\rem We have made a canonical choice for the conditional distribution $\check\Theta^{(1)}_z$ of 
$(Y,\wt\mm,\wt\mm')$ knowing $\z^Y=z$, and so the conditional distribution of $\Gamma^\ve$ knowing that
$\z^Y=z$ is also well defined for every $z$.

\proof The equality $\z^Y=\z^{Y,\ve}+ \wh\z^{Y,\ve}$ gives
\begin{equation}
\label{var-dis-tec2}
\varphi(x)=\ve^3 \varphi_\ve(x)+ \int_0^x \varphi_\ve(y)\wh\varphi_\ve(x-y)\,\dd y
\end{equation}
for every $x>0$.
 Let $G$ and $g$ be measurable functions defined respectively on
$\W\times M_p(\R_+\times\S)^2$ and on $\R_+$, such that $0\leq G\leq 1$
and $0\leq g\leq 1$.
Then,
\begin{align*}
\E[G(\Gamma^\ve)g(\z^Y)]&=\E[G(\Gamma^\ve)g(\z^{Y,\ve}+\wh\z^{Y,\ve})]\\
&=\int \dd z\,\wh\varphi_\ve(z)\,\E[G(\Gamma^\ve)g(\z^{Y,\ve}+z)] +\ve^3\,\E[G(\Gamma^\ve)g(\z^{Y,\ve})]\\
&=\E\Big[G(\Gamma^\ve)\int_{\z^{Y,\ve}}^\infty \dd z\,g(z)\,\wh\varphi_\ve(z-\z^{Y,\ve})\Big]
+\ve^3\int \dd z\,g(z)\,\varphi_\ve(z)\,\E[G(\Gamma^\ve)\,|\, \z^{Y,\ve}=z]\\
&=\int \dd z\,g(z)\,\Big( \E[G(\Gamma^\ve)\mathbf{1}_{\{\z^{Y,\ve}<z\}}\,\wh\varphi_\ve(z-\z^{Y,\ve})]
+\ve^3\,\varphi_\ve(z)\,\E[G(\Gamma^\ve)\,|\, \z^{Y,\ve}=z]\Big).
\end{align*}
Recalling that the density of $\z^Y$ is $\varphi$, it follows that we have $\dd z$ a.e.,
\begin{equation}
\label{var-dis-tec10}
\E[G(\Gamma^\ve)\,|\, \z^{Y}=z]
=  \E\Bigg[G(\Gamma^\ve)\mathbf{1}_{\{\z^{Y,\ve}<z\}}\,\frac{\wh\varphi_\ve(z-\z^{Y,\ve})}{\varphi(z)}\Bigg]
+ \ve^3 \E[G(\Gamma^\ve)\,|\, \z^{Y,\ve}=z]\,\frac{\varphi_\ve(z)}{\varphi(z)},
\end{equation}
where we observe that $\E[G(\Gamma^\ve)\mid \z^{Y,\ve}=z]=\check\Theta^{(\ve)}_z(G)$ is well defined for every $z$. 
We now want to argue that \eqref{var-dis-tec10} holds for {\it every} $z>0$ and not only $\dd z$ a.e.
To this end, it is enough to consider the special case $G(\w,\mu,\mu')=\exp(-f(\w)-\langle \mu,h\rangle-\langle\mu',h'\rangle)$ where, $f,h,h'$ are nonnegative functions, $f$ is bounded and continuous on $\W$, $h$ and $h'$ are bounded and continuous
on $\R_+\times \S$ and both $h$ and $h'$ vanish on $\{(t,\omega):\sigma(\omega)\leq\delta\}$ for some $\delta>0$. In that case,
using a scaling argument and Corollary
 \ref{finite-height-cor}, one checks that both sides of  \eqref{var-dis-tec10} are continuous functions of
 $z$, so that they must be equal for every $z>0$. 
 
From \eqref{density-exit}, we have $\varphi_\ve(z)=O(\ve)$ as $\ve\to 0$, hence, for every fixed $z>0$,
\begin{equation}
\label{var-dis-tec11}
\lim_{\ve\to 0}  \ve^3 \E[G(\Gamma^\ve)\mid \z^{Y,\ve}=z]\,\frac{\varphi_\ve(z)}{\varphi(z)} =0,
\end{equation}
uniformly in the choice of $G$. On the other hand,
using Lemma \ref{var-dis0} and the fact that $\z^{Y,\ve}\la 0$ as $\ve\to 0$, we
have
\begin{equation}
\label{var-dis-tech}
\lim_{\ve\to 0} \mathbf{1}_{\{\z^{Y,\ve}<z\}}\,\frac{\wh\varphi_\ve(z-\z^{Y,\ve})}{\varphi(z)} =1
\end{equation}
almost surely. Moreover, using \eqref{var-dis-tec2}, we have
$$\E\Bigg[\mathbf{1}_{\{\z^{Y,\ve}<z\}}\,\frac{\wh\varphi_\ve(z-\z^{Y,\ve})}{\varphi(z)}\Bigg]
=\frac{1}{\varphi(z)}\int_0^z \dd y\,\varphi_\ve(y)\wh\varphi_\ve(z-y)
=\frac{1}{\varphi(z)}(\varphi(z)-\ve^3 \varphi_\ve(z)),$$
which tends to $1$ as $\ve\to 0$. By Scheff\'e's lemma, the convergence \eqref{var-dis-tech} also
holds in $L^1$. The statement of the lemma then follows from \eqref{var-dis-tec10}
and \eqref{var-dis-tec11}. \endproof

 \noindent
 {\it Proof of Theorem \ref{conv-hp}.} The proof is based on a coupling argument
 relying on Lemma \ref{var-dis}. 
  If $E$ is a pointed metric space, we use the notation $B_r(E)$ for the closed ball of radius $r$ centered
 at the distinguished point. The theorem will follow if we can prove that, for every $K>0$ and every $\delta>0$,
 if $\lambda$ is large enough we can couple $\H_\infty$ and $\bar\D^{\bullet,1}_z$ in such a way that
 the balls $B_K(\lambda\cdot \bar\D^{\bullet,1}_z)$ and $B_K(\H_\infty)$ are isometric with probability at least $1-\delta$
 (with an isometry preserving the volume measure and the distinguished point).
 Equivalently, recalling that $\lambda\cdot \H_\infty$ has the same distribution as $\H_\infty$, it suffices to prove that, for $\eta>0$ small enough, 
$\H_\infty$ and $\bar\D^{\bullet,1}_z$ can be coupled so that 
$B_{\eta}(\bar\D^{\bullet,1}_z)$ and $B_{\eta}(\H_\infty)$ are isometric with probability at least $1-\delta$
(again with an isometry preserving the volume measure and the distinguished point). 

As explained at the end 
of Section \ref{sec:point-disk}, we may and will assume that $\bar\D^{\bullet,1}_z$ is constructed from
a coding triple $(Y^{(z)},\wt\mm^{(z)},\wt\mm'^{(z)})$ distributed according to $\check\Theta^{(1)}_z$. 
The labeled tree associated with $(Y^{(z)},\wt\mm^{(z)},\wt\mm'^{(z)})$ is denoted by  $\t^{(z)}$, 
and we write $\Delta^{(z),\circ}$ and 
$\Delta^{(z)}$ for the pseudo-distance functions on $\t^{(z)}$, so that $\Delta^{(z)}$ induces the metric on $\bar\D^{\bullet,1}_z$.
The set of all points of $\t^{(z)}$ with positive label is denoted by $\t^{(z),\circ}$.

For $\ve\in(0,1)$, let $\Gamma^{(z),\ve}$ be defined as $\Gamma^\ve$ but replacing the triple $(Y,\wt\mm,\wt\mm')$ by $(Y^{(z)},\wt\mm^{(z)},\wt\mm'^{(z)})$
(so $\Gamma^{(z),\ve}$ is distributed as $\Gamma^\ve$ conditioned on $\z^Y=z$). 
Let $\check\Gamma^{(z),\ve}$, resp. $\check\Gamma^\ve$, denote the image of $\Gamma^{(z),\ve}$, resp. $\Gamma^\ve$, under the time reversal operation $\mathbf{TR}$ in \eqref{time-rev}. 
We fix $\delta>0$ and claim that:
\begin{itemize}
\item[$1.$] For $\ve\in(0,1)$ small enough, the triples $(Y^{(z)},\wt\mm^{(z)},\wt\mm'^{(z)})$ and $(R,\pp,\pp')$ can be coupled in such a way that
the equality
\begin{equation}
\label{coup-id}
\check\Gamma^{(z),\ve}= \Big((R_t)_{0\leq t\leq \mathbf{L}_\ve}, \mathbf{1}_{[0,\mathbf{L}_\ve]}(t)\,\wt\pp(\dd t\dd\omega), \mathbf{1}_{[0,\mathbf{L}_\ve]}(t)\,\wt\pp'(\dd t\dd\omega)\Big)
\end{equation}
holds with probability at least $1-\frac{\delta}{2}$. 
\item[$2.$] For $\ve\in(0,1)$ small enough, we can choose $\eta_0>0$ so that for every $0<\eta\leq \eta_0$, we have
$$B_\eta(\bar\D^{\bullet, 1}_z)=B_\eta(\H_\infty)$$
on the event where \eqref{coup-id} holds, except possibly on an event of probability at most $\frac{\delta}{2}$. 
\end{itemize} 

Clearly the theorem follows from Properties 1 and 2. 
Property 1 is a consequence of Lemma \ref{var-dis}: just note that the distribution of the coding triple in the right-hand side of 
\eqref{coup-id} is the (unconditional) distribution of $\check\Gamma^\ve$. 

It remains to verify Property 2. We fix $\ve>0$ small enough so that we can apply Property 1. We then assume that 
the triples $(Y^{(z)},\wt\mm^{(z)},\wt\mm'^{(z)})$ and $(R,\pp,\pp')$ have been coupled in such a way that the event where
\eqref{coup-id} holds has probability greater than $1-\frac{\delta}{2}$, and we
denote the latter event by $\f$. 
We argue on the intersection $\f\cap\f'$, where $\f'$ denotes
 the event where $W_*(\omega_i)=0$ for at least one atom $(t_i,\omega_i)$ of $\wt\mm^{(z)}$ or of $\wt\mm'^{(z)}$ such that $t_i<T_\ve^{(z)}:=\inf\{t\geq 0:Y^{(z)}_t=\ve\}$. Clearly we can also assume that $\f'$ has probability greater than $1-\frac{\delta}{6}$
by choosing $\ve$ even smaller if necessary.

Recall the notation $\t^{hp}_{\infty,r}$ introduced before Lemma \ref{ball-root}. From this lemma, we know
that, for $\eta>0$ small enough, the set $\{v\in\t^{hp}_\infty: \Delta^{hp}(\rho,v)\leq 4\eta\}$ will be contained 
in $\t^{hp}_{\infty,\ve}$, except on an event of probability at most $\frac{\delta}{6}$. Moreover, if the latter property holds, we claim that we have also, for every $u,v\in\t^{hp,\circ}_\infty$
such that $\Delta^{hp}(\rho,u)\leq \eta$ and $\Delta^{hp}(\rho,v)\leq \eta$,
\begin{equation}
\label{tech-dist}
\Delta^{hp}(u,v) = \build{\inf_{u_0=u,u_1,\ldots,u_p=v}}_{u_1,\ldots,u_{p-1}\in \t^{hp}_{\infty,\ve}\cap\t^{hp,\circ}_\infty}^{} \sum_{i=1}^p\Delta^{hp,\circ}(u_{i-1},u_i).
\end{equation}
In other words, in formula \eqref{delta-infty} applied to $\Delta^{hp}(u,v)$, we may restrict the infimum to the case where
all $u_i$'s belong to  $\t^{hp}_{\infty,\ve}$. Let us justify \eqref{tech-dist}. Assume that $\Delta^{hp}(\rho,u)\leq \eta$ and $\Delta^{hp}(\rho,v)\leq \eta$ (so that in particular $\Delta^{hp}(u,v)\leq 2\eta$) and $u_0=u,u_1,\ldots,u_q
\in \t^{hp,\circ}_\infty $ are such that 
$$\sum_{i=1}^q \Delta^{hp,\circ}(u_{i-1},u_i)<\Delta^{hp}(u,v)+\eta.$$
It then follows that $\Delta^{hp}(u,u_q)< 3\eta$ and $\Delta^{hp}(\rho,u_q)<4\eta$ which implies $u_q\in  \t^{hp}_{\infty,\ve}$.

Furthermore, when applying formula \eqref{delta0} to compute the quantities $\Delta^{hp,\circ}(u_{i-1},u_i)$ in the right-hand side 
of \eqref{tech-dist}, it is enough to consider the case when the interval $[u_{i-1},u_i]$ (resp.~$[u_i,u_{i-1}]$) is 
contained in $\t^{hp}_{\infty,\ve}$, because otherwise this interval contains 
$\t^{hp}_\infty\backslash \t^{hp}_{\infty,\ve}$ and then the infimum of labels on $[u_{i-1},u_i]$ is $0$. 
To summarize, on the event where \eqref{tech-dist} holds for every $u,v\in\t^{hp,\circ}_\infty$
such that $\Delta^{hp}(\rho,u)\leq \eta$ and $\Delta^{hp}(\rho,v)\leq \eta$, we get that the value of
$\Delta^{hp}(u,v)$ for such points $u$ and $v$ is determined by the tree $\t^{hp}_{\infty,\ve}$ and the labels on this tree.

On the event where \eqref{tech-dist} holds, we thus get that the ball $B_\eta(\H_\infty)$ can be written as a function of the coding triple
$$\Big((R_t)_{0\leq t\leq \mathbf{L}_\ve}, \mathbf{1}_{[0,\mathbf{L}_\ve]}(t)\,\wt\pp(\dd t\dd\omega), \mathbf{1}_{[0,\mathbf{L}_\ve]}(t)\,\wt\pp'(\dd t\dd\omega)\Big)$$
since the tree $\t^{hp}_{\infty,\ve}$ and the labels on this tree are functions of this triple
(and also the distinguished point of $B_\eta(\H_\infty)$ corresponds to the root
of this coding triple). To complete the argument (recalling that
we assume \eqref{coup-id}), we need to justify that $B_\eta(\bar\D^{\bullet,1}_z)$ is given by the same function applied to the triple $\check\Gamma^{(z),\ve}$,
except possibly on a set of small probability. To get this, recall that $\D^{\bullet,1}_z$ is obtained by applying $\Xi^\bullet$ to the
snake trajectory $\Omega(Y^{(z)},\wt\mm^{(z)},\wt\mm'^{(z)})$. With the coding triple $\check\Gamma^{(z),\ve}$ 
we associate a labeled tree 
$\t^{(z)}_\ve$, which is identified to a subtree of the labeled tree $\t^{(z)}$, and, modulo this identification, $\t^{(z)}_\ve$ is rooted at the
top of the spine of the tree $\t^{(z)}$, which corresponds to the distinguished point  $\alpha$ of $\bar\D^{\bullet,1}_z$. 
We claim that the image of 
$\t^{(z)}_\ve$ (viewed as a subset of $\t^{(z)}$) under the canonical projection from $\t^{(z)}$ onto $\bar\D^{\bullet,1}_z$ must contain a neighborhood of $\alpha$.
As in the proof of Lemma \ref{ball-root}, this property holds because the equivalence class of $\alpha$ in $\t^{(z)}/\{\Delta^{(z)}=0\}$ is a singleton, which is
a consequence of the fact that two points $u$ and $v$ of 
$\t^{(z)}$ with zero label are identified in $\bar\D^{\bullet,1}_z$ if and only if labels stay positive on the interval $]u,v[$, or on the interval $]v,u[$
(see the discussion after Proposition \ref{pointed-disk}).

It follows from the preceding claim that,
for $\eta$ small enough, we have $\Delta^{(z)}(\alpha,v)>4\eta$ whenever $v\notin \t^{(z)}_\ve$, except on an event
of probability at most $\frac{\delta}{6}$. Discarding the latter event of small probability, the same argument as above shows that the analog of \eqref{tech-dist} 
holds for every $u,v\in\t^{(z),\circ}$ such that $\Delta^{(z)}(\alpha,u)<\eta$ and 
$\Delta^{(z)}(\alpha,v)<\eta$, provided we replace $\Delta^{hp}$ by
$\Delta^{(z)}$, $\t^{hp,\circ}_\infty$ by $\t^{(z),\circ}$, and $\t^{hp}_{\infty,\ve}$  by $\t^{(z)}_\ve$. Furthermore, the quantities 
$\Delta^{(z),\circ}(u_{i-1},u_i)$ appearing in this analog can be computed from the labeled tree $\t^{(z)}_\ve$ (here we use
our definition of $\f'$, which implies that $\t^{(z)}\backslash \t^{(z)}_\ve$ contains points with zero label). 

It follows from the preceding discussion that, on the event $\f$ that has probability at least $1-\frac{\delta}{2}$, and except on an event of probability at most $\frac{\delta}{2}$, the ball $B_\eta(\bar\D^{\bullet,1}_z)$
is obtained from the triple $\check\Gamma^{(z),\ve}$ by applying the same function that can be used to
get the ball $B_\eta(\H_\infty)$ from the triple in the right-hand side of \eqref{coup-id}.
The desired result follows. \hfill$\square$

\section{Applications}
\label{appli}

\subsection{Infinite Brownian disks in the Brownian plane}
\label{iBdBp}

Recall the construction of the Brownian plane $(\B\P_\infty,\Delta^p)$ from the coding triple $(X,\ll,\rr)$ in Section \ref{Br-plane} and note that the same triple was also considered in Section \ref{cod-disk}.
We use the notation $(\t^p_\infty,(\Lambda_v)_{v\in\t^p_\infty})$ for the labeled tree associated with 
the triple $(X,\ll,\rr)$, and we
write $\Pi^p$ for the canonical projection from $\t^p_\infty$ onto $\B\P_\infty$. The distinguished point $\rho$ of 
$\B\P_\infty$ is the image of the root of $\t^p_\infty$ under $\Pi^p$. 

To simplify notation, for every $r>0$, we write  $B(r)=B_r(\BP_\infty)$ for the closed ball of radius $r$ centered at $\rho$ in $\BP_\infty$. The hull 
$B^\bullet(r)$ is then the subset of $\BP_\infty$ defined by saying that
$\BP_\infty \backslash B^\bullet(r)$ is the unique unbounded connected component of $\BP_\infty\backslash B(r)$ (this component is unique since 
$\BP_\infty$ is homeomorphic to the plane \cite{Plane}). 
Informally, $B^\bullet(r)$ is obtained by filling in the (bounded) holes in $B(r)$. As in the introduction, it will be convenient to use the notation
$$\check B^\bullet(r)= \BP_\infty \backslash B^\bullet(r).$$
One can give an explicit description of $\check B^\bullet(r)$ in terms of the labeled tree $(\t^p_\infty,(\Lambda_v)_{v\in\t^p_\infty})$.
For $v\in \t^p_\infty$, we recall that $\llbracket v,\infty\llbracket$ is the geodesic ray from $v$
in $\t^p_\infty$. Then, $\check B^\bullet(r)=\Pi^p(F_r)$, where
\begin{equation}
\label{hull-expl}
F_r:=\{v\in\t^p_\infty: \Lambda_w>r\hbox{ for every }w\in \llbracket v,\infty\llbracket\}.
\end{equation}
Similarly, the topological boundary of
$\check B^\bullet(r)$ (or of $B^\bullet(r)$) is $\partial\check B^\bullet(r)=\partial B^\bullet(r)= \Pi^p(\partial F_r)$, with
\begin{equation}
\label{bdry-expl}
\partial F_r=\{v\in\t^p_\infty: \Lambda_w=r\hbox{ and }\Lambda_w>r\hbox{ for every }w\in \rrbracket v,\infty\llbracket\},
\end{equation}
with the obvious notation $\rrbracket v,\infty\llbracket$. See formulas (16) and (17) in \cite{CLG}. We note that the intersection 
of the set $F_r$ with the spine of $\t^p_\infty$
is just the interval $(\mathrm{L}_r,\infty)$, where we
recall the notation $\mathrm{L}_r$ in \eqref{last-pass} (as in Section \ref{coding-infinite}, the spine is identified to $\R_+$). Following \cite{CLG},
we define the boundary size of $B^\bullet(r)$ to be $|\partial B^\bullet(r)|=\z^{(r)}$, where the quantity 
$\z^{(r)}$ is defined in \eqref{peri-hull}: $\z^{(r)}$ is the sum over all atoms $(t,\omega)$ of 
$\ll$ and $\rr$ such that $t>\mathrm{L}_r$ of the exit measures $\z_r(\omega)$ at level $r$ --- see formula (18) in \cite{CLG}.
We write $\mathrm{cl}(\check B^\bullet(r))=\check B^\bullet(r)\cup \partial B^\bullet(r)$ for the closure of
$\check B^\bullet(r)$, and similarly $\mathrm{cl}(F_r)=F_r\cup \partial F_r$. 

Recall that the intrinsic metric on an open subset $O$ of $\BP_\infty$ is defined by declaring that
the distance between two points $x$ and $y$ of $O$ is the infimum of the lengths of all
continuous curves $\gamma:[0,1]\la O$ such that $\gamma(0)=x$ and $\gamma(1)=y$. Here the
lengths are of course computed with respect to the metric $\Delta^p$ of $\BP_\infty$. 

\begin{theorem}
\label{infBdBp}
Let $r>0$. Then a.s.~the intrinsic metric on $\check B^\bullet(r)$ has a unique continuous extension to 
$\mathrm{cl}(\check B^\bullet(r))$, which is a metric on this set. We write $(\D^{\infty,(r)},\Delta^{\infty,(r)})$ for
the resulting random locally compact metric space, which is
equipped with the restriction of the volume measure on $\BP_\infty$ and pointed at $\Pi^p(\mathrm{L}_r)$.
Then, conditionally on $|\partial B^\bullet(r)|$,  $(\D^{\infty,(r)},\Delta^{\infty,(r)})$
is an infinite-volume Brownian disk with perimeter $|\partial B^\bullet(r)|$.
\end{theorem}

\proof Recall the notation $(X^{(r)},\wt\ll^{(r)},\wt\rr^{(r)})$ introduced in Section \ref{cod-disk},
and the fact that, conditionally on $\z^{(r)}=z$, this coding triple is distributed according to $\Theta_z$
(Proposition \ref{condi-dis2}). 
%Write $\Delta^{\infty,(r),\circ}$ and $\Delta^{\infty,(r)}$ for the 
%pseudo-metric functions associated with this coding triple. 
The construction of Section \ref{inf-Bdisk} produces,
from the triple $(X^{(r)},\wt\ll^{(r)},\wt\rr^{(r)})$, a random measure 
metric space $(\D^{\infty,(r)},\Delta^{\infty,(r)})$ such that, conditionally on $\z^{(r)}=z$, $(\D^{\infty,(r)},\Delta^{\infty,(r)})$
is an infinite-volume Brownian disk with perimeter $z$. Furthermore, $\D^{\infty,(r)}$ is obtained as a quotient
space of the labeled tree $\t^i_\infty$ coded by the  triple $(X^{(r)},\wt\ll^{(r)},\wt\rr^{(r)})$. Here we use the 
same notation $\t^i_\infty$ as in Section \ref{inf-Bdisk}, where we were dealing with a different triple 
distributed according to $\Theta_z$, but this should create no confusion. We write
$\Pi^i$ for the canonical projection from $\t^i_\infty$ onto $\D^{\infty,(r)}$.

It is easy to verify that the tree $\t^i_\infty$ can be identified with $\mathrm{cl}(F_r)$. The spine 
of $\t^i_\infty$ is identified with the part $[\mathrm{L}_r,\infty)$ of the spine of $\t^p_\infty$, and we observe that, for each atom
$(t_i,\omega_i)$ of $\ll$ or $\rr$ such that $t_i>\mathrm{L}_r$ (so that $\tr_r(\omega_i)$ shifted by $-r$ corresponds to an 
atom of $\wt\ll^{(r)}$ or $\wt\rr^{(r)}$), the genealogical tree $\t_{(\tr_r(\omega_i))}$ is
identified with $\{v\in\t_{(\omega_i)}: \Lambda_w>r\hbox{ for every }w\in\llbracket\rho_{(\omega_i)},v\llbracket\}$
(see the end of Section \ref{sna-tra}). 
The identification of $\t^i_\infty$ with $\mathrm{cl}(F_r)$ preserves labels, provided labels on $\mathrm{cl}(F_r)$ are shifted by $-r$. With a slight abuse of 
notation, if $u\in\t^i_\infty$,  we will also write $\Lambda_u$ for the label of the point 
of $\mathrm{cl}(F_r)$ corresponding to $u$ in the identification of 
$\t^i_\infty$ with $\mathrm{cl}(F_r)$ (so the label of $u$ in $\t^i_\infty$ is $\Lambda_u-r$).
 
Furthermore, two vertices of $\mathrm{cl}(F_r)$ are identified in the quotient space $\BP_\infty$
if and only if the corresponding vertices of $\t^i_\infty$ are identified in the quotient $\D^{\infty,(r)}$: to check this property in the case where the vertices 
belong to the boundary (the other case is immediate) we use the fact that
two vertices $u$ and $v$ of $\t^i_\infty$ with zero label are identified if and only if labels remain positive 
on one of the two intervals  $]u,v[$ and $]v,u[$ of the tree $\t^i_\infty$ (Lemma \ref{cons-inf-tech} (iii)). Thus we can 
identify $\D^{\infty,(r)}$ with the set $\Pi^p(\mathrm{cl}(F_r))= \mathrm{cl}(\check B^\bullet(r))$, in such a way that
$\partial \D^{\infty,(r)}$ is identified with $\partial \check B^\bullet(r)$, and this 
identification preserves the volume measures.

Modulo the preceding identification, both assertions of the theorem follow from the next lemma.

\begin{lemma}
\label{compar-dist}
Let $x$ and $y$ be two points of 
$\D^{\infty,(r)}\backslash \partial \D^{\infty,(r)}$, and let $x'$ and $y'$ be the corresponding points
in $\check B^\bullet(r)$. Then the intrinsic distance (relative to the open set $\check B^\bullet(r)$) between $x'$ and $y'$ coincides with 
$\Delta^{\infty,(r)}(x,y)$. 
\end{lemma}

\proof 
Let $\Delta^{\infty,(r),\circ}(v,w)$ be defined as in \eqref{delta0} for the labeled tree 
$(\t^i_\infty,(\Lambda_u-r)_{u\in \t^i_\infty})$ (recall that the label in $\t^i_\infty$ of a point $u\in\t^i_\infty$
is equal to $\Lambda_u-r$), so that $\Delta^{\infty,(r)}(v,w)$ is then given from $\Delta^{\infty,(r),\circ}(v,w)$ by formula
\eqref{delta-infty}. 

We first prove that the intrinsic distance between $x'$ and $y'$ is bounded above by $\Delta^{\infty,(r)}(x,y)$. 
To this end, let $v$ and $w$ be points of 
$\t^i_\infty$ such that $\Pi^i(v)=x$ and $\Pi^i(w)=y$. We claim
that, if $\Delta^{\infty,(r),\circ}(v,w)<\infty$, then $\Delta^{\infty,(r),\circ}(v,w)$ is the length of 
a continuous curve in $\check B^\bullet(r)$ that connects $x'$ to $y'$. 
Let us explain this. Without loss of generality, we may assume that
$$\Delta^{\infty,(r),\circ}(v,w)=\Lambda_v+\Lambda_w - 2\inf_{u\in[v,w]}\Lambda_u,$$
with $\inf_{u\in[v,w]}\Lambda_u>r$.
We let $v'$ and $w'$ be the points of $F_r$ corresponding to
$v$ and $w$ in the identification of $\t^i_\infty$ with $\mathrm{cl}(F_r)$ (in particular $\Lambda_{u'}=\Lambda_u$
and $\Lambda_{v'}=\Lambda_v$). We note that
the condition $\inf_{u\in[v,w]}\Lambda_u>r$ implies that the interval $[v,w]$ of $\t^i_\infty$
is also identified with the interval $[v',w']$ of $\t^p_\infty$ (in particular we have
$\Delta^{p,\circ}(v',w')=\Delta^{\infty,(r),\circ}(v,w)$), and furthermore $\Pi^p([v',w'])$ is contained
in $\check B^\bullet(r)$. By concatenating two simple geodesics starting from 
$\Pi^p(v')=x'$ and $\Pi^p(w')=y'$ respectively up to their merging time, as explained at the end of Section \ref{sec:pseudo}, we construct 
a path from $x'$ to $y'$ whose length is equal to $\Delta^{p,\circ}(v',w')$, and 
which stays in $\Pi^p([v',w'])\subset\check B^\bullet(r)$. This gives our claim. 
%As in Section \ref{coding-infinite}, write $(\ee_s)_{s\in\R}$ for the
%clockwise exploration of $\t^p_\infty$, and pick $s$ and $\tilde s$ such that $\ee_s=v'$,
%$\ee_{\tilde s}=w'$, and $[v',w']=\{\ee_h:h\in[s,\tilde s]\}$. For simplicity, we assume that $s\leq \tilde s$ (the
%other case is treated in a similar manner). We construct a path from $x'$ to $y'$
%by concatenating two pieces of ``simple geodesics'' (see e.g.~\cite[Section 2.6]{Uniqueness} for
%a definition of these geodesics in the Brownian map). For every $t$ such that $0\leq t\leq \Lambda_v-\inf_{[v,w]}\Lambda_u$, set
%$$\xi_t:= \inf\{h\geq s: \Lambda_{\ee_h}=\Lambda_v-t\}.$$
%Then $(\Pi^p(\ee_{\xi_t}),0\leq t\leq \Lambda_v-\inf_{[v,w]}\Lambda_u)$ gives a geodesic path 
%(with respect to $\Delta^p$)
%from $x'$ to a point $a$ of $\Pi^p([v',w'])$ with label $\inf_{[v,w]}\Lambda_u $, and this geodesic path stays
%in $\check B^\bullet(r)$. Similarly, by
%considering for $0\leq t\leq \Lambda_w-\inf_{[v,w]}\Lambda_u$,
%$$\tilde \xi_t:= \sup\{h\leq s': \Lambda_{\ee_h}=\Lambda_w-t\},$$
%we get a geodesic path from $y'$ to the same point $a$ that stays
%in $\check B^\bullet(r)$.
%The concatenation of these two
%geodesic paths gives a continuous curve from $x'$ to $y'$ in $\check B^\bullet(r)$
%with length $\Delta^{\infty,(r),\circ}(v,w)$ as required. 

From the definition of
$\Delta^{\infty,(r)}$ as an infimum, we now get that $\Delta^{\infty,(r)}(x,y)$ is bounded below by the infimum 
of lengths of continuous curves connecting $x'$ and $y'$ that stay in $\check B^\bullet(r)$.
We thus
obtain that the intrinsic distance between $x'$ and $y'$ (with respect to the open set $\check B^\bullet(r)$)
is bounded above by $\Delta^{\infty,(r)}(x,y)$.

It remains to prove the reverse bound. To this end, we need to verify that,
if $\gamma:[0,1]\to \check B^\bullet(r)$ is a continuous curve such that $\gamma(0)=x'$ and $\gamma(1)=y'$,
then the length of $\gamma$ is bounded below by $\Delta^{\infty,(r)}(x,y)$.
We write $\ov\gamma(t)$ for the point of $\D^{\infty,(r)}$ corresponding to $\gamma(t)$ in the identification of
$\mathrm{cl}(\check B^\bullet(r))$ with $\D^{\infty,(r)}$.
We may find $\delta>0$ such that $\Lambda_{\gamma(t)}>r+\delta$ for every $t\in[0,1]$ (recall that 
$\Lambda_z=\Delta^p(\rho,z)$ for every $z\in \BP_\infty$). Then, we may choose $n$ large enough 
so that $\Delta^p(\gamma(\frac{i-1}{n}),\gamma(\frac{i}{n}))<\delta/2$ for every $1\leq i\leq n$. 
The length of $\gamma$ is bounded below by
$\sum_{i=1}^n \Delta^p(\gamma(\frac{i-1}{n}),\gamma(\frac{i}{n}))$,
and so to get the desired result it suffices to verify that, for every $1\leq i\leq n$,
$$\Delta^p\Big(\gamma\Big(\frac{i-1}{n}\Big),\gamma\Big(\frac{i}{n}\Big)\Big)\geq \Delta^{\infty,(r)}\Big(\ov\gamma\Big(\frac{i-1}{n}\Big),\ov\gamma\Big(\frac{i}{n}\Big)\Big).
$$
Fix $1\leq i\leq n$, and recall the definition \eqref{delta-infty} of $\Delta^p(\gamma(\frac{i-1}{n}),\gamma(\frac{i}{n}))$ as an infimum over possible
choices of $u_0=\gamma(\frac{i-1}{n}),u_1,\ldots,u_p=\gamma(\frac{i}{n})$ in $\t^p_\infty$, where
we may restrict our attention to choices
of $u_0,u_1,\ldots,u_p$ such that $\Lambda_{u_j}>r+\delta/2$ (use $\Delta^{p}(u,v)\geq |\Lambda_v-\Lambda_u|$) and 
$\Delta^{p,\circ}(u_{j-1},u_j)<\delta/2$ for every $1\leq j\leq p$. It suffices to consider one such choice and to prove that
\begin{equation}
\label{tech-iBdBp}
\sum_{j=1}^p \Delta^{p,\circ}(u_{j-1},u_j)\geq \Delta^{\infty,(r)}\Big(\ov\gamma\Big(\frac{i-1}{n}\Big),\ov\gamma\Big(\frac{i}{n}\Big)\Big).
\end{equation}
For every $1\leq j\leq p$, the properties $\Lambda_{u_j}>r+\delta/2$ and $\Delta^{p,\circ}(u_{j-1},u_j)<\delta/2$ 
imply that the minimal 
label on $[u_{j-1},u_j]$ is greater than $r$ (or the same holds with $[u_{j-1},u_j]$ replaced by $[u_j,u_{j-1}]$).
This shows in particular that there is a continuous curve from $\gamma(\frac{i}{n})$ to $\Pi^p(u_j)$ that stays 
in the complement of $B(r)$, so that $\Pi^p(u_j)$ belongs to $\check B^\bullet(r)$ and $u_j$
must belong to $F_r$, which allows us to define $\ov u_j$ as the point of $\t^i_\infty$ corresponding to $u_j$. Furthermore the fact that the minimal 
label on $[u_{j-1},u_j]$ is greater than $r$ also implies that the interval $[u_{j-1},u_j]$ in $\t^p_\infty$
is identified to the interval $[\ov u_{j-1},\ov u_j]$ in $\t^i_\infty$, and then that
 $\Delta^{p,\circ}(u_{j-1},u_j)=\Delta^{\infty,(r),\circ}(\ov u_{j-1},\ov u_j)$. 
 The bound \eqref{tech-iBdBp} follows, which completes the proof
of the lemma and of Theorem \ref{iBdBp}. \endproof

In view of applications to isoperimetric inequalities in the Brownian plane \cite{Rie}, we state
another result which complements Theorem \ref{infBdBp} by showing that, in some sense, the exterior of
the hull $B^\bullet(r)$ is independent of this hull, conditionally on its boundary size. We keep the notation introduced at the
beginning of this section, and in particular, we recall that the Brownian plane $\BP_\infty$
is constructed from the labeled tree $(\t^p_\infty,(\Lambda_v)_{v\in\t^p_\infty})$ associated with 
the coding triple $(X,\ll,\rr)$. We fix $r>0$ and write $K_r$ for the complement of the set
$F_r$ defined in \eqref{hull-expl},
$$K_r:=\{v\in\t^p_\infty: \Lambda_w\leq r\hbox{ for some }w\in \llbracket v,\infty\llbracket\}.$$
 We have then $B^\bullet(r)=\Pi^p(K_r)$ (cf.~formulas (16) and (17) in \cite{CLG}). Recall that, for every 
$u,v\in \BP_\infty$, 
\begin{equation}
\label{pseudo-hull1}
\Delta^{p,\circ}(u,v)=\Lambda_u + \Lambda_v -2 \max\Big( \inf_{w\in[u,v]} \Lambda_w,\inf_{w\in[v,u]} \Lambda_w\Big).
\end{equation}
We then set, for every $u,v\in K_r$,
\begin{equation}
\label{pseudo-hull2}
\Delta^{p,(r)}(u,v) = \build{\inf_{u_0=u,u_1,\ldots,u_p=v}}_{u_0,u_1,\ldots,u_p\in K_r}^{} \sum_{i=1}^p \Delta^{p,\circ}(u_{i-1},u_i)
\end{equation}
where the infimum is over all choices of the integer $p\geq 1$ and of the
finite sequence $u_0,u_1,\ldots,u_p$ in $K_r$ such that $u_0=u$ and
$u_p=v$. For every $u,v\in K_r$, we have $\Delta^p(u,v)\leq \Delta^{p,(r)}(u,v)$ (just note that $\Delta^p(u,v)$ is defined by the
same formula \eqref{pseudo-hull2} without the restriction to $u_0,\ldots,u_p\in  K_r$) and we also know 
that $\Delta^p(u,v)=0$  implies $\Delta^{p,\circ}(u,v)=0$ and a fortiori $\Delta^{p,(r)}(u,v)=0$. It follows that $\Delta^{p,(r)}$
induces a metric on $\Pi^p(K_r)=B^\bullet(r)$, and we keep the notation $\Delta^{p,(r)}$ for this metric. 
For future use, we also observe that, in the right-hand side of formula \eqref{pseudo-hull1}
applied to $u,v\in K_r$, we may replace the 
infimum over $w\in[u,v]$ by an infimum over $w\in[u,v]\cap K_r$: The point is that, if the clockwise exploration
going from $u$ to $v$ (or from $v$ to $u$) intersects $F_r$, then it necessarily visits a point with label 
at most $r$, because otherwise $u$ and $v$ would have to be in $F_r$.

\begin{theorem}
\label{indep-hull}
Conditionally on $|\partial B^\bullet(r)|$, the random compact measure metric space $(B^\bullet(r),\Delta^{p,(r)})$
and the space $(\D^{\infty,(r)},\Delta^{\infty,(r)})$ in Theorem \ref{infBdBp} are independent. Furthermore, the restriction of
the metric $\Delta^{p,(r)}$ to the interior of $B^\bullet(r)$ coincides with the intrinsic metric induced by $\Delta^p$
on this open set.
\end{theorem}

\proof  The general idea is to show that the space $(B^\bullet(r),\Delta^{p,(r)})$ can be constructed from random quantities 
that are independent of $(\D^{\infty,(r)},\Delta^{\infty,(r)})$ conditionally on $|\partial B^\bullet(r)|$.
We start by introducing the labeled tree $\t^{p,(r)}$ which consists of the part $[0,\mathrm{L}_r]$
of the spine of
$\t^p_\infty$, and of the subtrees branching off $[0,\mathrm{L}_r]$. Equivalently, $\t^{p,(r)}$ is associated with the finite spine coding triple
\begin{equation}
\label{trunc-coding}
\mathfrak{X}^{(r)}:=\Big((X_t)_{0\leq t\leq \mathrm{L}_r}, \mathbf{1}_{[0,\mathrm{L}_r]}(t)\,\ll(\dd t\dd\omega),  \mathbf{1}_{[0,\mathrm{L}_r]}(t)\,\rr(\dd t\dd\omega)\Big).
\end{equation}
Clearly, $\t^{p,(r)}$ viewed as a subset of $\t^p_\infty$ is contained in $K_r$. If $(\ee_s)_{s\in\R}$
is the clockwise exploration of $\t^p_\infty$, $\t^{p,(r)}$ corresponds to the points visited by
$(\ee_s)_{s\in\R}$ during an interval of the form $[-\sigma_r,\sigma'_r]$, with $\sigma_r,\sigma'_r>0$. 
We also let $\check{\mathfrak{X}}^{(r)}$ be the image of $\mathfrak{X}^{(r)}$ under the time-reversal operation
\eqref{time-rev} and denote the associated labeled tree by $\check\t^{p,(r)}$ (replacing 
$\t^{p,(r)}$ by $\check\t^{p,(r)}$ just amounts to interchanging the roles of the root and the top of the spine).

We then consider all subtrees branching off the spine of $\t^p_\infty$ at a level higher than $\mathrm{L}_r$, and, for each such 
subtree whose labels hit $[0,r]$, the ``excursions outside'' $(r,\infty)$. To make this precise, write
$$\ll=\sum_{i\in I} \delta_{(t_i,\omega^i)}\;,\ \rr=\sum_{i\in J} \delta_{(t_i,\omega^i)},$$
where the indexing sets $I$ and $J$ are disjoint. In the time scale of the clockwise exploration, each $\omega^i$
corresponds to an interval $[\alpha_i,\beta_i]$ contained in $(-\infty,0)$ if $i\in J$, or in $(0,\infty)$ if $i\in I$, and
$\sigma(\omega^i)=\beta_i-\alpha_i$. 
Set $I_r:=\{i\in I:t_i>\mathrm{L}_r\hbox{ and }W_*(\omega^i)\leq r\}$ and
$J_r:=\{i\in J:t_i>\mathrm{L}_r\hbox{ and }W_*(\omega^i)\leq r\}$. For each $i\in I_r\cup J_r$, we can make sense of the 
exit local time of $\omega^i$ at level $r$, as defined in Section \ref{sna-mea}, and we 
denote this local time by $(L^{i,r}_t)_{t\in[0,\sigma(\omega^i)]}$. We then set, for every $t\in \R$,
$$L^{*,r}_t= \sum_{i\in I\cup J} L^{i,r}_{t\wedge \beta_i - t\wedge \alpha_i},$$
so that, in some sense, $L^{*,r}_t$ represents the total exit local time accumulated at $r$
by the clockwise exploration up to time $t$. We note that
$$L^{*,r}_\infty= \sum_{i\in I\cup J} L^{i,r}_{\sigma(\omega^i)}=\z^{(r)},$$
and $|\partial B^\bullet(r)|=\z^{(r)}$ by definition.

Then, for every $i\in I_r\cup J_r$, we consider the excursions $(\omega^{i,k})_{k\in\N}$
of $\omega^i$ outside $(r,\infty)$ (we refer to \cite[Section 2.4]{ALG} for more information about
such excursions). These excursions $\omega^{i,k}$, $k\in\N$ are in one-to-one correspondence with the
connected components $(a_{i,k},b_{i,k})$, $k\in \N$, of the open set $\{s\in[0,\sigma(\omega^i)]:\tau_r(\omega^i_s)<\zeta_s(\omega^i)\}$, in such a way that, for every $s\geq 0$, 
$$\omega^{i,k}_s(t):=\omega^i_{(a_{i,k}+s)\wedge b_{i,k}}(\zeta_{a_{i,k}}(\omega^i)+t),\hbox{ for }
0\leq t\leq \zeta_s(\omega^{i,k}):=\zeta_{(a_{i,k}+s)\wedge b_{i,k}}(\omega^i)-\zeta_{a_{i,k}}(\omega^i).$$
In the time scale of the clockwise exploration, $\omega^{i,k}$ corresponds to the interval $[\alpha_{i,k},\beta_{i,k}]$,
where $\alpha_{i,k}=\alpha_i+a_{i,k}$ and $\beta_{i,k}=\alpha_i + b_{i,k}$. In particular,
 the (labeled) tree $\t_{(\omega^{i,k})}$ coincides with the subtree of $\t^p_\infty$ consisting
of the descendants of $\ee_{\alpha_{i,k}}=\ee_{\beta_{i,k}}$ (this set of descendants is $\{\ee_s:s\in[\alpha_{i,k},\beta_{i,k}]\}$). 

Recall the coding triple $(X^{(r)},\wt\ll^{(r)},\wt\rr^{(r)})$ which is used to construct the space
$(\D^{\infty,(r)},\Delta^{\infty,(r)})$. An application of the special Markov property, in the form 
given in the appendix of \cite{subor}, shows that, conditionally on $\z^{(r)}$, the 
point measure
$$\nn_{(r)}(\dd t\dd \omega):=\sum_{i\in I_r\cup J_r}\sum_{k\in\N} \delta_{(L^{*,r}_{\alpha_{i,k}},\omega^{i,k})}(\dd t\dd\omega)$$
is Poisson with intensity $\mathbf{1}_{[0,\z^{(r)}]}(t)\,\dd t\,\N_r(\dd\omega)$, and is independent of 
$(X^{(r)},\wt\ll^{(r)},\wt\rr^{(r)})$. Note in particular that $\z^{(r)}$ is a measurable function 
of $\nn_{(r)}$. On the other hand, the coding triple $\mathfrak{X}^{(r)}$ in \eqref{trunc-coding} is
clearly independent of the pair $(\nn_{(r)},(X^{(r)},\wt\ll^{(r)},\wt\rr^{(r)}))$. So the first assertion of the theorem
would follow if we could prove that the space $(B^\bullet(r),\Delta^{p,(r)})$ is a function of $\nn_{(r)}$ and $\mathfrak{X}^{(r)}$. This is not correct, but we will see that $(B^\bullet(r),\Delta^{p,(r)})$ is a
function of $(\nn_{(r)},L^{*,r}_0,\mathfrak{X}^{(r)})$, whose conditional distribution given 
$(\z^{(r)},L^{*,r}_0)$ only depends on $\z^{(r)}$. This will suffice to get the first assertion of 
Theorem \ref{indep-hull}. 

Let us explain how the space $(B^\bullet(r),\Delta^{p,(r)})$ can be written as a function of 
$(\nn_{(r)},L^{*,r}_0)$ and $\mathfrak{X}^{(r)}$. To begin with, we introduce the right-continuous inverse
of the process $(L^{*,r}_t)_{t\in\R}$: for every $s\in[0,\z^{(r)})$,
$$\tau^{*,r}_s:=\inf\{t\in\R: L^{*,r}_t>s\},$$
and we also make the convention that $\tau^{*,r}_{\z^{(r)}}$ is the left limit of 
$s\mapsto \tau^{*,r}_s$ at $s= \z^{(r)}$. Then one verifies that $(\Pi^p(\ee_{\tau^{*,r}_s}))_{0\leq t\leq \z^{(r)}}$
is an injective loop whose range is precisely $\partial B^\bullet(r)$. Let us briefly justify this. Recall that
$\partial B^\bullet(r)= \Pi^p(\partial F_r)$, with
$\partial F_r$ given by 
\eqref{bdry-expl}. We first observe that the mapping $s\mapsto \Pi^p(\ee_{\tau^{*,r}_s})$ is continuous. Indeed,
we already know that the function $s\mapsto \Pi^p(\ee_s)$ is continuous. Furthermore,
if $\tau^{*,r}_{s-}< \tau^{*,r}_s$, the support property of the exit local time implies that either 
all points of the form $\ee_u$ with $u\in(\tau^{*,r}_{s-}, \tau^{*,r}_s)$ 
 are descendants of $\ee_{\tau^{*,r}_{s-}}$ and necessarily $\ee_{\tau^{*,r}_{s-}}=\ee_{ \tau^{*,r}_s}$, or the labels of all such points $\ee_u$ are greater than $r$.
In both cases, we have $\Pi^p(\ee_{\tau^{*,r}_{s-}})
=\Pi^p(\ee_{\tau^{*,r}_s})$. 
Then one easily deduces from the same support property that any point of
the form $\Pi^p(\ee_{\tau^{*,r}_s})$ belongs to $\partial B^\bullet(r)$. Conversely, using \eqref{bdry-expl},
any point $x$ of $\partial B^\bullet(r)$, with the exception of the point $\mathrm{L}_r$ of the spine, must be of the form $\Pi^p(v)$ where $v$ 
belongs to a subtree $\t_{(\omega^i)}$ with $i\in I\cup J$, and labels along the line segment between $v$ 
and the root of $\t_{(\omega^i)}$ are greater then $r$ except at $v$. From the support property of the exit local
time, it follows that $v=\ee_{\tau^{*,r}_s}$ for some $s\in[\alpha_i,\beta_i]$.
The formula $v=\ee_{\tau^{*,r}_s}$ also holds for $v=\mathrm{L}_r$ with $s=L^{*,r}_0$.  Finally, from the description of
the distribution of $\nn_{(r)}$ and the fact that $\Delta^p(u,v)=0$ holds if and only if 
$\Delta^{p,\circ}(u,v)=0$, one checks that, for every $0\leq s<s'\leq \z^{(r)}$, the points $\Pi^p(\ee_{\tau^{*,r}_s})$
and $\Pi^p(\ee_{\tau^{*,r}_{s'}})$ are distinct, except in the case $s=0$ and $s'=\z^{(r)}$. 

We let $\mathfrak{H}$ be derived from the disjoint union 
$$[0,\z^{(r)}] \cup \Big(\build{\bigcup_{i\in I_r\cup J_r}}_{k\in\N}^{} \t_{(\omega^{i,k})}\Big) \cup \check\t^{p,(r)}$$
by identifying $0$ with $\z^{(r)}$, the root of $\check\t^{p,(r)}$ with the point $L^{*,r}_0$ of $[0,\z^{(r)}]$, and,
for every $i\in I_r\cup J_r$ and $k\in\N$, the root of $\t_{(\omega^{i,k})}$ with the point $L^{*,r}_{\alpha_{i,k}}$ of $[0,\z^{(r)}]$.
We assign labels $(\Lambda^{(r)}_x)_{x\in \mathfrak{H}}$ to the points of $\mathfrak{H}$: the label of
any point of $[0,\z^{(r)}]$ is equal to $r$, and points of the labeled trees $\t_{(\omega^{i,k})}$ and $\check\t^{p,(r)}$
keep their labels. We also define a volume measure on $\mathfrak{H}$ by summing the volume measures
of the trees $\t_{(\omega^{i,k})}$ and of $\check\t^{p,(r)}$. The total volume of $\mathfrak{H}$ is
$$\Sigma^{(r)}:=|\check\t^{p,(r)}|+\sum_{i\in I_r\cup J_r}\sum_{k\in\N} \sigma(\omega^{i,k}),$$
using the notation $|\check\t^{p,(r)}|$ for the total volume of $\check\t^{p,(r)}$.

We need to define a cyclic clockwise exploration of $\mathfrak{H}$, which will be denoted
by $(\ee^{(r)}_s)_{s\in[0,\Sigma^{(r)}]}$. Roughly speaking this exploration corresponds
to concatenating the clockwise explorations of the trees $\t_{(\omega^{i,k})}$ and $\check\t^{p,(r)}$ in the order 
prescribed by the exploration of $\t^p_\infty$. To give a more precise definition, we first observe that we can 
write $K_r={K^\circ_r} \cup \partial F_r$, where $\partial F_r$ is as in \eqref{bdry-expl}, and 
$$K^\circ_r:=\{v\in \t^p_\infty : \Lambda_w\leq r\hbox{ for some }w\in\rrbracket v,\infty \llbracket\}.$$
If $v\in \partial F_r$, we know that $\Pi^p(v)\in \partial B^\bullet(r)$, so that,  by previous observations, 
there is a unique $s\in[0,\z^{(r)})$ such that $v=\Pi^p(\ee_{\tau^{*,r}_{s}})$, and we set $\Phi_{(r)}(v):=s$.

We then define, for every $s\in\R$,
$$A^{(r)}_s:=\int_{-\infty}^s \dd t\,\mathbf{1}_{\{\ee_t\in K_r\}}.$$
Note that $A^{(r)}_\infty=\Sigma^{(r)}$. We set $\eta^{(r)}_t:=\inf\{s\in\R: A^{(r)}_s>t\}$ for every $t\in[0,\Sigma^{(r)})$. 
Then, for every $t\in[0,\Sigma^{(r)})$, either $\ee_{\eta^{(r)}_t}$ belongs to $K_r^\circ$, which implies that
$\ee_{\eta^{(r)}_t}$ is a point of one of the trees $\t_{(\omega^{i,k})}$ or of $\check\t^{p,(r)}$, and we let 
$\ee^{(r)}_t$ be the ``same'' point in $\mathfrak{H}$, or $\ee_{\eta^{(r)}_t}$ belongs to $\partial F_r$, and
we set $\ee^{(r)}_t=\Phi_{(r)}(\ee_{\eta^{(r)}_t})\in [0,\z^{(r)})$. Finally we take $\ee^{(r)}_{\Sigma^{(r)}}=\ee^{(r)}_0=0$. Although this is not apparent in the preceding
presentation, the reader will easily check that this exploration process $\ee^{(r)}$ 
only depends on $(\nn_{(r)},L^{*,r}_0)$ and $\mathfrak{X}^{(r)}$ (the reason why we need $L^{*,r}_0$
is because we have to rank the tree $\t^{p,(r)}$ among the trees $\t_{(\omega^{i,k})}$ --- of course the
order between the different trees $\t_{(\omega^{i,k})}$ is prescribed by the point measure $\nn_{(r)}$). 

The clockwise exploration of $\mathfrak{H}$ allows us to make sense of intervals on $\mathfrak{H}$.
In turn, we can then define the function $D'^{(r),\circ}(u,v)$, for $u,v\in \mathfrak{H}$, by the
right-hand side of \eqref{pseudo-hull1}, where we simply replace $\Lambda$ by $\Lambda^{(r)}$. 
Similarly, we define $D'^{(r)}(u,v)$, for $u,v\in \mathfrak{H}$, by replacing $\Delta^{p,\circ}$ with $D'^{(r),\circ}$ in
the
right-hand side of \eqref{pseudo-hull2} (and of course replacing $u_0,u_1,\ldots,u_p\in K_r$
by $u_0,u_1,\ldots,u_p\in \mathfrak{H}$). We now claim that the quotient space 
$\mathfrak{H}/\{D'^{(r)}=0\}$, equipped with the metric induced by $D'^{(r)}$ and with the volume 
measure which is the pushforward of the volume measure on $\mathfrak{H}$, coincides 
with $(B^\bullet(r),\Delta^{p,(r)})$. This is a straightforward consequence of our construction
(using the fact that one can replace $[u,v]$ by $[u,v]\cap K_r$
in the right-hand side of \eqref{pseudo-hull1} when $u,v\in K_r$), and we
omit the details. 

We also observe that the conditional distribution of the space $\mathfrak{H}/\{D'^{(r)}=0\}$ given 
$(\z^{(r)},L^{*,r}_0)$ does not depend on $L^{*,r}_0$. This follows from the fact that the law 
of a Poisson point measure on $[0,z]\times \S$ with intensity $\dd t\,\N_{r}(\dd \omega)$
is invariant under the shift $t\mapsto t+a\hbox{ mod.}\; z$, for any fixed $a\in[0,z]$. 

Finally, we can write $(B^\bullet(r),\Delta^{p,(r)})=(\mathfrak{H}/\{D'^{(r)}=0\},D'^{(r)})=\Psi(\nn_{(r)},L^{*,r}_0,\mathfrak{X}^{(r)})$
with a $\mathbb{K}$-valued function $\Psi$, and we have for every nonnegative
measurable function $F$ on $\mathbb{K}$,
\begin{align*}\E[F(\Psi(\nn_{(r)},L^{*,r}_0,\mathfrak{X}^{(r)}))\mid (X^{(r)},\wt\ll^{(r)},\wt\rr^{(r)})]&=\E[F(\Psi(\nn_{(r)},L^{*,r}_0,\mathfrak{X}^{(r)}))\mid (\z^{(r)},L^{*,r}_0)]\\
&=\E[F(\Psi(\nn_{(r)},L^{*,r}_0,\mathfrak{X}^{(r)}))\mid \z^{(r)}].
\end{align*}
The second equality follows from the preceding observation, and the first one holds because $\nn_{(r)}$ and $(X^{(r)},\wt\ll^{(r)},\wt\rr^{(r)})$
are conditionally independent given $\z^{(r)}$ (and $L^{*,r}_0$ is a measurable function of $(X^{(r)},\wt\ll^{(r)},\wt\rr^{(r)})$). Since 
$(\D^{\infty,(r)},\Delta^{\infty,(r)})$ is a function of $(X^{(r)},\wt\ll^{(r)},\wt\rr^{(r)})$, this gives the first
% We have thus proved that $(B^\bullet(r),\Delta^{p,(r)})$ is 
% a function of $(\nn_{(r)},L^{*,r}_0,\mathfrak{X}^{(r)})$ whose conditional distribution 
%given $(\z^{(r)},L^{*,r}_0)$ only depends on $\z^{(r)}$. As we already explained, this gives the first
assertion of Theorem \ref{indep-hull}. 

The proof of the second assertion of Theorem \ref{indep-hull} is very similar to the
proof of Lemma \ref{compar-dist}, and we leave the details to the reader. \endproof

\rem The preceding proof gives a description of the distribution of the hull $B^\bullet(r)$ equipped 
with its intrinsic metric in terms of the space $\mathfrak{H}$. We note that the labeled tree 
$\check\t^{p,(r)}$ has the same distribution as the tree $\t_{(\omega)}$ under
$\N_r(\dd\omega\,|\, W_*(\omega)=0)$ (see \cite{Bessel}). So the conditional distribution of $B^\bullet(r)$
knowing $\z^{(r)}=z$ could as
well be defined from a Poisson point measure $\sum_{i\in I} \delta_{(t_i,\omega_i)}$ with intensity $\mathbf{1}_{[0,z]}(t)\,\dd t\,\N_r(\dd\omega)$ conditioned on the event $\{\inf_{i\in I}W_*(\omega_i) = 0\}$. In this form, there is a striking analogy
with the construction of the (free) Brownian disk with perimeter $z$ found in \cite{Bet} or \cite{BM} --- see Section
\ref{sec:consist} below for a presentation within the formalism of the present work. The essential difference
comes from the fact that the construction of the hull assigns constant labels equal to $r$ to points of $\mathfrak{H}$ that belong to $[0,z]$, whereas,
in the construction of the Brownian disk, labels along $[0,z]$ evolve like a scaled Brownian bridge.

\subsection{Horohulls in the Brownian plane}
\label{sec:horo}

In this section, we explain how pointed Brownian disks with a given height appear as horohulls 
in the Brownian plane. Let us first recall the definition of these horohulls. We consider the
Brownian plane $(\BP_\infty,\Delta^p)$, with the distinguished point $\rho$. One can prove
\cite{Plane} that, a.s. for every $a,b\in\BP_\infty$, the limit
$$\lim_{x\to\infty} (\Delta^p(a,x)-\Delta^p(b,x))$$
exists in $\R$. Here the limit as $x\to\infty$ means that $x$ tends to the point at infinity in the
Alexandroff compactification of $\BP_\infty$. Clearly, the limit in the preceding display can be
written in the form $\mathcal{H}_a-\mathcal{H}_b$, where the ``horofunction'' $a\mapsto \mathcal{H}_a$ is uniquely defined
if we impose $\mathcal{H}_\rho=0$. We interpret $\mathcal{H}_a$ as a (relative) distance from $a$ to $\infty$, and call
$\mathcal{H}_a$ the horodistance from $a$. Note the bound $|\mathcal{H}_a-\mathcal{H}_b|\leq \Delta^p(a,b)$. 

For every $r>0$, let $\mathfrak{B}^\circ(r)$ be the connected component of the open set $\{x\in\BP_\infty:\mathcal{H}_x>-r\}$
that contains $\rho$. So a point $x$ belongs to $\mathfrak{B}^\circ(r)$ if and only if there is a continuous path 
from $\rho$ to $r$ that stays at horodistance greater than $-r$. The horohull $\mathfrak{B}^\bullet(r)$ is defined 
as the closure of $\mathfrak{B}^\circ(r)$. We view $\mathfrak{B}^\bullet(r)$ as a pointed compact measure metric
space with distinguished point $\rho$. Note that the compactness of $\mathfrak{B}^\bullet(r)$ is not obvious
a priori, but will follow from the description that we give in the proof of the next statements.

We write $\mathrm{Vol}(\cdot)$ for the volume measure on $\BP_\infty$.
In the following two statements, we fix $r>0$.

\begin{proposition}
\label{bdry-horo}
The limit 
$$\lim_{\ve\downarrow 0} \ve^{-2}\,\mathrm{Vol}(\{x\in \mathfrak{B}^\bullet(r):\mathcal{H}_x<-r+\ve\})$$
exists a.s. This limit is called the boundary size of $\mathfrak{B}^\bullet(r)$ and 
denoted by $|\partial \mathfrak{B}^\bullet(r)|$.
\end{proposition}

\begin{theorem}
\label{law-horo}
The intrinsic metric on $\mathfrak{B}^\circ(r)$ has a.s. a continuous extension to $\mathfrak{B}^\bullet(r)$,
which is denoted by $\Delta^{horo, r}_\infty$. Then conditionally on $|\partial \mathfrak{B}^\bullet(r)|=z$, the pointed
measure metric space $(\mathfrak{B}^\bullet(r),\Delta^{horo, r}_\infty)$ is a pointed Brownian disk 
with perimeter $z$ and height $r$.
\end{theorem}

The proof of both Proposition \ref{bdry-horo} and Theorem \ref{law-horo} relies 
on the construction of the Brownian plane found in \cite{Plane}, which is 
different from the one given in Section \ref{Br-plane}. Let us recall the
construction of \cite{Plane} using our formalism of coding triples
(the presentation therefore seems to differ from the one in \cite{Plane}, but
the relevant random objects are the same). We consider a coding triple
$(B,\pp,\pp')$, such that:
\begin{itemize}
\item[$\bullet$]
$B=(B_t)_{t\geq 0}$ is a linear Brownian motion started from $0$.
\item[$\bullet$]
Conditionally on $B$, $\pp$ and $\pp'$ are independent Poisson point measures on $\R_+\times \S$
with intensity
$2\,\dd t\,\N_{B_t}(\dd \omega)$.
\end{itemize}
Following Section \ref{coding-infinite}, we then consider the infinite labeled
tree $(\t'^p_\infty,(\Lambda'_v)_{v\in\t'^p_\infty})$ associated with this 
coding triple. We define the functions $D^{\infty,\circ}(u,v)$
and $D^\infty(u,v)$, for $u,v\in\t'^p_\infty$, in a way similar to Section \ref{sec:pseudo}
(note however that labels are here of arbitrary sign):
\begin{equation}
\label{pseudo-BP1}
D^{\infty,\circ}(u,v)=\Lambda'_u + \Lambda'_v -2 \max\Big( \inf_{w\in[u,v]} \Lambda'_w,\inf_{w\in[v,u]} \Lambda'_w\Big),
\end{equation}
and
\begin{equation}
\label{pseudo-BP2}
D^\infty(u,v) = \inf_{u_0=u,u_1,\ldots,u_p=v} \sum_{i=1}^p D^{\infty,\circ}(u_{i-1},u_i)
\end{equation}
where the infimum is over all choices of the integer $p\geq 1$ and of the
finite sequence $u_0,u_1,\ldots,u_p$ in $\t'^p_\infty$ such that $u_0=u$ and
$u_p=v$. We have $D^\infty(u,v)=0$
if and only if $D^{\infty,\circ}(u,v)=0$ \cite[Proposition 11]{Plane}.

We let $\BP'_\infty$ be the quotient space $\t'^p_\infty/\{D^\infty=0\}$,
which is equipped with the distance induced by $D^\infty(u,v)$ and the volume 
measure which is the pushforward of the volume measure on $\t'^p_\infty$,
and with the distinguished point which is the equivalence class of 
the root of $\t'^p_\infty$. We also let
$\Pi'^p$ stand for the canonical projection from $\t'^p_\infty$
onto $\BP'_\infty$. 

Then the pointed measure metric space $\BP'_\infty$ is a Brownian plane,
that is, it
has the same distribution as $\BP_\infty$ (see \cite[Theorem 3.4]{CLG}).
Therefore, we can replace $\BP_\infty$ by $\BP'_\infty$ in the proof
of both Proposition \ref{bdry-horo} and Theorem \ref{law-horo}. The point of this replacement
is the fact that the horodistance from a point $a$ of $\BP'_\infty$
is now equal to its label $\Lambda'_a$ \cite[Proposition 17]{CLG}. Indeed, we can
summarize the difference between the two constructions of the Brownian plane
by saying that labels correspond to distances from the distinguished point 
in the first construction, and to horodistances in the second one. In the proofs below, we
assume that $\mathfrak{B}^\circ(r)$ and $\mathfrak{B}^\bullet(r)$ are defined
in $\BP'_\infty$, and without risk of confusion we use the notation $\rho$
both for the root of $\t'^p_\infty$ and for the distinguished point of $\BP'_\infty$.

\medskip
\noindent{\it Proof of Proposition \ref{bdry-horo} and Theorem \ref{law-horo}.} The first step is to observe 
that we have $\mathfrak{B}^\circ(r)=\Pi'^p(G_r)$, where
\begin{equation}
\label{horo-expl}
G_r:=\{v\in\t'^p_\infty: \Lambda'_w>-r\hbox{ for every }w\in \llbracket \rho,v\rrbracket\},
\end{equation}
and $\mathfrak{B}^\bullet(r)=\mathfrak{B}^\circ(r)\cup \partial\mathfrak{B}^\circ(r)$, 
with $\partial\mathfrak{B}^\circ(r)=\Pi'^p(\partial G_r)$, and 
\begin{equation}
\label{horo-bdry}
\partial G_r=\{v\in\t'^p_\infty: \Lambda'_v=-r\hbox{ and }\Lambda'_w>-r\hbox{ for every }w\in \llbracket \rho,v\llbracket\}.
\end{equation}
Notice the similarity with \eqref{hull-expl} and \eqref{bdry-expl}. Let us justify the equality 
$\mathfrak{B}^\circ(r)=\Pi'^p(G_r)$. The inclusion $\mathfrak{B}^\circ(r)\supset\Pi'^p(G_r)$
is easy, because, if $v\in G_r$, the image under $\Pi'^p$ of the geodesic segment 
$\llbracket\rho,v\rrbracket$ yields a continuous path from $v$ to $\rho$ along which labels
(horodistances) stay greater than $-r$. The reverse inclusion comes from the so-called 
``cactus bound'' which says than any continuous path between $\rho$ and $\Pi'^p(v)$ must 
visit a point whose label is smaller than or equal to $\min_{u\in\llbracket\rho,v\rrbracket}\Lambda'_u$
(see formula (4) in \cite{Plane} for a short proof
in the case of the Brownian map, which is immediately extended to the present setting). Once the equality $\mathfrak{B}^\circ(r)=\Pi'^p(G_r)$ is established,
the property $\partial\mathfrak{B}^\circ(r)=\Pi'^p(\partial G_r)$ is easy and we omit the details.

Write $\mathrm{cl}(G_r)=G_r\cup \partial G_r$, which we can view as a (compact) subtree of
the tree $\t'^p_\infty$. In a way very similar to the proof of Theorem \ref{infBdBp}, we may
interpret $\mathrm{cl}(G_r)$ as the (labeled) tree associated with a coding triple derived from
the triple $(B,\pp,\pp')$. To this end, we set
$$T_r:=\inf\{t\geq 0: B_t=-r\},$$
and we note that $\mathrm{cl}(G_r)$ consists of the union of the part
$[0,T_r]$ of the spine of $\t'^p_\infty$ with the subtrees branching off the spine
between levels $0$ and $T_r$ and truncated at label $-r$. To make this more precise,
if
$$\pp=\sum_{i\in I} \delta_{(t_i,\omega_i)},$$
we define
$$\pp^{(r)}=\sum_{i\in I, t_i<T_r} \delta_{(t_i,\tr_{-r}(\omega_i))}$$
and we similarly define $\pp'^{(r)}$ from $\pp'$. Let
$B^{(r)}$ stand for the stopped path $(B_t)_{0\leq t\leq T_r}$. Then $\mathrm{cl}(G_r)$
is canonically and isometrically identified with the (labeled) tree coded by
the triple $(B^{(r)},\pp^{(r)},\pp'^{(r)})$. This identification preserves the labels and
the volume measures. The fact that the limit in Proposition \ref{bdry-horo} exists,
and is in fact given by
$$|\partial\mathfrak{B}^\bullet(r)|= \int \z_{-r}(\omega)\,\pp^{(r)}(\dd\omega)
+ \int \z_{-r}(\omega)\,\pp'^{(r)}(\dd\omega)$$
now follows from the approximation formula \eqref{formu-exit} for exit measures, using also Proposition \ref{Palm}
and \eqref{tech-exit1}.

Recall the notation $\vartheta_r$ introduced before Proposition \ref{ident-cond-dis}.
In order to derive the statement of Theorem \ref{law-horo}, we now notice that
the triple $\mathfrak{T}^{(r)}:=(B^{(r)}+r,\vartheta_{-r}\pp^{(r)},\vartheta_{-r}\pp'^{(r)})$ has the same 
distribution as the coding triple $(Y,\wt\mm,\wt\mm')$ considered at the beginning of Section
\ref{sec-ident}, provided we take $a=r$. It follows that the conditional distribution 
of  $\mathfrak{T}^{(r)}$ knowing that 
$|\partial\mathfrak{B}^\bullet(r)|=z$ is $\check\Theta^r_z$. Recall the mapping
$\Omega$ defined in Section \ref{coding-infinite}.  Then $\Omega(\mathfrak{T}^{(r)})$
is a random snake trajectory which, conditionally on $|\partial\mathfrak{B}^\bullet(r)|=z$, is
distributed according to $\N^{(z)}_r$. Furthermore, by Proposition \ref{pointed-disk}, 
the random metric space
$\D^{\bullet,r}:=\Xi^\bullet(\Omega(\mathfrak{T}^{(r)}))$
is a pointed Brownian disk with perimeter $z$ and height $r$, conditionally on $|\partial\mathfrak{B}^\bullet(r)|=z$. 
To complete the proof, we just need to identify $\mathfrak{B}^\bullet(r)$ (equipped with the intrinsic distance) with
$\D^{\bullet,r}$.

By a preceding observation, $\mathrm{cl}(G_r)$ is identified to the genealogical tree of 
$\Omega(\mathfrak{T}^{(r)})$ (which is the labeled tree associated with $\mathfrak{T}^{(r)}$) and this identification
preserves labels, provided labels on $\mathrm{cl}(G_r)$ are shifted by $r$. One then verifies that
two points of $\mathrm{cl}(G_r)$ are identified in $\Pi'^p(\mathrm{cl}(G_r))=\mathfrak{B}^\bullet(r)$
if and only if the corresponding points of the genealogical tree of $\Omega(\mathfrak{T}^{(r)})$ are identified in $\D^{\bullet,r}$.
It follows that $\mathfrak{B}^\bullet(r)$ and $\D^{\bullet,r}$ can be identified as sets.
To complete the proof of Theorem \ref{law-horo}, it then remains to show that 
the intrinsic distance between two points of 
$\mathfrak{B}^\circ(r)$ coincides with the distance between the corresponding points 
of the interior of $\D^{\bullet,r}$ (from the discussion in Section \ref{sec:point-disk}, this will imply
first that the intrinsic distance on $\mathfrak{B}^\circ(r)$ can be extended to the boundary, and then
that $\mathfrak{B}^\bullet(r)$ is isometric to $\D^{\bullet,r}$ as desired). This is derived 
by arguments very similar to the end of the proof of Theorem \ref{infBdBp}, and we omit
the details. 
\hfill$\square$

\medskip
We conclude this section with some explicit distributional properties of the process of horohulls.
It will be convenient to use the Skorokhod space $\D(\R_+,\R)$ of c\`adl\`ag functions 
from $\R_+$ into $\R$. We write $(Z_t)_{t\geq 0}$ for the canonical process 
on $\D(\R_+,\R)$, and $(\f_t)_{t\geq 0}$ for the canonical filtration. We then introduce, for every $x\geq 0$, the probability measure $P_x$
which is the law of the continuous-state branching process with branching mechanism $\Phi$ 
(in short, the $\Phi$-CSBP), where $\Phi(\lambda)=\sqrt{8/3}\,\lambda^{3/2}$. 
We refer to \cite[Section 2.1]{CLG} for the definition and some properties of the $\Phi$-CSBP. 

The $\Phi$-CSBP is critical, meaning that $E_x[Z_t]=x$ for every $t\geq 0$ and $x\geq 0$. 
Then, for every $x>0$, we can define the law $P^\uparrow_x$ of the $\Phi$-CSBP started
from $x$ and conditioned to non-extinction via the $h$-transform
\begin{equation}
\label{htrans}
{\frac{\dd P^\uparrow_x}{\dd P_x}}_{\Big|\f_t} = \frac{Z_t}{x}.
\end{equation}
See \cite[Section 4.1]{Lam} for a discussion of continuous-state branching
processes conditioned on non-extinction. The preceding formula 
does not make sense for $x=0$. However, \cite[Theorem 2]{KP} shows that 
the laws $P^\uparrow_x$ converge weakly as $x\downarrow 0$ to a limiting
law denoted by $P^\uparrow_0$, which is characterized by the following
two properties:
\begin{itemize}
\item[(i)] for every $t>0$, the law of $Z_t$ under $P^\uparrow_0$ is given by
$$E^\uparrow_0[e^{-\lambda Z_t}]= \Big(1+t\,\sqrt{{2\lambda}/{3}}\,\Big)^{-3},\quad \lambda\geq0,$$
so that in particular $Z_t>0$, $P^\uparrow_0$ a.s.;
\item[(ii)] for every $t>0$, under $P^\uparrow_0$, conditionally on $(Z_u)_{0\leq u\leq t}$,
the process $(Z_{t+s})_{s\geq 0}$ is distributed according to $P^\uparrow_{Z_t}$.
\end{itemize}
From \eqref{Laplace-ex}, property (i) is equivalent to  saying that the density of $\z_t$ is $z\,h_t(z)$.
%We also note that the process $(Z_t)_{t\geq 0}$ is Markov under $P^\uparrow_x$, for any $x\geq 0$, and its
%transition kernels are given by
%$$E^\uparrow_x[\exp(-\lambda Z_{t+s})\,|\, Z_t] =  \Big(1+t\,\sqrt{{2\lambda}/{3}}\Big)^{-3}\, \exp\Big(-Z_t\,\Big(\lambda^{-1/2}+s\sqrt{2/3}\Big)^{-2}\Big),\quad \lambda >0,\;s,t\geq 0.$$

In the next proposition, we take $|\partial \mathfrak{B}^\bullet(0)|=0$ by convention.

\begin{proposition}
\label{law-horohull-bdry}
The process $(|\partial \mathfrak{B}^\bullet(r)|)_{r\geq 0}$ has a c\`adl\`ag modification, which is distributed 
according to $P^\uparrow_0$.
\end{proposition}

\proof
As a preliminary observation, we recall from \cite[Section 2.2]{CLG} that
 the exit measure process $(\z_{-r})_{r>0}$ is Markovian under 
$\N_0$, with the transition kernels of the $\Phi$-CSBP. In other words, we can find a 
c\`adl\`ag modification of $(\z_{-r})_{r>0}$ such that, for every $t>0$, 
the conditional distribution of $(\z_{-t-r})_{r\geq 0}$ under $\N_0$ and knowing $(\z_{-u})_{0<u\leq t}$
is $P_{\z_{-t}}$.  

In order to get the statement of the 
proposition, it suffices to verify that the finite-dimensional distributions of the process
$(|\partial \mathfrak{B}^\bullet(r)|)_{r\geq 0}$ coincide with the finite-dimensional marginals under $\P^\uparrow_0$.
So we need to verify that, for every $0<t_1<\cdots<t_p$, for every nonnegative measurable functions 
$\varphi_1,\ldots,\varphi_p$ on $\R_+$, we 
have
\begin{equation}
\label{hull-bdry1}
\E\Big[\varphi_1(|\partial \mathfrak{B}^\bullet(r_1)|)\cdots \varphi_1(|\partial \mathfrak{B}^\bullet(r_p)|)\Big]
= E_0^\uparrow[\varphi_1(Z_{r_1})\cdots \varphi_p(Z_{r_p})].
\end{equation}
Now recall from the preceding proof that, for every $j=1,\ldots,p$, 
$$|\partial \mathfrak{B}^\bullet(r_j)|= \int \mathbf{1}_{\{t<T_{r_j}\}}\,\z_{-r_j}(\omega)\,\pp(\dd t\,\dd\omega) + 
\int \mathbf{1}_{\{t<T_{r_j}\}}\,\z_{-r_j}(\omega)\,\pp'(\dd t\,\dd\omega).$$
It then follows from \eqref{tech-exit2} that
\begin{align*}
\E\Big[\varphi_1(|\partial \mathfrak{B}^\bullet(r_1)|)\cdots \varphi_1(|\partial \mathfrak{B}^\bullet(r_p)|)\Big]
&=\N_0\Big(\z_{-r_p}\,\varphi_1(\z_{-r_1})\cdots\varphi_p(\z_{-r_p})\Big)\\
&=\N_0\Big(\varphi_1(\z_{-r_1})\,E_{\z_{-r_1}}\Big[Z_{r_p-r_{1}}\,\varphi_2(Z_{r_2-r_1})\cdots \varphi_p(Z_{r_p-r_{1}})\Big]\Big),
\end{align*}
where we use the first observation of the proof in the last equality. Thanks to the $h$-transform relation \eqref{htrans}, the
right-hand side is also equal to 
$$
\N_0\Big(\varphi_1(\z_{-r_1})\z_{-r_1}E^\uparrow_{\z_{-r_1}}\Big[\varphi_2(Z_{r_2-r_1})\cdots \varphi_p(Z_{r_p-r_{1}})\Big]\Big)\!
=\!\!\int_0^\infty\! \dd z\,zh_{r_1}(z)\varphi_1(z)E^\uparrow_{z}\Big[\varphi_2(Z_{r_2-r_1})\cdots \varphi_p(Z_{r_p-r_{1}})\Big]
$$
since the density of $\z_{-r_1}$ under $\N_0(\cdot\cap\{\z_{-r_1}\not=0)$ is $h_{r_1}$ (Proposition \ref{densityuniform}). Finally,
properties (i) and (ii) above show that the right-hand side of the last display equals $E_0^\uparrow[\varphi_1(Z_{r_1})\cdots \varphi_p(Z_{r_p})]$, which completes the proof. \endproof

In the next proposition, we compute the joint distribution of the boundary size and the volume of the horohull $\mathfrak{B}^\bullet(r)$.

\begin{proposition}
\label{bdry-vol}
Let $r>0$. We have, for every $\lambda\geq 0$ and $\mu>0$,
$$\E\Big[\exp\Big(-\lambda |\partial\mathfrak{B}^\bullet(r)|-\mu \mathrm{Vol}(\mathfrak{B}^\bullet(r))\Big)\Big]
= \frac{(\frac{2}{3} + \frac{\lambda}{3} \sqrt{2/\mu}\,)^{-1/2} \sinh((2\mu) ^{1/4}r)+\cosh((2\mu)^{1/4}r)}
{\Big((\frac{2}{3} + \frac{\lambda}{3} \sqrt{2/\mu}\,)^{1/2} \sinh((2\mu) ^{1/4}r)+\cosh((2\mu)^{1/4}r)\Big)^3}.
$$
\end{proposition}

\proof From the fact that  $\mathfrak{B}^\bullet(r)=\Pi'^p(G_r\cup \partial G_r)$, with $G_r$ and $\partial G_r$ given by \eqref{horo-expl} and \eqref{horo-bdry}, we
easily obtain that
$$\mathrm{Vol}(\mathfrak{B}^\bullet(r)))=\int \mathbf{1}_{\{t<T_{r}\}}\,\sigma(\tr_{-r}(\omega))\,\pp(\dd t\,\dd\omega) + 
\int \mathbf{1}_{\{t<T_{r}\}}\,\sigma(\tr_{-r}(\omega))\,\pp'(\dd t\,\dd\omega).$$
Thanks to the similar formula for $|\partial\mathfrak{B}^\bullet(r)|$, and to Proposition \ref{Palm}, we get
$$\E\Big[\exp(-\lambda |\partial\mathfrak{B}^\bullet(r)|-\mu \mathrm{Vol}(\mathfrak{B}^\bullet(r)))\Big]
= \N_0\Big(\z_{-r}\,\exp(-\lambda\z_{-r}-\mu \mathcal{Y}_{-r})\Big),$$
where $\mathcal{Y}_{-r}=\int_0^\sigma \dd s\,\mathbf{1}_{\{\tau_{-r}(W_s)=\infty\}}$. By Lemma 4.5 in \cite{CLG}, we have, for $\lambda>\sqrt{2\mu}$,
$$\N_0\Big(1-\exp(-\lambda\z_{-r}-\mu \mathcal{Y}_{-r})\Big)
= \sqrt{\frac{\mu}{2}}\Bigg( 3\Bigg(\coth\Big((2\mu)^{1/4}r + \coth^{-1}\Big(\sqrt{\textstyle{\frac{2}{3}} +
\frac{\lambda}{3}\sqrt{2/\mu}}\,\Big)\Big)\Bigg)^2 - 2\Bigg).$$
By differentiating with respect to $\lambda$, we get the formula of the proposition. The restriction 
to $\lambda>\sqrt{2\mu}$ can be removed by an argument of analytic continuation. \endproof

\medskip
\rem Up to unimportant scaling constants, the formula of Proposition \ref{bdry-vol} already appears in \cite[Proposition 4]{CM},
which deals with asymptotics for the boundary size and volume of the (discrete) horohulls in the UIPT. This should not come
as a surprise since the Brownian plane is known to be the scaling limit of the UIPT \cite{Bud}. Note however that it
would not be easy to deduce Proposition \ref{bdry-vol} from the corresponding discrete result. 

\medskip
Our last proposition characterizes the distribution of the process $(|\partial\mathfrak{B}^\bullet(r)|, \mathrm{Vol}(\mathfrak{B}^\bullet(r)))_{r>0}$. This is an analog of \cite[Theorem 1.3]{CLG}, which is concerned with the usual hulls in the
Brownian plane.

\begin{proposition}
\label{process-bdry-vol}
Let $U=(U_t)_{t\geq 0}$ be a random process distributed according to $P^\uparrow_0$, and let
$s_1,s_2,\ldots$ be a measurable enumeration of jump times of $U$. Let $\xi_1,\xi_2,\ldots$
be an independent sequence of positive random variables distributed according to the
density $(2\pi x^5)^{-1/2}\exp(-1/2x)$. Assume that the sequence $(\xi_1,\xi_2,\ldots)$ is
independent of the process $U$. Then,
$$\Big(|\partial\mathfrak{B}^\bullet(r)|, \mathrm{Vol}(\mathfrak{B}^\bullet(r))\Big)_{r>0}
\build{=}_{}^{\rm(d)} \Bigg(U_r, \sum_{i:s_i\leq r} \xi_i\,(\Delta U_{s_i})^2\Bigg)_{r>0}.$$
\end{proposition}

From our presentation of the Brownian plane 
in terms of the triple $(B,\pp,\pp')$, and using Proposition \ref{Palm}
to relate this triple to the Brownian snake excursion measure, Proposition \ref{process-bdry-vol}
follows as a straightforward application of the excursion theory developed in \cite{ALG}
(see in particular Theorem 40 and Proposition 32 in \cite{ALG}). We omit the details of the proof. 

\subsection{Removing a strip from the Brownian half-plane}
\label{sec:strip}

In this section, we give an analog of Theorem \ref{infBdBp} showing that, if one removes a 
strip of width $r$ from the Brownian half-plane, the resulting space equipped with its
intrinsic metric is again a Brownian half-plane. We let $(\H_\infty,\Delta^{hp})$ stand for 
the Brownian half-plane constructed from a coding triple $(R,\wt\pp,\wt\pp')$ as explained in Section \ref{half-plane}.
Recall that $\H_\infty$ is obtained as a quotient space of the labeled tree $\t^{hp}_\infty$
associated with $(R,\wt\pp,\wt\pp')$, and 
that every $x\in \H_\infty$ thus has a label $\Lambda_x$, which is equal to the distance from
$x$ to the boundary $\partial\H_\infty$. 

We fix $r>0$ and set
$$\H^{(r)}_\infty:=\{x\in\H_\infty: \Lambda_x\geq r\}.$$
The interior $\H^{(r),\circ}_\infty$ is $\{x\in\H_\infty:\Lambda_x>r\}$. We distinguish a special point
$x_{(r)}$ of the boundary of $\H^{(r)}_\infty$, which corresponds to the point of the spine
of $\t^{hp}_\infty$ at height $\mathbf{L}_r=\sup\{t\geq 0: R_t= r\}$.

\begin{theorem}
\label{strip}
The intrinsic metric on $\H^{(r),\circ}_\infty$ has a unique continuous extension to $\H^{(r)}_\infty$, which is
a metric on this space. Furthermore, the resulting random measure metric space pointed at $x_{(r)}$ 
is a Brownian half-plane. 
\end{theorem}

\proof This is very similar to the proof of Theorem \ref{infBdBp}, and we only sketch the arguments.
We first introduce the process $(R^{(r)}_t)_{t\geq 0}$ defined by
$$R^{(r)}_t:=R_{\mathbf{L}_r+t}-r,$$
and we note that $(R^{(r)}_t)_{t\geq 0}$ is also a three-dimensional Bessel process started at $0$. 
Recalling the point measures $\wt\pp$ and $\wt\pp'$ used in the construction of $\H_\infty$, we 
define two other point measures $\wt\pp^{(r)}$ and $\wt\pp'^{(r)}$ on $\R_+\times\S$
by setting, for every nonnegative measurable function $\Phi$ on $\R_+\times \S$,
$$\langle \wt\pp^{(r)},\Phi\rangle=\int \pp(\dd t\dd \omega)\,\mathbf{1}_{(\mathbf{L}_r,\infty)}(t)\,\Phi(t-\mathbf{L}_r,
\vartheta_r(\tr_r\omega)),$$
and similarly for $\wt\pp'$, where we recall the notation $\vartheta_r$ for the shift on snake trajectories. Then it is straightforward
to verify that the coding triples $(R,\wt\pp,\wt\pp')$ and $(R^{(r)},\wt\pp^{(r)},\wt\pp'^{(r)})$ have the same
distribution.

Consequently, the construction of Section \ref{half-plane} applied to the triple $(R^{(r)},\wt\pp^{(r)},\wt\pp'^{(r)})$
yields a pointed measure metric space $(\H'^{(r)}_\infty,\Delta'^{hp,(r)})$ which is a Brownian half-plane.
To complete the proof we just have to identify $(\H'^{(r)}_\infty,\Delta'^{hp,(r)})$ with the space
$\H^{(r)}_\infty$ equipped with its intrinsic metric. This is done in the same way as in the
proof of Theorem \ref{infBdBp} and we omit the details. \hfill$\square$

\smallskip
\rem We could also have derived an analog of Theorem \ref{indep-hull} showing that
the space $\H^{(r)}_\infty$ in Theorem \ref{strip} is independent of the 
strip $\H_\infty\backslash \H^{(r)}_{\infty}$ equipped with its intrinsic metric. 
We leave the precise formulation and proof of this result to the reader.

\section{Consistency with previous definitions}
\label{sec:consist}

In this section, we show that our definitions of the infinite-volume Brownian disk and of the Brownian half-plane are
consistent with the previous definitions in \cite{BMR} and \cite{GM0}.
This is relatively easy for the Brownian half-plane but somewhat more delicate for the infinite-volume Brownian disk.

We start by recalling the definition of the free pointed Brownian disk that can be found in
\cite{BMR,Bet,BM}. Our presentation uses the notation introduced in the preceding sections
and is therefore slightly different from the one in the previous papers.

We fix $z>0$ and consider a Poisson point measure $\nn=\sum_{i\in I} \delta_{(t_i,\omega_i)}$ on
$\R_+\times \S$ with intensity
$$2\,\mathbf{1}_{[0,z]}(t)\,\dd t\,\N_0(\dd \omega).$$
We then introduce the compact metric space $\t'$, which is obtained from
the disjoint union
\begin{equation}
\label{tree-disk}
[0,z] \cup \Big(\bigcup_{i\in I} \t_{(\omega_i)}\Big)
\end{equation}
by identifying $0$ with $z$ and, for every $i\in I$, the root $\rho_{(\omega_i)}$ of $\t_{(\omega_i)}$
with the point $t_i$ of $[0,z]$. The metric on $\t'$ is defined in a very similar manner to 
Section \ref{coding-infinite}. For instance, if $v\in \t_{(\omega_i)}$ and $w\in \t_{(\omega_j)}$, with $j\not =i$, the
distance between $v$ and $w$ is
$$d_{(\omega_i)}(v,\rho_{(\omega_i)})+ \min\{(t_i\vee t_j)-(t_i\wedge t_j), z-(t_i\vee t_j)+(t_i\wedge t_j)\}+ d_{(\omega_j)}(\rho_{(\omega_j)},w),$$
and the reader will easily guess the formula in other cases. 
The volume measure on $\t'$ is just the sum of the volume measures on the trees $\t_{(\omega_i)}$, $i\in I$.

Set $\sigma':=\sum_{i\in I}\sigma(\omega_i)$. We can  define a clockwise exploration $(\ee'_t)_{0\leq t\leq\sigma'}$
of $\t'$, basically by concatenating the mappings $p_{(\omega_i)}:[0,\sigma(\omega_i)]\la \t_{(\omega_i)}$ in the 
order prescribed by the $t_i$'s. Note that, as in the finite spine case of Section \ref{coding-infinite}, this exploration
is cyclic (because $0$ and $z$ have been identified in $\t'$). The clockwise exploration allows us to define intervals 
in the space $\t'$, exactly as in Section \ref{coding-infinite}.

We next assign real labels to the points of $\t'$. To this end we let $(\beta_t)_{0\leq t\leq z}$ be a standard
Brownian bridge (starting and ending at $0$) over the time interval $[0,z]$, which is independent of $\nn$. For $t\in[0,z]$,
we set $\ell'_t=\sqrt{3}\,\beta_t$, and for $v\in\t_{(\omega_i)}$, $i\in I$, 
$$\ell'_v=\sqrt{3}\,\beta_{t_i} + \ell_v(\omega_i),$$
where $\ell_v(\omega_i)$ denotes the label of $v$ in $\t_{(\omega_i)}$, as in Section \ref{sna-tra}. 
Then, $\min\{\ell'_v:v\in \t'\}$ is attained at a unique point $v_*$ of $\t'$.

We may now define the pseudo-metric functions $D'^\circ$ and $D'$ exactly as 
in \eqref{pseudo-BP1} and \eqref{pseudo-BP2},
\begin{equation}
\label{pseudo-BD1}
D'^{\circ}(u,v)=\ell'_u + \ell'_v -2 \max\Big( \inf_{w\in[u,v]} \ell'_w,\inf_{w\in[v,u]} \ell'_w\Big),
\end{equation}
and
\begin{equation}
\label{pseudo-BD2}
D'(u,v) = \inf_{u_0=u,u_1,\ldots,u_p=v} \sum_{i=1}^p D'^{\circ}(u_{i-1},u_i)
\end{equation}
where the infimum is over all choices of the integer $p\geq 1$ and of the
finite sequence $u_0,u_1,\ldots,u_p$ in $\t'$ such that $u_0=u$ and
$u_p=v$. It is immediate to verify that, for every $u\in\t'$, $D'^\circ(u,v_*)=D'(u,v_*)= \ell'_u-\ell'_{v_*}$.

Let $\D'^\bullet_z$ denote the space $\t'/\{D'=0\}$, which is equipped with the metric induced by
$D'$, with the pushforward of the volume measure on $\t'$, and with the distinguished point which
is the equivalence class of $v_*$ (without risk of confusion, we will also write $v_*$ for this
equivalence class). Then $\D'^\bullet_z$ is a free pointed Brownian disk with perimeter $z$ whose boundary 
$\partial\D'^\bullet_z$ is the image of $[0,z]$ under the canonical projection. This
construction is basically the one in \cite[Section 2.3]{BM}, and, as already noted, it is consistent with Definition \ref{def-disk}. 
We set 
\begin{equation}
\label{dist-bdry-v*}
H'_z=D'(v_*,\partial \D'^\bullet_z)=\min\{\ell'_v:v\in[0,z]\} - \ell'_{v_*}.
\end{equation}

A variant of the preceding construction yields the infinite-volume Brownian disk with perimeter $z$
as considered\footnote{Unfortunately, it seems that the definition given in \cite{BMR} is slightly incorrect. We believe that
the construction below is the correct way to define the infinite-volume Brownian disk as it appears in the
limit theorems proved in \cite{BMR}.} in \cite{BMR}.
We keep the same notation as before, and we introduce an infinite labeled tree $\t'_\infty$ which has the same 
distribution as the tree $\t'^p_\infty$ of Section \ref{sec:horo} (so this is the labeled tree associated with a triple
$(B,\pp,\pp')$ whose distribution is specified in Section \ref{sec:horo}). We assume that $\nn$ and $\t'_\infty$ are independent, and we also consider a
random variable $U$ uniformly distributed over $[0,z]$ and independent of the pair $(\nn,\t'_\infty)$. Then 
we let $\t'^{(\infty)}$ be derived from the disjoint union
\begin{equation}
\label{tree-infdisk}
[0,z] \cup \Big(\bigcup_{i\in I} \t_{(\omega_i)}\Big)\cup \t'_\infty
\end{equation}
by the same identifications as in \eqref{tree-disk}, and furthermore by identifying the root
of $\t'_\infty$ with the point $U$ of $[0,z]$. The metric on $\t'^{(\infty)}$ is defined as in
the case of $\t'$. The clockwise exploration $(\ee'^{(\infty)}_t)_{t\in\R}$ of $\t'^{(\infty)}$ is then defined in much the same
way as in the infinite spine case of Section \ref{coding-infinite}: We have
$\ee'^{(\infty)}_0=0=z$, and the points $(\ee'^{(\infty)}_t)_{t<0}$ now correspond to
the right side of the tree $\t'_\infty$ , to the trees $\t_{(\omega_i)}$ with $t_i>U$ and to the interval $[U,z)$, and 
similarly for the points $(\ee'^{(\infty)}_t)_{t>0}$. The labels $\ell'^{(\infty)}_v$ on $\t'^{(\infty)}$ are obtained exactly as in the case of $\t'$,
using the same Brownian bridge $\beta$ and taking $\ell'^{(\infty)}_v=\beta_U+ \Lambda'_v$ when $v\in \t'_\infty$, where 
$\Lambda'_v$ stands for the label of $v$ in $\t'_\infty$.

We may now define the pseudo-metric functions $D'^{\circ(\infty)}(u,v)$ and $D'^{(\infty)}(u,v)$
on $\t'^{(\infty)}$ by the very same formulas as in \eqref{pseudo-BD1} and \eqref{pseudo-BD2},
just replacing the labels $\ell'_v$ by $\ell'^{(\infty)}_v$ (and noting that the clockwise exploration
$(\ee'^{(\infty)}_t)_{t\in\R}$ allows us to define intervals on $\t'^{(\infty)}$, exactly as in
Section \ref{coding-infinite}). 

We then define $\D'^\infty_z$ as the quotient space $\t'^{(\infty)}/\{D'^{(\infty)}=0\}$, which is equipped 
with the metric induced by $D'^{(\infty)}$, with the volume measure which is the pushforward of the volume
measure on $\t'^{(\infty)}$ and with the distinguished point which is the equivalence class of $0$. 
In the terminology of \cite{BMR}, $\D'^\infty_z$ is an infinite-volume Brownian disk with perimeter $z$. The next
proposition shows that this is
consistent with Definition \ref{def-infinite-disk}. 

\begin{proposition}
\label{ident-infinite-BD}
The pointed locally compact measure metric spaces $\D^\infty_z$ and $\D'^\infty_z$
have the same distribution.
\end{proposition}

We will deduce Proposition \ref{ident-infinite-BD} from Proposition \ref{local-G-H-bis} below, which
shows that  
$\D'^\infty_z$ is a limit of conditioned Brownian disks, in a way similar to Theorem \ref{local-G-H} for $\D^\infty_z$. We note that \cite{BMR} proves that 
the space $\D'^\infty$ is the limit in distribution of Brownian disks with perimeter $z$ conditioned to
have a large volume, but it is not so easy to verify that this conditioning has the same effect
as the one in Theorem \ref{local-G-H}, which involves the height of the distinguished point. 

Let us start with some preliminary observations. Since $\D'^\bullet_z$
has the same distribution as $\D^\bullet_z$, we know from the discussion after
Proposition \ref{pointed-disk} that there exists 
a  measure $\mu'_z$ on $\partial \D'^\bullet_z$ with total mass $z$, such that, a.s. for any continuous function
$\varphi$ on $\D'^\bullet_z$, we have
$$\langle \mu'_z,\varphi\rangle =\lim_{\ve\to 0} \frac{1}{\ve^2}\int_{\D'^\bullet_z} \mathrm{Vol}(\dd x)\,\mathbf{1}_{\{D(x,
\partial \D'^\bullet_z)<\ve\}}\,\varphi(x)$$
where $\mathrm{Vol}(\dd x)$ denotes the volume measure on $\partial \D'^\bullet_z$. For our purposes, it is important to know that $\mu'_z$ is also
the pushforward of Lebesgue measure on $[0,z]$ under the canonical projection 
from $\t'$ onto $\D'^\bullet_z$. This is proved in 
\cite[Theorem 9]{Disk-Bdry}. 

We now note that, in addition to $v_*$, $\D'^\bullet_z$ has another distinguished point (belonging to its boundary)
namely the point $v_\partial$ which is the equivalence class of $0$ in the quotient $\t'/\{D'=0\}$. 
We note that $v_\partial$ is uniformly distributed over $\partial\D'^\bullet_z$, in the following sense.
Similarly as in Proposition \ref{doubly-point}, we introduce the 
the doubly pointed 
measure metric space $\D'^{\bullet\bullet}_z$ which is obtained by viewing $v_\partial$ as a second
distinguished point of $\D'^{\bullet}_z$. We have then, for any nonnegative measurable function $F$ on the space of all doubly pointed compact
measure metric spaces,
\begin{equation}
\label{bdry-uni}
\E[F(\D'^{\bullet\bullet}_z)]=\frac{1}{z}\,\E\Big[ \int \mu'_z(\dd x)\,F\Big([\D'^{\bullet}_z,x]\Big)\Big],
\end{equation}
with the same notation as in Proposition \ref{doubly-point}. The proof of \eqref{bdry-uni} is straightforward: 
For $r,t\in[0,z]$ use the notation $t\oplus r=t+r$ if $t+r\leq z$ and $t\oplus r=t+r-z$ if $t+r>z$, and 
note that, for every $r\in[0,z]$, the point measure $\sum_{i\in I} \delta_{(t_i\oplus r,\omega_i)}$ has the same 
distribution as $\nn$, whereas $(\beta_{r\oplus t}-\beta_r)_{0\leq t\leq z}$ has the same distribution as
$(\beta_t)_{0\leq t\leq z}$. 

Consider a random doubly pointed space $\D^{\bullet\bullet}_z$ whose distribution is obtained by
integrating the distribution of $\D^{\bullet\bullet,a}_z$ with respect to the probability density $p_z(a)$ in Proposition \ref{densityuniform}. By integrating
the formula of Proposition \ref{doubly-point} with respect to $p_z(a)\,\dd a$, we get
\begin{equation}
\label{bdry-uni-bis}
\E[F(\D^{\bullet\bullet}_z)]=\frac{1}{z}\,\E\Big[ \int \mu_z(\dd x)\,F\Big([\D^{\bullet}_z,x]\Big)\Big].
\end{equation}
Since the pairs $(\D'^{\bullet}_z,\mu'_z)$ and 
$(\D^{\bullet}_z,\mu_z)$ have the same distribution, we obtain by comparing \eqref{bdry-uni}
and \eqref{bdry-uni-bis} that $\D'^{\bullet\bullet}_z$ and $\D^{\bullet\bullet}_z$ have the
same distribution. Let $\bar\D^\bullet_z$, resp. $\bar\D'^\bullet_z$, be the pointed space obtained from 
$\D^{\bullet\bullet}_z$, resp. from $\D'^{\bullet\bullet}_z$, by forgetting the first distinguished point.
Then $(\bar\D^\bullet_z,H_z)$ and $(\bar\D'^\bullet_z,H'_z)$ also have the same distribution.

\begin{proposition}
\label{local-G-H-bis}
For every $a>0$, let $\bar\D'^{\bullet,(a)}_z$ be distributed as $\bar\D'^\bullet_z$ conditioned 
on the event $\{H'_z\geq a\}$. Then
$$\bar\D'^{\bullet,(a)}_z \build{\la}_{a\to\infty}^{\rm(d)} \D'^\infty_z$$
in distribution in the sense of the local Gromov-Hausdorff-Prokhorov convergence. 
\end{proposition}

Before proving Proposition \ref{local-G-H-bis}, let us explain why the statement of Proposition \ref{ident-infinite-BD}
follows from this proposition. Recall from Section \ref{sec:point-disk} that $H_z$ denotes the distance from the distinguished 
point of $\D^\bullet_z$ to the boundary. 
Since $(\bar\D^\bullet_z,H_z)$ and $(\bar\D'^\bullet_z,H'_z)$ have the same distribution, $\bar\D'^{\bullet,(a)}_z$
has the same distribution as $\bar\D^\bullet_z$ conditioned on $H_z\geq a$,
whereas the pointed space $\bar\D^{\bullet,a}_z$ in Theorem \ref{local-G-H}
has the distribution of $\bar\D^\bullet_z$ conditioned on $H_z= a$.
Hence, by comparing the convergences in Theorem \ref{local-G-H} and
in Proposition \ref{local-G-H-bis}, we conclude that $\D^\infty_z$
and $\D'^\infty_z$ have the same distribution. 

\medskip
\noindent{\it Proof of Proposition \ref{local-G-H-bis}.} Let $E_a$ stand for the event $\{H'_z\geq a\}$. The idea of the proof is 
to study the effect on the pair $(\beta,\nn)$ (which determines $\D'^{\bullet\bullet}_z$) of 
conditioning on $E_a$. To this end, it will be useful to replace $E_a$ by another 
event for which the conditioning will be easier to study. We first note that, by \eqref{dist-bdry-v*}
and the definition of labels on $\t'$, we have
\begin{equation}
\label{altern-dist}
H'_z= \min_{0\leq t\leq z}(\sqrt{3}\,\beta_t) - \inf_{i\in I} \Big(\sqrt{3}\,\beta_{t_i} + W_*(\omega_i)\Big).
\end{equation}
Set $\|\beta\|=\sup\{|\beta_t|:0\leq t\leq z\}$, and consider the events
$$\wt E_a:=\Big\{\inf_{i\in I} W_*(\omega_i)\leq -a\Big\},\,E'_a:=\Big\{\inf_{i\in I} W_*(\omega_i)\leq -a -2\sqrt{3}\|\beta\|\Big\},\,
E''_a:=\Big\{\inf_{i\in I} W_*(\omega_i)\leq -a +2\sqrt{3}\|\beta\|\Big\}.$$
From \eqref{altern-dist}, we have $E'_a\subset E_a\subset E''_a$. On the other hand, it is an easy exercise 
to check that the ratio $\P(E'_a)/\P(E''_a)$ tends to $1$ as $a\to\infty$
(in fact, it follows from \eqref{hittingpro} that both $\P(E'_a)$ and $\P(E''_a)$ are 
asymptotic to $3z/a^2$). Since we have also $E'_a\subset \wt E_a\subset E''_a$, we may
condition on $\wt E_a$ instead of conditioning on $E_a=\{H'_z\geq a\}$ in order to get
the convergence of the proposition.

Conditioning on $\wt E_a$ does not affect $\beta$. On the other hand, when $a$ is large, the conditional distribution 
of $\nn$ knowing $\wt E_a$ is close in total variation to the law of
$$\nn'^{(a)} + \delta_{(\wt U,\omega^{(a)})},$$
where $\nn'^{(a)}$ is a Poisson point measure with intensity 
$2\,\mathbf{1}_{[0,z]}(t)\,\mathbf{1}_{\{W_*(\omega)>-a\}}\,\dd t\,\N_0(\dd \omega)$,
$\omega^{(a)}$ is distributed according to $\N_0(\cdot\mid W_*\leq -a)$, and $\wt U$
is uniformly distributed over $[0,z]$ (and $\nn'$, $\omega^{(a)}$ and $\wt U$
are independent). When $a$ is large, $\nn$ and $\nn'^{(a)}$ can be coupled 
so that they are equal with high probability. 

We then want to argue that, when $a$ is large, we can couple $\omega^{(a)}$
and the labeled tree $\t'_\infty$ used to define $\D'_\infty$ so that $\t_{(\omega^{(a)})}$ and $\t'_\infty$,
both viewed as labeled trees, are close
in some appropriate sense. Recall that $\t'_\infty$ was constructed from 
a coding triple $(B,\pp,\pp')$ such that
$B=(B_t)_{t\geq 0}$ is a linear Brownian motion started from $0$ and,
conditionally on $B$, $\pp$ and $\pp'$ are independent Poisson point measures on $\R_+\times \S$
with intensity
$2\,\dd t\,\N_{B_t}(\dd \omega)$. On the other hand, the main results of \cite{Bessel} give the
distribution of $\omega^{(a)}$. If $b\in[a,\infty)$, the conditional distribution of
$\omega^{(a)}$ knowing that $W_*(\omega^{(a)})=-b$ is that of the snake trajectory corresponding to
a coding triple $(V,\mm,\mm')$ such that $V=(V'_t-b)_{0\leq t\leq T^{V'}}$, where $(V'_t)_{0\leq t\leq T^{V'}}$
is a Bessel process of dimension $-5$ started from $b$ and stopped when it hits $0$, and,
conditionally on $V$, $\mm$ and $\mm'$ are independent Poisson measures 
on $\R_+\times \S$ with intensity 
$$2\,\mathbf{1}_{[0,T^{V'}]}(t)\,\mathbf{1}_{\{W_*(\omega)>-b\}}\,\dd t\,
\N_{V_t}(\dd \omega).$$
From this description, we easily get that, for every $h>0$ and $\ve\in(0,1)$, we can for $a$ large enough 
couple the coding triples $(B,\pp,\pp')$ and $(V,\mm,\mm')$ in such a way that the following two properties hold
except on a set of probability smaller than $\ve$:
\begin{itemize}
\item $V_t=B_t$ for $0\leq t\leq h$;
\item the restriction of $\pp$, resp.~of $\pp'$, to $[0,h]\times \S$ coincides with 
the restriction of $\mm$, resp.~of $\mm'$, to $[0,h]\times \S$. 
\end{itemize}
Now recall that the construction of $\D'_\infty$ relies on the 4-tuple $(\beta,\nn,U,\t'_\infty)$, whereas,
up to an event of small probability when $a$ is large, the space $\bar\D'^{\bullet,(a)}_z$ (which is $\bar\D'^\bullet_z$ conditioned on 
$\{H_z\geq a\}$) may be obtained from the 4-tuple $(\beta,\nn'^{(a)},\wt U, \t_{(\omega_{(a)})})$.
It follows from the preceding considerations that, up to a set of small probability when $a$ is large,
we can couple these two 4-tuples in such a way that their first three components coincide and moreover
the labeled trees $\t'_\infty$ and $ \t_{(\omega_{(a)})}$ with their spines ``truncated at height $h$'' also
coincide (in the case of $ \t_{(\omega_{(a)})}$, the spine corresponds to the line segment between the 
root and the vertex with minimal label). Given $r>0$, we deduce from the preceding observation
(by choosing $h$ large enough) that
we can couple the spaces $\D'^\infty_z$ and $\bar\D'^{\bullet,(a)}_z$ so that the balls of radius $r$ centered at the distinguished point are the same
in both spaces, except on an event of small probability when $a$ is large. We omit the detailed 
verification of this last coupling, which is very similar to the proof of Proposition \ref{cons-infi} above
or Theorem 1 in \cite{Plane}. The convergence in distribution stated in 
Proposition \ref{local-G-H-bis} follows. \hfill$\square$

\medskip
\rem The quantities $H_z$ and $H'_z$ have the same distribution, and thus the density of
the random variable in the right-hand side of \eqref{altern-dist} is equal to $p_z(a)$.
The reader is invited to give a direct proof of this fact, as a verification of the consistency
of our definition of the free pointed Brownian disk with the one in \cite{BM}.

\medskip

To conclude this section, we explain why our definition of the Brownian half-plane $\H_\infty$ is
consistent with the one given in \cite{BMR} or \cite{GM0}. We use the notation $\H'_\infty$ for the
Brownian half-plane as defined in \cite{BMR}. Then, from \cite[Corollary 3.9]{BMR}, we get that 
$\H'_\infty$ is the tangent cone in distribution of the Brownian disk with perimeter $z$ and volume $r$
at a point uniformly distributed over its boundary --- here ``uniformly distributed'' refers to the 
analog of the measure $\mu'_z$  (the construction of $\bar\D'^\bullet_z$
given above also works for the Brownian disk with fixed volume $r$, just by conditioning 
$\sum_{i\in I}\sigma(\omega_i)$ to be equal to $r$). By randomizing the volume $r$, we infer that we have also
$$\lambda \cdot \bar \D'^{\bullet}_z \build{\la}_{\lambda\to\infty}^{\rm(d)} \H'_\infty,$$
in distribution in the sense of the local Gromov-Hausdorff-Prokhorov convergence. On the other hand, 
it follows from Theorem \ref{conv-hp} that we have 
$$\lambda \cdot \bar \D^{\bullet}_z \build{\la}_{\lambda\to\infty}^{\rm(d)} \H_\infty.$$
Since $\D^{\bullet}_z$ and $\D'^{\bullet}_z$ have the same distribution, we conclude 
that $\H_\infty$ and $\H'_\infty$ also have the same distribution as desired.

\section*{Appendix: Some Laplace transforms}

Recall the standard notation
$$\hbox{erfc}(x)=\frac{2}{\sqrt{\pi}}\int_x^\infty e^{-t^2}\,\dd t.$$
Then the function $\chi_1$ defined for $x>0$ by
$$\chi_1(x)=\frac{1}{\sqrt{\pi}}x^{-1/2} - e^x\,\mathrm{erfc}(\sqrt{x})
=\frac{1}{\sqrt{\pi}}\,e^x \int_{\sqrt{x}}^\infty \frac{1}{t^2}\,e^{-t^2}\,\dd t,\eqno{\rm(A.0)}$$
satisfies, for every $\lambda>0$,
$$\int_0^\infty \dd x\,e^{-\lambda x}\,\chi_1(x)= (1+\sqrt{\lambda})^{-1}.\eqno{\rm(A.1)}$$
This is easily verified via an integration by parts which gives for $\lambda>0$,
$$\int_0^\infty  \mathrm{erfc}(\sqrt{x})e^x\,e^{-\lambda x}\,\dd x = \frac{1}{\sqrt{\lambda}(1+\sqrt{\lambda})}.$$
From the last two displays and an integration by parts, one checks that
the function $\chi_2=\chi_1*\chi_1$, which
satisfies
$$\int_0^\infty \dd x\,e^{-\lambda x}\,\chi_2(x)= (1+\sqrt{\lambda})^{-2},\eqno{\rm(A.2)}$$
is given for $x>0$ by
$$\chi_2(x)= e^x\mathrm{erfc}(\sqrt{x})-2x\,\chi_1(x)= (2x+1)e^x\mathrm{erfc}(\sqrt{x}) -\frac{2}{\sqrt{\pi}}x^{1/2}.$$ 
Similar manipulations show that the function $\chi_3=\chi_1*\chi_1*\chi_1$ satisfying
$$\int_0^\infty \dd x\,e^{-\lambda x}\,\chi_3(x)= (1+\sqrt{\lambda})^{-3}.\eqno{\rm(A.3)}$$
is given by
$$\chi_3(x)=\frac{2}{\sqrt{\pi}}(x^{3/2} + x^{1/2}) - 2x(x+\frac{3}{2})\,e^x\,\hbox{erfc}(\sqrt{x}).$$
We observe that $\chi_1(x)>0$ for every $x>0$ (this is obvious from (A.0)) and thus we have also $\chi_3(x)>0$
for every $x>0$.
Finally, we note that
$$\int_0^{\infty} \frac{1-e^{-\lambda x}}{x}\,\chi_3(x)\,\dd x
= \int_0^\lambda \dd \mu\int_0^\infty e^{-\mu x}\,\chi_3(x)\,\dd x = \int_0^\lambda \dd \mu(1+\sqrt{\mu})^{-3} 
= (1+\lambda^{-1/2})^{-2}. \eqno{\rm(A.4)}$$

\end{document}